\documentclass[12pt]{amsart}
\usepackage{amsmath,amssymb,amscd}
\usepackage{stmaryrd}
\usepackage{breqn}

\numberwithin{equation}{section}

\newtheorem{lemma}{Lemma}[section]
\newtheorem{theorem}{Theorem}[section]
\newtheorem{corollary}{Corollary}[section]
\newtheorem{proposition}{Proposition}[section]

\newtheorem{prop}{Proposition}[section]

\usepackage{tikz}
\usetikzlibrary{arrows}
\usetikzlibrary{decorations.markings}

\usetikzlibrary{decorations.pathreplacing,decorations.markings}
\newcommand{\bl}{\begin{lemma}}
\newcommand{\el}{\end{lemma}}
\newcommand{\bp}{\begin{proposition}}
\newcommand{\ep}{\end{proposition}}
\newcommand{\bcor}{\begin{corollary}}
\newcommand{\ecor}{\end{corollary}}
\newcommand{\bth}{\begin{theorem}}
\newcommand{\et}{\end{theorem}}
\newcommand{\be}{\begin{equation}}
\newcommand{\ee}{\end{equation}}
\newcommand{\bal}{\begin{align}}
\newcommand{\eal}{\end{align}}
\newcommand{\bi}{\begin{itemize}}
\newcommand{\ei}{\end{itemize}}
\newcommand{\la}{\label}
\newcommand{\ed}{\,{\buildrel d \over =}\,}

\newcommand{\bs}{{\bigskip}}
\newcommand{\ms}{{\medskip}}
\newcommand{\noi}{\noindent}

\newcommand{\bq}{{\bold q}}

\newcommand{\bn}{{\bold n}}

\newcommand{\bX}{{\bold X}}

\newcommand{\brho}{{\pmb{\rho}}}

\renewcommand{\a}{\alpha}
\renewcommand{\b}{\beta}
\renewcommand{\d}{\delta}
\newcommand{\D}{\Delta}

\newcommand{\g}{\gamma}
\newcommand{\G}{\Gamma}
\renewcommand{\l}{\lambda}
\renewcommand{\L}{\Lambda}

\newcommand{\var}{\varphi}
\newcommand{\s}{\sigma}

\renewcommand{\th}{\theta}
\renewcommand{\O}{\Omega}
\renewcommand{\o}{\omega}

\newcommand{\x}{\times}
   
\renewcommand{\i}{\infty}
\newcommand{\p}{\partial}

\newcommand{\bE}{{\mathbb E}}
\newcommand{\bN}{{\mathbb N}}
\newcommand{\bP}{{\mathbb P}}
\newcommand{\bR}{{\mathbb R}}

\newcommand{\bZ}{{\mathbb Z}}

\newcommand{\bz}{{\bf{z}}}

\newcommand{\cK}{{\mathcal K}}

\newcommand{\cV}{{\mathcal V}}

\newcommand{\cL}{{\mathcal L}}

\newcommand{{\cA}}{{\mathcal A}}
\newcommand{\cM}{{\mathcal M}}

\newcommand{\cQ}{{\mathcal Q}}
\newcommand{\cF}{{\mathcal F}}
\newcommand{\cH}{{\mathcal H}}

\newcommand{\cX}{{\mathcal X}}

\newcommand{\cC}{{\mathcal C}}
\newcommand{\cE}{{\mathcal E}}
\newcommand{\cS}{{\mathcal S}}
\newcommand{\cT}{{\mathcal T}}

\def\lf{\left\lfloor}   
\def\rf{\right\rfloor}
\newcommand{\1}{1\!\!1}
\newcommand{\mf}{\mathfrak}

\DeclareSymbolFont{bbold}{U}{bbold}{m}{n}
\DeclareSymbolFontAlphabet{\mathbbold}{bbold}
\newcommand{\ind}{\mathbbold{1}}

\begin{document}
\author{Mehdi Ouaki}
\address{Department of Statistics, University of California, Berkeley, CA 94720-3840}
\email{mouaki@berkeley.edu}

\author{Fraydoun Rezakhanlou}
\address{Department of Mathematics, University of California, Berkeley, CA 94720-3840}
\email{rezakhan@math.berkeley.edu}

\title{ Random Tessellations and Gibbsian Solutions of Hamilton-Jacobi Equations}

\maketitle

\begin{abstract}
We pursue two goals in this article. As our first goal, we construct a family $\cM_G$ of 
Gibbsian measures on the set of
piecewise linear convex functions $g:\bR^2\to\bR$. It turns out that there is a one-to-one correspondence between the gradient of such convex functions and {\em Laguerre tessellations}. Each cell in a Laguerre tessellation is a convex polygon that is marked by a vector $\rho\in\bR^2$. Each measure $\nu^f\in\cM_G$ in our family is uniquely characterized by a kernel $f(x,\rho^-,\rho^+)$, which represents the rate at which a line separating two cells associated with marks $\rho^-$ and $\rho^+$ passes through $x$. To construct our measures, we give a precise recipe for the law of the restriction of our tessellation to a box. This recipe involves a boundary condition, and a dynamical description of our random tessellation inside the box. As we enlarge the box, the consistency of these random tessellations requires that the kernel satisfies a suitable kinetic like PDE. 
Our interest in the family $\cM_G$ stems from its close connection to Hamilton-Jacobi PDEs of the form $u_t=H(u_x)$, 
with $u:\G\to\bR$, for a convex set 
$\G\subset\bR^2\x[0,\i)$, and $H:\bR^2\to\bR$ a strictly convex function. As our second goal, we study the invariance of the set $\cM_G$ 
with respect to the dynamics of such Hamilton-Jacobi PDEs. 
In particular we {\em conjecture} the invariance of 
a suitable subfamily $\widehat \cM_G$ of $\cM_G$.
More precisely,
 we expect that if the initial slope $u_x(\cdot,0)$ is selected according to a measure $\nu^{f}\in\widehat \cM_G$, then at a later time
 the law of $u_x(\cdot, t)$
   is  given by a measure  
$\nu^{\Theta_t(f)}\in\widehat\cM_G$,
for a suitable kernel $\Theta_t(f)$. As we vary $t$, the kernel 
$\Theta_t(f)$ must satisfy a suitable kinetic equation. 
We remark that the function $u$ is also piecewise linear convex function in $(x,t)$, and its law is an example of a Gibbsian measure
on the set of Laguerre tessellations of certain convex subsets of $\bR^3$.
\end{abstract}

\section{Introduction}
\label{sec1}

In numerous models of nonequilibrium statistical mechanics we  
encounter an interface that separates different phases and is evolving with time.  It is often the case that the evolution of such an interface depends on the location $x$,  the time $t$, and the inclination $\rho$ of the interface at $x$. 
If the interface is represented by a graph of a {\em height} function 
$u:\mathbb{R}^d\x[0,\i)\to\mathbb{R}$,  then a natural model for its evolution is a Hamilton-Jacobi PDE:
\be\la{eq1.1}
u_t=H(x,t,u_x).
\ee
Since in practice the exact form of the Hamiltonian 
function $H$ is not known to us,
 it is a common practice to assume that $H$ is random. A natural question is whether or not we can describe the stochastic law $\nu_t$ of the height function $u(\cdot,t)$ as $t$ varies. Ideally we would like to derive a
tractable/explicit evolution equation for $\nu_t$.  Alternatively, we may keep track of the inclination $\rho=u_x$, and wonder whether or not the law of $\rho(\cdot,t)$ follows an explicitly describable evolution equation. This is indeed the case for a small number of
exactly solvable one dimensional discrete models. In this article however, we pursue a very different strategy:
we search for a natural class of stochastic laws that is invariant with respect to the evolution of the 
Hamilton-Jacobi PDE \eqref{eq1.1}. This strategy has already been tested in dimension one: If initially 
the process $x\mapsto \rho(x,0)$ evolves as an ODE that is interrupted by Markovian jumps, then the same is true at later times, and the evolution of the jump rates can be described by a kinetic equation. 
In the present article we examine this strategy in higher dimensions. It turns out that the evolution of the height function is significantly more complex when $d>1$.  Fortunately, when 
 $H$ is independent of $(x,t)$, and convex in the {\it momentum} variable $\rho$, the dynamics simplify and 
the classical formulas of Hopf, Lax and Oleinik lead to variational representations of solutions. A particularly tractable case is when the height function is piecewise linear and convex. As our first step, we offer a recipe
for a Gibbsian measure on the set of piecewise linear and convex functions. For our second step, we study the evolution of this measure with respect to the Hamilton-Jacobi dynamics.

\subsection{Hamilton-Jacobi semigroup}
As a preparation for the statement of our main results, we first recall two classical variational formulas for solutions of \eqref{eq1.1} when $H$ is convex and independent of $(x,t)$:

\ms\noi
{\bf(1)} {\em (Hopf Formula)} If $g$ is convex, then 
\be\la{eq1.2}
u(x,t)=(g^*-tH)^*(x),
\ee
where $g(x)=u(x,0)$ is the initial condition, and $g^*$ denotes the Legendre transform of $g$.
More explicitly,
\begin{align}
u(x,t)&=\sup_{\rho}\big(x\cdot \rho-g^*(\rho)+tH(\rho)\big).
\la{eq1.3}
\end{align}

\ms\noi
{\bf(2)} {\em (Hopf-Lax-Oleinik Formula)} If $H$ is convex, then 
\be\la{eq1.4}
u(x,t)=\sup_{y}\left(g(y)-tL\left(\frac{y-x}t\right)\right),
\ee
where $L=H^*$ is the Legendre transform of $H$ (see for example [E]).

\ms
In this article we assume that both the Hamiltonian function and the initial data are convex.
An immediate consequence of the convexity of $H$ is that the flow of \eqref{eq1.1} is {\em strongly monotone}. More precisely, if we  write 
$\Phi_t$ for the flow of our PDE:
\[
\Phi_t(g)(x)=u(x,t)\ \ \ \Leftrightarrow\ \ \ 
u_t=H(u_x),\ \ {\text{and}}\ \ u(x,0)=g(x),
\]
Then $\Phi_t\big(\sup_\rho h^\rho\big)=\sup_\rho\Phi_t\big( h^\rho\big)$, which would follow from \eqref{eq1.4}. In particular, if we choose $h^\rho$ to be a linear function of the form
$h^\rho(x)=x\cdot\rho-g^*(\rho)$, then 
\[
\Phi_t(h^\rho)(x)=x\cdot \rho-g^*(\rho)+tH(\rho),
\]
which in turn implies \eqref{eq1.2}-\eqref{eq1.3}.
For our purposes, it is more convenient to keep track of 
the slope $\rho(x,t)=u_x(x,t)$. Its evolution with respect to time can be represented by a semigroup $\widehat\Phi_t$;
\[
\widehat\Phi_t(\nabla g)(x)=\rho(x,t).
\]

\bs\noi
{\bf Definition 1.1(i)}
Given a convex set $\L$, we write $\cC(\L)$ for the set
of convex functions $g:\L\to\bR$. The set of piecewise 
linear functions $g\in\cC(\L)$ is denoted by
$\cC_0(\L)$.

\ms\noi
{\bf(ii)}
We write 
$\widehat \cC_0(\L)$ for the set of functions $\rho:\L\to\bR^d$ such that $\rho=\nabla g$, for some $g\in\cC_0(\L)$.
\qed

\bs
Observe that $\Phi_t\big(\cC(\bR^d)\big)\subset\cC(\bR^d)$ by \eqref{eq1.2}. Moreover, the set of piecewise linear 
convex functions $\cC_0(\bR^d)$ is  
also invariant with respect to $\Phi_t$. 
 As one of our main contribution, we construct 
a family $\cM_G$ of Gibbsian measures on $\widehat \cC_0$ that is expected to be invariant 
under the flow $\widehat\Phi_t$.

\subsection{Tessellations}

It turns out that there is a one-to-one correspondence between the members of $\widehat\cC_0(\bR^d)$, and the {\em Laguerre tessellations} of $\bR^d$. 
Henceforth our aforementioned set $\cM_G$ offers a 
natural family of Gibbsian measures on the set of tessellations. To explain this further, let us remark that $g\in\cC_0$ means that there exists 
a {\em discrete set} $ \mathcal{S} \subset\bR^d$ such that
\[
g(x)=\sup_{\rho\in \mathcal{S}}\big(x\cdot\rho-g^*(\rho)\big).
\]
If we set 
\[
X(\rho)=\big\{x\in\bR^d:\ g(x)=x\cdot\rho-g^*(\rho)\big\},
\]
then the {\em cell} 
$X(\rho)$ is a {\em convex polytope} and the collection
\[
\big\{\big(\rho,X(\rho)\big):\ \rho\in \mathcal{S} \big\},
\]
is a {\em Laguerre tessellation} of $\bR^d$. Moreover, the function $\nabla g$ is piecewise constant, and has a representation of the form
\be\la{eq1.5}
\rho(x):=\nabla g(x)=\sum_{\rho\in \mathcal{S}}\rho\ \1\big(x\in X(\rho)\big).
\ee
It is this geometric interpretation that is at the heart of our strategy for constructing our Gibbsian measures.

 \ms
Before embarking on our construction, we first need to come up with criteria that would guarantee that a polytope tessellation
does come from a function $g\in \cC_0$. 
Since we will be mostly studying planar tessellations in this article, let us assume that $d=2$. 
Indeed if 
\[
X(\rho^-,\rho^+):=X(\rho^-)\cap X(\rho^+)\neq \emptyset,
\]
for a planar tessellation,
then generically the set $X(\rho^-,\rho^+)$ is a line segment,
 and if $\tau(\rho^-,\rho^+)$
is a vector that is parallel to this line segment, then we must have
\be\la{eq1.6}
\tau(\rho^-,\rho^+)\cdot (\rho^+-\rho^-)=0.
\ee
We can readily verify this using the fact that the linear functions
$h^\pm(x)=x\cdot\rho^\pm-g^*(\rho^\pm)$ must agree on
the set $X(\rho^-,\rho^+)$. It is worth mentioning that
if a function $\rho$ is given by \eqref{eq1.5}, then its {\em weak
derivative} $D\rho$ is a matrix measure that is concentrated on the union of edges
$X(\rho^-,\rho^+),\ \rho^\pm\in \mathcal{S}$. Indeed,
\[ 
D\rho(dx)=\sum_{\rho^\pm\in \mathcal{S}}\1\big(x\in X(\rho^-,\rho^+)\big)
\ \big[(\rho^+-\rho^-)\otimes n(\rho^-,\rho^+)\big]\ m(dx),
\]
where $dm$ denotes the one-dimensional Lebesgue measure (on the union
of the edges), 
 and $n(\rho^-,\rho^+)$ is a unit normal that is orthogonal to
$\tau(\rho^-,\rho^+)$, and is pointing from the $X(\rho^-)$
side to the $X(\rho^+)$ side of $X(\rho^-,\rho^+)$. For $\rho$ 
to be a gradient $\nabla g$, the matrix $D\rho$ must be symmetric. The matrix $D\rho$ is symmetric if and only if  \eqref{eq1.6} holds. For the convexity of 
$g$, we need $D\rho\ge 0$, which is equivalent to saying that
$\rho^+-\rho^-$ is pointing from the $X(\rho^-)$
side to the $X(\rho^+)$ side of $X(\rho^-,\rho^+)$. In other words,
\[
n(\rho^-,\rho^+)=\frac{\rho^+-\rho^-}{|\rho^+-\rho^-|}.
\]
We summarize our discussion in the next definition.

\bs\noi
{\bf Definition 1.2(i)} Let $\L$ be a convex polytope in $\bR^2$. By a {\em (generic Laguerre) tessellation} of $\L$ we mean
a countable collection $\bX=\big\{\big(\rho,X(\rho)\big):\ \rho\in \mathcal{S}\big\}$ such that 
\bi
\item Each $X(\rho)$ is a convex polytope, and 
\[
\bigcup_{\rho\in \mathcal{S}}X(\rho)=\L.
\]
We refer to each $X(\rho)$ as a {\em cell} of $\bX$.
\item If $\rho^\pm\in \mathcal{S}$ are distinct, and $X(\rho^-,\rho^+)\neq\emptyset$, then 
$X(\rho^-,\rho^+)$ is a line segment orthogonal 
to $\rho^+-\rho^-$. We refer to such $X(\rho^-,\rho^+)$ as an
{\em edge} of $\bX$.
\item If $\rho^\pm,\rho^*\in \mathcal{S}$ are distinct, and 
\[
X(\rho^-,\rho^*,\rho^+):=X(\rho^-)\cap X(\rho^*)\cap
X(\rho^+)\neq\emptyset,
\]
 then $X(\rho^-,\rho^*,\rho^+)$ consists of a single point. We refer 
to this point as a {\em vertex} of the tessellation $\bX$.
\item For each edge $X(\rho^-,\rho^+)$,
the vector $\rho^+-\rho^-$ is pointing from the $X(\rho^-)$
side to the $X(\rho^+)$ side of $X(\rho^-,\rho^+)$.
\ei
We write $\cX(\L)$ for the set of all generic tessellations of $\L$.

\ms\noi
{\bf(ii)} We set
\[
\G=\big\{(\rho^-,\rho^+)\in \bR^2\x \bR^2:\ \rho^-\neq\rho^+\big\}.
\]
By an {\em orientation} $\tau$, 
we  mean a continuous function $\tau:\G\to\bR^2$, such that \eqref{eq1.6} holds. 
\qed

\subsection{Gibbsian measures on $\cX$ or $\cC_0$}

Our definition of generic tessellations can be readily extended to any dimension. Observe that when $d=1$, each {\em cell} $X(\rho)$ is an interval on which the nondecreasing function $\rho(\cdot)$ is constant.
As a natural candidate for a measure on  $\widehat\cC_0$, we may pick a continuous {\em kernel}
$f(x, \rho^-,d\rho^+)$ (a measure in $\rho^+$ for every $(x,\rho^-)\in
\bR^2$) with
\[
\lambda(x,\rho^-):=\int_{\rho^-}^\i f(x,\rho^-,d\rho^+)<\i,
\]
and set $\nu^{f}$ to be the law of an 
{\em inhomogeneous Markov process} $x\mapsto \rho(x)$
with a jump rate given by $f$. In other words, as $x$ increases,
the Markov process $\rho(x)$ has an infinitesimal generator
\[
\cL_x F(\rho^-)=\int \big(F(\rho^+)-
F(\rho^-)\big)\ f(x,\rho^-,d\rho^+).
\]
Given an initial law $\ell^0(d\rho)$, and $a^-\in\bR$,
we may construct a measure
on piecewise constant functions $\rho:[a^-,\i)\to\bR$,
 so that $\rho(a^-)$
is selected according to $\ell^0$, and evolves in a Markovian fashion 
with the generator $\cL_x$. Note that if the law of $\rho(x)$ is given by
$\ell(x,d\rho)$, then $\ell$ satisfies the forward equation
$\ell_x=\cL^*_x\ell$, subject to the initial condition $\ell(a^-,d\rho)
=\ell^0(d\rho)$. Given $a^+>a^-$, we may {\em also} interpret
$\big(\rho(x):\ x\in[a^-,a^+]\big)$, as a Markov process 
that starts at time $a^+$ with law $\ell(a^+,d\rho)$, and evolves 
backward in a Markovian fashion as we {\em decrease} $x$. 
In order to have the same law on $\big(\rho(x):\ x\in[a^-,a^+]\big)$, the jump kernel of this (backward) Markov process must be selected appropriately (see \eqref{eq5.3} below). 

\ms
In the same manner, we wish to construct a Gibbsian measure
$\nu^f$ on $\widehat\cC_0$ for a bounded continuous kernel
$f(x,\rho^-,d\rho^+)$, 
which is a measure in $\rho^+$ for any given
$(x,\rho^-)$, and depends continuously on
$x$. We carry out this construction when $d=2$
in the present article.
Though, as our method of construction suggests, it seems plausible that one can carry out similar constructions in higher dimensions in an inductive manner. Indeed our method of construction takes advantage of the fact that we already have a natural candidate for such measures in dimension one, namely Markov jump processes.
Once our measures are constructed for $d=2$, we may use them 
to construct Gibbsian measures in dimension $3$ in a similar manner.
We should mention that if $\nabla g\in\widehat\cC_0(\bR^2)$, and
$\hat\rho(x,t)=(u_x,u_t)(x,t)$, then 
$\hat\rho\in\widehat\cC_0(\bR^{2}\times[0,\i))$.
A Gibbsian choice of $\nabla g$ leads to a probability measure on
$\widehat\cC_0(\bR^{d}\times[0,\i))$, which has a similar flavor as our construction 
when $d=2$.  We expand on this in Section 1.6.

Given a kernel $f(x,\rho^-,d\rho^+)$ in $\bR^2$, we wish to construct a measure on $\widehat\cC_0$ so that $f(x,\rho^-,d\rho^+)$
represents the rate at which $\rho^-$ changes to $\rho^+$ as we cross an edge of $X(\rho^-)$ at a point $x$. To achieve this we adopt the following strategy:

\ms\noi
{\bf(i)}
We take a convex planar set $\L$
(for example a box), and construct a measure $\nu^{f,\L}$ on the set of tessellations of $\L$. This is carried out by first constructing
the one dimensional tessellation
\[
\big\{\big(\rho,X(\rho)\cap \p\L\big):\ \rho\in \mathcal{S}\big\},
\]
in a Markovian fashion with a jump rate that is expressed in terms of $f$. We then use the information coming from the boundary to build the tessellation inside. 
More edges will be added inside $\L$ in a Markovian fashion.

\ms\noi
{\bf(ii)} We then show that when $f$ satisfies a suitable kinetic
 PDE,
the measures $\nu^{f,\L}$ are consistent as we enlarge $\L$.

 \ms
 As we attempt  to carry out the above strategy, we encounter two 
 problems. Our treatment of these problems are responsible
 for our final recipe of our Gibbsian measures. 
 
 \ms
 \noi
 {\bf Problem 1.} According to our strategy, we would like to construct $\rho(x)$ on the boundary as a Markov process. Imagine that we start from a point $a\in\p\L$ and select
 a slope $\rho^-$ for $\rho(a)$. We then move counterclockwise on
 the boundary and change $\rho(x)$ as a jump process. If a jump from $\rho^-$ to $\rho^+$ occurs at a point $b\in\p\L$ to the right of $a$, then there will be an edge of our tessellation separating
 $X(\rho^-)$ from $X(\rho^+)$. If no other jump occurs as we traverse the whole boundary, the restriction of the desired tessellation to the set $\L$ consists of exactly two cells. However 
 the edge $X(\rho^-,\rho^+)$ intersects the boundary at a second point $b'$ that was not a jump point of our Markov process.
  In other words, the restriction of the desired tessellation to $\p \L$ cannot be realized as a Markov process.
 
 \ms\noi
 {\bf Our Remedy:} We resolve this problem by giving an orientation to the edges. This orientation will be used to decide what points of the boundary tessellation will be created as a Markov jump process. In other words, given an orientation $\tau$, we consider a Markov process on the boundary with the jump rate  
 \[
 \big(\tau(\rho^-,\rho^+)\cdot n(x)\big)^+\ f(x,\rho^-,\rho^+),
  \]
 where $n(x)$ denotes the inward unit normal to $\p\L$ at $x$.
This Markov process would allow us to determine the {\em entering edges} only, not the {\em exiting edges}. 
  
 \ms
 \noi
 {\bf Problem 2.} After determining all the entering edges, we need 
 to use them to build our tessellation inside $\L$. These edges may intersect to produce vertices. How can this be done in an orderly manner?

\ms\noi
 {\bf Our Remedy:} Given a fixed direction $v$, we expect/insist 
 \[
 \rho^\pm\in R,\ \ \ \rho^+\neq\rho^-\ \ \ \implies\ \ \ \  (\rho^+-\rho^-)\cdot v\neq 0,
 \]
 in the support of our measure. This is equivalent to saying that 
$\tau(\rho^-,\rho^+)$
 is not parallel to $v$. Without loss of generality, we may choose $v=e_1=(1,0)$, so that $\tau$ is never horizontal. 
A continuous choice of $\tau$ forces a 
 fixed sign for $\tau(\rho^-,\rho^+)\cdot e_2$. 
Without loss of generality, we
require that $\tau(\rho^-,\rho^+)
 \cdot e_2>0$. For the sake of definiteness, we choose
\be\la{eq1.7}
\tau(\rho^-,\rho^+)=\big(-[\rho^-,\rho^+],1\big),\ \ \ \
{\text{where}}\ \ \ [\rho^-,\rho^+]:=\frac {\rho^+_2-\rho^-_2}
{\rho^+_1-\rho^-_1}.
\ee
We may treat $\tau$ as a velocity for a point/particle that is created
at a point $b$ on the boundary, and its trajectory determines the edge emanating from $b$. Here, we are treating $x_2$
as a time parameter. We will use this time parameter to order the creation of particles, and the occurrence of particle collisions 
as we increase $x_2$.
\qed

\bs 
We now have all the ingredients to describe our Gibbsian measure
in a convex set $\L$. For this construction, we use the orientation
\eqref{eq1.7}. 
 Note that because of our choice of $\tau$, we may talk about a cell $X(\rho^-)$ (respectively $X(\rho^+)$) that lies on the left (respectively right) of the edge $\tau(\rho^-,\rho^+)$. With this convention, the 
fourth condition
in Definition 1.2{\bf(i)} is equivalent to saying that for each
$(x,\rho^-)$, the support of the measure $f(x,\rho^-,d\rho^+)$
is contained in the set
\be\la{eq1.8}
R(\rho^-):=\big\{\rho^+=(\rho^+_1,\rho^+_2):\ \rho^+_1>\rho^-_1\big\}.
\ee
To simplify our presentation, we give a precise recipe for our measure when $\L$ is a box:
\[
\L=\L(a^-,a^+,t_0,t_1):=
[a^-,a^+]\x [t_0,t_1].
\]
Because of our choice of $\tau$, there will be no entering edge from the top side of the 
box. To ease the notation, we simply write $t$ for the second coordinate. Pictorially,  we may  decorate
 each edge with an arrow that always points upward (see Figure 1 below). Our Gibbsian measure will be supported on  generic tessellations
of which each vertex is of degree 3. With an orientation at our disposal, 
we may interpret each vertex as either a {\em coagulation point} (when two intersecting edges are replaced with one edge as time increases), or a {\em fragmentation point} (when an edge splits into two edges as time increases). With these interpretations, 
we may fully determine the tessellation inside $\L$ in terms a collection of particles that travel according to their velocities, and may experience coagulation and fragmentation. 
More precisely, the  function 
$x_1\mapsto \rho(x_1,t)$
can be expressed as
\[
 \rho(x_1,t)=
 \sum_{i\in J(t)}\rho^i(t)\ \1\big(x_1\in(z_{i}(t),z_{i+1}(t)\big)\big),
 \]
 with the interpretation that $z_i(t)$ represents the position of the $i$-th particle. Writing $q_i=(z_i,\rho^i)$, and 
 $\bq(t)=\big(q_i(t):\ i\in J(t)\big)$, the dynamics of $\bq$ can be conveniently described as a Markov process. 

\bs\noi
{\bf Definition 1.3(i)} Given a pair $\rho^\pm\in \bR^2$, we write
$\rho^-\prec\rho^+$ if $\rho^+\in R(\rho^-)$, where $R(\rho^-)$
was defined in \eqref{eq1.8}. Similarly, we define the set $L(\rho^{+})$ as 
\[
L(\rho^{+}):=\{ \rho^{-} \in \bR^2 : \rho^{-} \prec \rho^{+}\},
\]
and the set $D(\rho^{-},\rho^{+})$ for $\rho^{-} \prec \rho^{+}$ to be
\[
D(\rho^{-},\rho^{+}):=\{ \rho^{*} : \rho^{-} \prec \rho^{*} \prec \rho^{+} \}.
\]
\ms\noi
{\bf(ii)} We write $\D=\bigcup_{n=0}^\i\D_n$, where $\D_n$ denotes
the set
of $\bq=\big(q_0,\dots,q_n\big)$ such that 
$q_i=(z_i,\rho^i)\in\bR^3$,
and
\[
z_0=a^-<z_1<\dots <z_n<z_{n+1}:=a^+,
\ \ \ \ \rho^0\prec\rho^1\prec\dots\prec\rho^n.
\]

\ms\noi
{\bf(iii)} We write $\cM$ for the set of measures on $\bR^2$, and equip $\cM$ with the topology of weak convergence. The set of probability measures is denoted by $\cM_1$.
We write $\cF(\L)$ for the set of kernels 
\[
f:\L\x\bR^2\to\cM,
\]
with the following properties:
\bi
\item The map $(x,\rho^-)\mapsto f(x,\rho^-,d\rho^+)$ is measurable, and 
\[
\sup_{(x,\rho^-)}\int \big|\tau(\rho^-,\rho^+)\big|\ f(x,\rho^-,d\rho^+)<\i.
\]
\item For every $(x,\rho^-)$,
\[
f\big(x,\rho^-,\bR^2\setminus R(\rho^-)\big)=0.
\]
 \ei
 
 \ms\noi
{\bf(iv)} Given $\b\in\cM$, we write $\cF(\b,\L)$ for the set of
$f\in\cF(\L)$ such that $f(x,\rho^-,d\rho^+)\ll  \b(d\rho^+)$.
With a slight abuse of notation, we write
$f(x,\rho^-,\rho^+)$ for the Radon-Nikodym derivative of
$f(x,\rho^-,d\rho^+)$ with respect to $\b$.
 
\ms\noi
{\bf(v)} Given  constants $P^\pm$, with $P^-<P^+$, and 
 positive constants $V_\i$ and $\d_0$, we put
\[
\G^{V_\i}=\left\{(\rho^-,\rho^+)\in [P^-,P^+]^2\x [P^-,P^+]^2:\ 
\rho^-\prec \rho^+,\ \big|[\rho^-,\rho^+]\big|\le V_\i\right\},
\]
and write
$\cF(\b,\L,V_\i,\d_0)$ 
for the set of 
$f\in\cF(\b,\L)$ such that 
\begin{align*} 
&x\in\bR^2,\ (\rho^-,\rho^+)\in \G^{V_\i}\ \ \
\implies\ \ \  f(x,\rho^-,\rho^+)\ge\d_0,\\
&x\in\bR^2,\ (\rho^-,\rho^+)\notin \G^{V_\i}\ \ \
\implies\ \ \  f(x,\rho^-,\rho^+)=0.
\end{align*}
\qed

\bs

We now introduce some definitions that will allow us to construct a Markov process taking values in $\Delta$.
 
 \ms\noi
{\bf Definition 1.4 (i)} Given $\ell^0\in\cM_1$, $f\in\cF(\b,\L)$, and $t_0\in\bR$,
 we write $\g(d\bq;t_0,\ell^0,f)$ for a probability measure on $\D$ that represents a Markov jump process with the jump rate 
 $ f((x_1,t_0),\rho^-,d\rho^+),$
 and the initial law  $\ell^0$. More precisely, if $\bq$ is selected
according to $\g(d\bq;t_0,\ell^0,f)$ , and the function $\rho(\cdot;\bq)$
is defined by
 \[
 \rho(x_1;\bq)=\sum_{i=0}^n\rho^i\
 \1\big(z_i\le x_1<z_{i+1}\big),
 \]
then $\rho(a^-;\bq)=\rho^0$ is distributed according to $\ell^0$,
and the Markov jump 
process $x_1\mapsto \rho(x_1;\bq)$ makes its $i$-th jump
from $\rho^{i-1}$ to $\rho^i$ at {\em time} $z_i$, with the rate 
$ f((z_i,t_0),\rho^{i-1},d\rho^{i}).$

\ms\noi
{\bf (ii)} Given $\ell^0\in \cM_1$ and $f\in\cF(\L,\b,V_\i,\d_0)$ as 
in the Definition 1.3{\bf(v)}, we define a Markov
process $\big(\bq(t):\ t\ge t_0\big)$ that takes value in
 the set $\D$. This Markov
process induces a function 
\[
\rho:\L\to\bR,\ \ \ \ \rho(x_1,x_2)=\rho(x_1,t):=\rho(x_1;\bq(t)),
\]
that belongs to $\widehat\cC_0(\L)$. The (initial) law of $\bq(t)$ at 
$t=t_0$ is given by $\g(d\bq;t_0,\ell^0,f)$. This process induces a
probability measure on $\widehat\cC_0(\L)$, that is denoted by 
$\nu^{f,\L}=\nu^{\ell^0,f,\L}$.  The Markovian dynamics of $\bq(t)$ is as follows:

\ms\noi
{\bf 1.} The particle $z_i$ travels with velocity
$-[\rho^{i-1},\rho^i]$. When $z_1$ 
reaches $a^-$, or  $z_n$ reaches respectively $a^+$, the number of particles reduces by one. In the former case, we relabel $(z_i,\rho^i)$ as
$(z_{i-1},\rho^{i-1})$ for $i\ge 1$.

\ms\noi
{\bf 2.} If at some time $t$, we have $z_i(t)=z_{i+1}(t)$, then
we remove the $i$-th particle from the system, and relabel
$(z_j,\rho^j)$ as
$(z_{j-1},\rho^{j-1})$ for $j>i$.

\ms\noi 
{\bf 3.} At the boundary point $a^-$, 
the function $t\mapsto \rho(a^-,t)$ can change from 
$\rho^0$ to $\rho^*$ with the rate
\[
[\rho^0,\rho^*]^-\frac{\ell(a^-,d\rho^*)f\big((a^-,t),\rho^*,d\rho^0\big)}{\ell(a^-,d\rho^0)},
\]
When this happens, we relabel $(z_i,\rho^i)$,   as
$(z_{i+1},\rho^{i+1})$, for $i\ge 0$, and declare $\rho^*$ to be our new
$\rho^0$.

\ms\noi
{\bf 4.} At the boundary point $a^+$, 
the function $t\mapsto \rho(a^+,t)$ can change from 
$\rho^n$ to $\rho^*$ with the rate
\[
[\rho^n,\rho^*]^+f\big((a^+,t),\rho^n,d\rho^*\big),
\]
When this happens, a new particle has been born
at $a^+$. That is, we now have $n+1$ many particles with
$z_{n+1}=a^+$, and $\rho^{n+1}=\rho^*$.

\ms\noi
{\bf 5.} The $i$-th particle can {\em fragment} into two particles.
This occurs at time $t$ with the rate density
\be\la{eq1.9}
\s\big(\rho^{i-1},\rho^*,\rho^{i}\big)^-\ \frac
{f\big((z_i,t),\rho^{i-1},\rho^{*}\big)f\big((z_i,t),\rho^{*},\rho^{i}\big)}{f\big((z_i,t),\rho^{i-1},\rho^{i}\big)},
\ee
where
\be\la{eq1.10}
\s(\rho^-,\rho^*,\rho^+)=[\rho^*,\rho^+]-[\rho^-,\rho^*].
\ee
By fragmentation we mean that the particle $(z_{j},\rho^{j})$ is relabeled
as $(z_{j+1},\rho^{j+1})$ for $j\ge i$, and the $i$-th particle at 
the location $z_i$
is associated with a new label
$\rho^{i}=\rho^*$. 

\ms\noi
{\bf(iii)} We write $\cM_G(\L)$ for the set of measures of the form
$\nu^{\ell^0,f,\L}$, as we vary $\ell^0\in\cM_1$ and 
$f\in \cF(\b,\L,V_\i,\d_0)$.
\qed

\bs
Because of our choice \eqref{eq1.9}, the fragmentation mechanism in {\bf 5} is the time reversal of the coagulation mechanism in {\bf 2}. The choice \eqref{eq1.9} plays an essential role in the validity of our first main result, namely,
 the consistency of our measures $\nu^{f,\L}$ as we vary $\L$.  We remark that the fragmentation occurs only when $\s<0$, so that the resulting particles move away from each other, with the $i+1$-th particle to the right of the newly born particle.

\tikzset{->-/.style={decoration={  
markings,
  mark=at position .5 with {\arrow[scale=2.5,>=stealth]{>}}},postaction={decorate}}}
\begin{figure}
\begin{center}
\begin{tikzpicture}
\centering
\draw (-5,0) -- (5,0)--(5,13)--(-5,13)--(-5,0);
\draw[->-] (-3.5,0)-- (-2.5,2.5);
\draw[->-] (-1,0)-- (-2.5,2.5);
\draw[->-] (-2.5,2.5)--(-2.75,4.5) ;
\draw[->-] (-2.75,4.5)--(0.5,6);
\draw[->-] (-2.75,4.5)--(-4.25,7.5);
\draw[->-] (-5,6.5)--(-4.25,7.5);
\draw[->-] (3.5,0)--(2.5,3.5);
\draw[->-] (5,2)--(2.5,3.5);
\draw[->-] (2.5,3.5)--(0.5,6);
\draw[->-] (0.5,6)--(1,8);
\draw[->-] (1,8)--(4,11.5);
\draw[->-] (1,8)--(-1,11);
\draw[->-] (-4.25,7.5)--(-3.5,10);
\draw[->-] (-3.5,10)--(-4.1,13);
\draw[->-] (-3.5,10)--(-1,11);
\draw[->-] ((-1,11)--(-0.2,13);
\draw[->-] (5,9.5)--(4,11.5);
\draw[->-] (4,11.5)--(3.75,13);
\node at (-5,-0.5) {$(a^{-},t_0)$};
\node at (5,-0.5) {$(a^{+},t_0)$};
\node at (5,13.5) {$(a^{+},t_1)$};
\node at (-5,13.5) {$(a^{-},t_1)$};
\node at (-3.5,-0.5) {$z_1$};
\node at (-1,-0.5) {$z_2$};
\node at (3.5,-0.5) {$z_3$};
\node at (-4,3.5) {$\rho_0$};
\node at (-2.5,1) {$\rho_1$};
\node at (0.3,3) {$\rho_2$};
\node at (4.2,1.5) {$\rho_3$};
\node at (-4.5,9.75) {$\rho_{-1}$};
\node at (3.2,6.5) {$\rho_4$};
\node at (-1.5,7.5) {$\rho_{02}$};
\node at (1.25,11.2) {$\rho_{024}$};
\node at (-2.5,11.5) {$\rho_{(-1)02}$};
\node at (4.5,11.6) {$\rho_{5}$};
\draw[blue,fill=blue] (2.5,3.5) circle (0.1);
\draw[red,fill=red] (-2.75,4.5) circle (0.1);
\end{tikzpicture}
\end{center}
\caption{The blue dot represents the coagulation of the particles with labels $(\rho_2,\rho_3)$ and $(\rho_3,\rho_4)$ into the particle with label $(\rho_2,\rho_4)$. The red dot represents the fragmentation of the particle with label $(\rho_0,\rho_2)$ into two particles of respective labels $(\rho_0,\rho_{02})$ and $(\rho_{02},\rho_2)$.}
\end{figure}
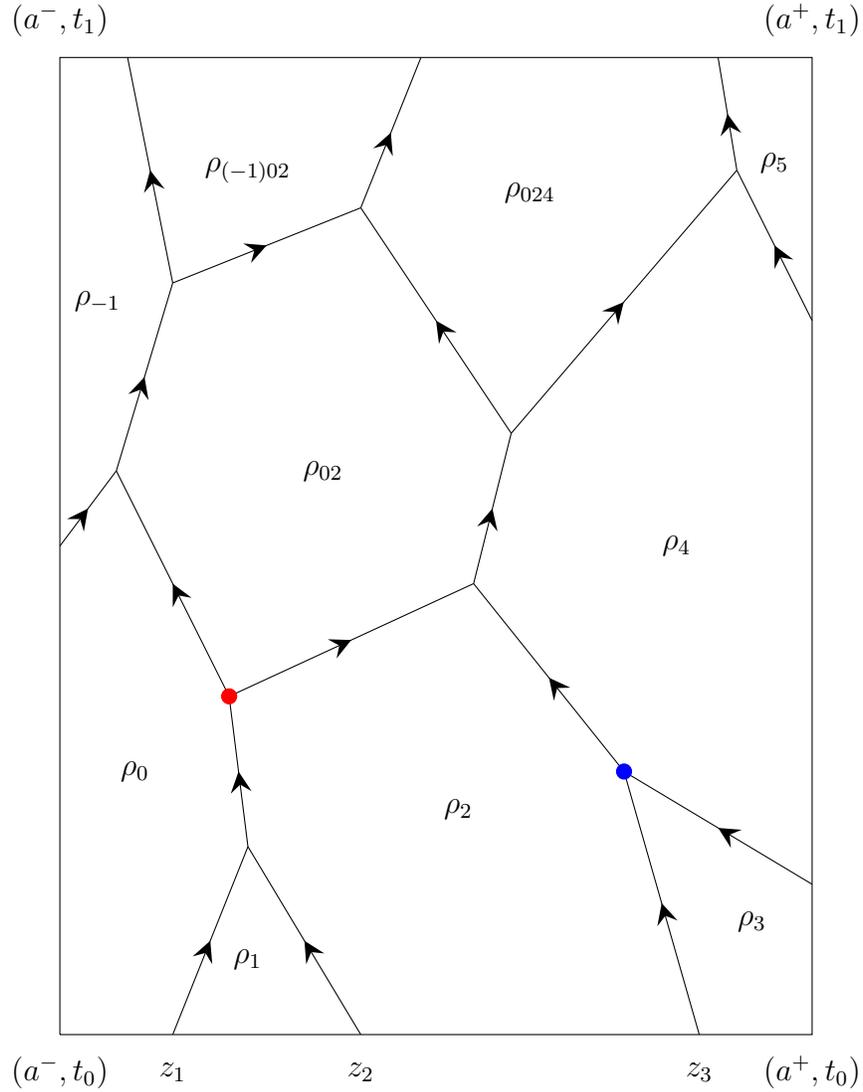

\bs

\subsection{Consistency}

We now turn our attention to the question of the consistency 
of our  measures $\nu^{f,\L}$, as we vary the sides of
the box $\L=\L(a^-,a^+,t_0,t_1)$. We first vary the horizontal 
sides. If 
\[
\L=\L(a^-,a^+,t_0,t_1), \ \ \ \ \ \L'=\L(a^-,a^+,t,t_1),
\]
 with $t\in(t_0,t_1)$, the consistency of
$\nu^{f,\L}$ and $\nu^{f,\L'}$, requires that the process
$x\mapsto \rho(x,t)$ under $\nu^{f,\L}$ to be a Markov process
with the jump rate 
\be\la{eq1.11}
f\big((x_1,t),\rho^-,d\rho^+\big).
\ee
This turns out to be equivalent to the requirement that the kernel
$f$  satisfies a kinetic type PDE  of the form
\be\la{eq1.12}
\tau(\rho^-,\rho^+)\cdot f_x(x,\rho^-,d\rho^+)=
Q(f)(x,\rho^-,d\rho^+),
\ee
for a suitable quadratic function $Q$. To ease the notation, 
we suppress the dependence on $x$ in our notations, and write
\be\la{eq1.13}
\tau f=(-\a f,f):=(-f^2,f^1),
\ee
where $\a(\rho^-,\rho^+)=[\rho^-,\rho^+]$.
With these conventions, the operator $Q$ equals 
\[
Q(f)=Q(f^1,f^2)=Q^+(f^1,f^2)-Q^-(f^1,f^2),
\]
 with
\begin{align}\nonumber
Q^+(f^1,f^2)(\rho^-,d\rho^+)=&(f^1*f^2)(\rho^-,d\rho^+)-(f^2*f^1)
(\rho^-,d\rho^+),\\
Q^-(f^1,f^2)(\rho^-,d\rho^+)=&\left(A(f^2)(\rho^+)-A(f^2)(\rho^-)\right)
f^1(\rho^-,d\rho^+)\la{eq1.14}\\
&-\left(A(f^1)(\rho^+)-A(f^1)(\rho^-)\right)f^2(\rho^-,d\rho^+).\nonumber
\end{align}
Here,
\begin{align*}
(h*k)(\rho^-,d\rho^+):&=\int h(\rho^-,d\rho^*)\
k(\rho^*,d\rho^+),\ \ \ \ \
A(h)(\rho):=\int h(\rho,d\rho^*).
\end{align*}
\qed

\bs\noi
{\bf Remark 1.1}
Given a pair of kernels $(f^1,f^2)$, define the quadratic operator
\[
\cQ(f^1,f^2)=f^1*f^2-A(f^1)\otimes f^2-f^1\otimes A(f^2),
\]
where
\[
(h\otimes k)(\rho^-,d\rho^+)=h(\rho^-)k(\rho^-,d\rho^+),\ \ \ \ 
(k\otimes h)(\rho^-,d\rho^+)=h(\rho^+)k(\rho^-,d\rho^+).
\]
Then \eqref{eq1.12} can be written as 
\[
div\ (\tau f)=f^1_{x_2}-f^2_{x_1}=\cQ(f_1,f_2)-\cQ(f_2,f_1).
\]
\qed

\bs
As our first main result, we verify the consistency of the measures
$\nu^{f,\L}$ as we vary the horizontal sides of $\L$.
\ms

\bth\la{th1.1} Assume that the kernel $f\in \cF(\b,\L,\L_\i,\d_0)$ is a $C^1$ function, and satisfies \eqref{eq1.12}, for
$\L=\L(a^-,a^+,t_0,t_1)$. Assume that $\ell^0\ge \d_0$.
	Then for every $x_2\in[t_0,t_1]$, the 
law of the function $x_1\mapsto \rho(x_1,x_2)$ with respect to
$\nu^{\ell^0,f,\L}$ coincides with the law of 
a Markov jump process with the jump
rate $f\big((x_1,x_2),\rho^-,d\rho^+\big)$.
\et

\bs\noi
{\bf Example 1.1} Given a continuous function $K:\bR\to\bR$, we may consider
\[
\b(d\rho)=\b(d\rho_1,d\rho_2)=\d_{K(\rho_1)}(d\rho_2)\ d\rho_1,
\]
which is a measure that is supported on the graph of the function $K$. In this case, the corresponding convex function $g$, with $\nabla g=\rho$ satisfies the PDE
\be\la{eq1.15}
g_{x_2}=K\big(g_{x_1}\big),
\ee
 inside the cells of the corresponding tessellation.
We put 
\begin{align*}
\tilde f^i(x_1,x_2,\rho^-_1,\rho^+_1)&=
f^i\big(x_1,x_2,\rho^-_1,K(\rho^-_1),\rho^+_1,K(\rho^+_2)\big),
\end{align*}
 write $\tilde f$ for $\tilde f^1$, and write
$\mu^{\ell^0,\L,\tilde f,K}$ for the corresponding measure
$\nu^{\ell^0,\L, f}$. We write $\widehat \cM_G$ for the set of
such measures as we vary $f$ and $K$.

\ms\noi
{\bf(i)} When $K$ is a convex function, then $g$ is a viscosity solution
of \eqref{eq1.15}. Moreover, there would be no fragmentation because
$\s\ge 0$. If we also assume that $f$ is independent of $x_2$ and that $K$ is increasing, then Theorem 1.1 was established in
[KR1], confirming affirmatively a conjecture of Menon and Srinivasan [MS]. We may rephrase Theorem 1.1 as follows:
If $g$ solves the Hamilton-Jacobi PDE \eqref{eq1.15} in $d=1$, and initially
$g_{x_1}(x_1,t_0)$ is a Markov jump process with the rate density
$\tilde f^1(x_1,t_0,\rho_1^-,\rho_1^+)$, then for every $t>t_0$,
the process $g_{x_1}(x_1,t)$ is also 
a Markov jump process with the rate density
$ \tilde f^1(x_1,t,\rho_1^-,\rho_1^+)$.

\ms\noi
{\bf(ii)} If we assume that $K$ is concave, then there would be
no collision as we increase $t=x_2$ because $\s\le 0$. In fact 
$g$ is a viscosity solution for a final value problem i.e., as we reverse (decrease) time  $x_2$. 
We may rephrase Theorem 1.1 in this case
as a statement for the reversed dynamics: 
If we reverse $x_2$ in part {\bf(i)}, then the dynamics can be described as a particle system with stochastic fragmentation; the fragmentation rate is given by \eqref{eq1.9}.

\ms\noi
{\bf(iii)} If $K$ is neither convex, nor concave, then $g$ is not a viscosity solution of the PDE \eqref{eq1.15}
no matter what direction for the coordinate $x_2$ is adopted.
\qed

\ms\noi
{\bf Remark 1.2}
Observe that a choice of an orientation for edges allowed us to have a natural time direction for the Markov processes on the boundary sides of the box $\L$, and the dynamics
 inside the box.  For the undirected tessellation, this choice is irrelevant, and there should be a formulation of the consistency criteria that is independent of the orientation. As an illustration,  we will
 demonstrate in Proposition 4.2 below how reversing a direction, or interchanging coordinates can be performed on the solutions of 
 \eqref{eq1.12}. 
\qed

\bs
As we mentioned before, the Markov process $\bq(t)=\big((z_i(t),\rho_i(t)):\ i\in J(t)\big)$,
yields a random tessellation 
\[
\bX_\L=\left\{(\rho,X(\rho)):\ \rho\in\cS_\L\right\},
\]
 of the box $\L$. The set $\cS_\L$ is simply defined by
\[
\cS_\L=\left\{\rho^i(t):\ t\in[t_0,t_1],\ i\in J(t)\right\},
\]
and the cells of $\bX_\L$ are the connected components of the set
\[
\L\setminus\left\{(z_i(t),t):\ t\in[t_0,t_1],\ i\in J(t)\right\}.
\]
The law of $\bX_\L$ is denoted by $\eta^{\ell^0,f,\L}$.

\ms
\bp\la{pro1.1} Under the assumptions of Theorem 1.1, we have that 
$\eta^{\ell^0,f,\L}\left(\chi_\L\right)=1.$
In words, the tessellation $\bX_\L$ is generic in the sense of Definition
{\bf{1.2(i)}}, with probability one with respect to $\eta^{\ell^0,f,\L}$.
\ep

The proof of this Proposition is rather straightforward and follows from
our Proposition 3.1 in Section 3.4.

\bs
We next examine the question of the consistency as we vary
the vertical sides of $\L$. Note however that although the boundary 
dynamics on the lower side is Markovian, the dynamics on the lateral sides may depend on the configuration inside $\L$. This can be avoided
if we assume that $\tau$ always points to the left. 

\ms

\bth\la{th1.2} Let $f$ and $\L$ be as in Theorem 1.1. Assume
that the measure $f(x,\rho^-,d\rho^+)$ is supported in the set
\be\la{eq1.16}
R_0(\rho^-)=\big\{\rho^+=(\rho^+_1,\rho^+_2):\ 
\rho^+_2> \rho^-_2,\ \rho^+_1> \rho^-_1\big\}.
\ee
for every $(x,\rho^-)$. Then  the law of the 
function $x_2\mapsto \rho(x_1,x_2)$ with respect to
$\nu^{\ell^0,f,\L}$ coincides with the law of 
a Markov jump process with the jump
rate given by $[\rho^-,\rho^+]f\big((x_1,x_2),\rho^-,d\rho^+\big)$.
\et

 Our assumptions on $f$ allow us to reduce Theorem 1.2 from 
Theorem 1.1. The details can be found in Section 5.

\bs\noi
{\bf Remark 1.3}
 More generally,  define 
\be\la{eq1.17}
R_c(\rho^-)=\big\{\rho^+=(\rho^+_1,\rho^+_2):\ 
\rho^+_2 -\rho^-_2+c(\rho^+_1- \rho^-_1)>0,\ \rho^+_1>
 \rho^-_1\big\},
\ee
and assume that there exists $c\ge 0$ such that  the measure $f(x,\rho^-,d\rho^+)$ is supported in the set
$R_c(\rho^-)$, for every $(x,\rho^-)$. Define 
$S_c(x_1,x_2)=(x_1+cx_2,x_2)$, and
\[
S_cg(x_1,x_2)=g(x_1+cx_2,x_2),\ \ \ \ T_c(\rho_1,\rho_2)=
(\rho_1,\rho_2+c\rho_1),
\]
so that if $\rho=\nabla g,$ then
\[
\nabla (S_cg) (x)=(T_c\rho)(x_1+cx_2,x_2)=:\rho'(x).
\]
We also define $f':=T_c^\sharp f$, i.e.,  for every bounded continuous function $\var :\bR\to\bR$, we have
\[
\int \var (\rho^+)\ f'(x,\rho^-,d\rho^+)=\int (\var\circ T_c)(\rho^+)
 f(x_1+cx_2,x_2, T_c\rho^-,d\rho^+).
\]
Then the kernel $f'(x,\rho^-,d\rho^+)$ is supported in the set $R_0(\rho^-)$, and Theorem 1.2 is applicable to $f'$. As a result,
under $\nu^{\ell^0,f',\L}$, 
 the law of $\rho$, restricted to a line of slope $c$ is a Markov jump process. This in turn implies the consistency for the measures
$\nu^{f,\L'}$, provided that $\L'=S_c\L$, for a box $\L$.
\qed

\bs\noi
{\bf Example 1.1(iv)} {\em (continued)} 
 When $K$ is increasing and convex, we claim
that for each $x_1$, the process $x_2\mapsto \rho(x_1,x_2)$
is a Markov jump process with the rate $\tilde f^2(x_1,t,\rho_1^-,\rho_1^+)$. 
To see this observe that we may write 
$g_{x_1}=K^{-1}\big(g_{x_2}\big)$, which suggests that we should regard $x_1$ as the time variable now. With this choice of time, we 
now have a scenario that resembles  Example 1.1(i), except for few non-essential differences: We are initially at $x=a^+$, and go backward by {\em decreasing} $x_1$. The function $K^{-1}$
is now concave, which implies that the convex function $g$
is a viscosity solution for the final-value Hamilton-Jacobi PDE
$g_{x_1}=K^{-1}\big(g_{x_2}\big)$. The initial jump process
with the jump rate density $\tilde f^2(a^+,t,\rho_1^-,\rho_1^+)$,  evolves
to a jump process with the jump rate density
$\tilde f^2(a,t,\rho_1^-,\rho_1^+)$
as we decrease $x_1$ from $a^+$ to $a$. 
\qed

\subsection{The invariance of $\cM_G$ and $\widehat \cM_G$}

We now examine the question of the invariance of the set $\cM_G$ under the flow 
$\widehat \Phi$ of
the Hamilton-Jacobi PDE
\be\la{eq1.18}
u_t=H(u_x),\ \ \ \ u(x,0)=g(x),
\ee
with $H:\bR^2\to\bR$ convex, and $g\in\cC_0(\L)$.

From Example 1.1{\bf(i)}, 
we already know that $\cM_G(\L)$ is invariant under 
$\widehat\Phi$, when $d=1$. In details, we choose $\b$ to be
the one-dimensional 
Lebesgue measure, and pick $f^0\in \cF(\L,\b)$, with $\L=[a^-,a^+]$
or $[a^-,\i)$. The measure $\nu^{\ell^0,f^0,\L}$ is the law of a 
Markov jump process $x\mapsto \rho^0(x)=g'(x), \ x\in \L$, with 
the rate density
$f^0$, and the initial $\rho(a^0,0)$ distributed according to $\ell^0$. Then 
$\widehat\Phi_t(\rho^0)$ is a Markov jump process associated 
with a rate $\Theta_t(f)$, where $f(x,t,\rho^-,\rho^+)
=\Theta_t(f)(x,\rho^-,\rho^+)$ solves the kinetic equation
\be\la{eq1.19}
f_t-v^Hf_x=Q^H(f),\ \ \ \ \ f(x,0,\rho^-,\rho^+)
=f^0(x,\rho^-,\rho^+),
\ee
where $v^H=v^H(\rho^-,\rho^+)=(H(\rho^-)-H(\rho^+))/(\rho^--\rho^+)$,
and $Q^H$ is as $Q$ of \eqref{eq1.14}, except that $\a$ is replaced with $v^H$. We write $\th^H_t$ for the flow of the kinetic equation
\eqref{eq1.19}:
\[
f(x,t,\rho^-,\rho^+)=:\th^H_t(f^0)(x,\rho^-,\rho^+).
\]

Before giving precise statements for our results in dimension $2$, we need to address a technical issue concerning the domain of the definition
of the function $u$. We remark that because of the quadratic nature
of the right-hand of our kinetic equation \eqref{eq1.12},  generically
 non-negative solutions are defined only locally in spatial variables (in the Appendix, we provide a local well-posedness for \eqref{eq1.12}
under some natural assumptions). Because of this, we will consider 
the PDE \eqref{eq1.12} in a convex domain $\L$ on which our 
Markovian kernels can be defined. In order to solve \eqref{eq1.18} in $\L$, we need to assign suitable boundary conditions. These boundary conditions are selected so that the law of the corresponding solutions are consistent as we vary $\L$.

As we mentioned earlier, since the function $u:\L\x[0,T]\to\bR$
is a piecewise linear convex function, it induces a tessellations
of $\hat\L=\L\x[0,T]$. The vector 
$\hat\rho=(u_x,u_t)=(\rho,H(\rho))$ lies on the graph of $H$, and if a 2-dimensional face $F$ separates 
$\hat\rho^+=(\rho^+,H(\rho^+))$ from 
$\hat\rho^-=(\rho^-,H(\rho^-)), $ then any vector $\hat v=(v,v_3)=
(v_1,v_2,v_3)\in\bR^2\x\bR$
parallel to $F$ must satisfy
\be\la{eq1.20}
v\cdot (\rho^+-\rho^-)+v_3\left(H(\rho^+)-H(\rho^-)\right)=0.
\ee

We address the question of invariance of $\cM_G$ in two settings. 

\bs\noi
{\bf Setting 1.1(i)} {\em (Hamiltonian Function)}
 We  assume that the convex 
function $H(\rho_1,\rho_2)=
H(\rho_1)$, depends on $\rho_1$ only. To simplify our presentation, we also assume that $H$ is an increasing function (see Remark 1.3 for 
general $H$). 

\ms\noi 
{\bf(ii)} {\em (Initial Condition)} Assume that $\L=[a^-,a^+]\x[t_0,t_1]$,
and $\nabla g=\rho(x,0)$ is distributed according to 
$\nu^{\ell^0,f,\L}$,
for a function $f$ that satisfies the kinetic equation \eqref{eq1.12}
in the set $\L$. We additionally assume that the measure 
$f(x,\rho^-,d\rho^+)$ 
is supported in the set $R_0(\rho^-)$, so that Theorem 1.2 is applicable.

\ms\noi 
{\bf(iii)} {\em (Kernel)} We define a kernel $\hat f (x,t,\rho^-,\rho^+)$
by
\[
\hat f(x_1,x_2,t,\rho^-,\rho^+)=
\Theta^H_t(f)(x_1,x_2,\rho^-,\rho^+):=
\th_t^H(h^{x_2})(x_1,\rho^-,\rho^+),
\]
where $h^{x_2}(x_1,\rho^-,\rho^+):=f(x_1,x_2,\rho^-,\rho^+)$.
We put
\[
\left(\hat f^1,\hat f^2,\hat f^3\right):=
\left(\hat f,\a\hat f,v^H\hat f\right),
\]
where 
\[
\a (\rho^-,\rho^+)=[\rho^-,\rho^+],\ \ \ \
v^H(\rho^-,\rho^+)=\left(H(\rho_1^-)-H(\rho_1^+)\right)/
\left(\rho_1^--\rho_1^+\right).
\]

\ms\noi
{\bf(iii)} {\em (Boundary Condition)} We assign a boundary condition at $x=a^+$.
The law of $(x_2,t)\mapsto u_x(a^+,x_2,t)=\rho(a^+,x_2,t)$
is denoted by $\mu$. We assume that under $\mu$, the process
$t\mapsto \rho(a^+,x_2,t)$ is a Markov jump process
with the jump rate density $\hat f^3(a^+,x_2,t)$ for every $x_2\in[t_0,t_1]$.
\qed

\bs
Note that for fixed
$x_2$, we may regard the equation \eqref{eq1.18} as a HJ equation
in dimension one in $x_1$ variable. The process 
\[
(x_1,t)\in [a^-,a^+]\x[0,T]\to
m(x_1,t)=m(x_1,t;x_2):=\rho(x_1,x_2,t),
\]
 is piecewise constant and induces a tessellation
in $[a^-,a^+]\x[0,T]$. 
The process $x_1\mapsto m(x_1,0)$
is a Markov jump process. Its discontinuity points $z_1<\dots<z_n$
travel with time $t$ with velocity
$-v^H(\rho_1^-,\rho^+_1)$.
We are tempted to use either Theorem 1.1 or Theorem 1.2
(or [KR1])
to determine the Markovian law of the process
 $(x_1,t)\mapsto \rho(x_1,x_2,t)$. However these theorems cannot be applied directly because the relation between the particle velocities $-v^H$  and the slopes $\rho^\pm$ is not exactly what we had in these theorems (namely $-[\rho^-,\rho^+]$).
We note that what we have is slightly different from the setting of [KR1], as the jump rate
depends on a vector $\rho$, not just its first coordinate 
$\rho_1=u_{x_1}$. Nonetheless a verbatim 
proof would allow us to have a similar result. In other words, given
a {\em velocity function} $v(\rho^-,\rho^+)$ satisfying some natural
conditions, we may consider a 
particle system as in the Definition 1.3{\bf(vii)} such that 
$[\rho^-,\rho^+]$ is replaced 
with $v(\rho^-,\rho^+)$, and the analogs of Theorems 1.1 and 1.2 
are still true.  
In summary we have the following result:
\bth\la{th1.3} Under the Setting 1.1, the following statements are true:

\ms\noi
{\bf(i)} For every $(x_2,t)\in [t_0,t_1]\x[0,T]$, the process 
$x_1\mapsto \rho(x_1,x_2,t)$ is a Markov jump process
with the jump rate density $\hat f(x_1,x_2,t,\rho^-,\rho^+)$.

\ms\noi
{\bf(ii)} For every $(x_1,x_2)\in \L$, the process 
$t\mapsto \rho(x_1,x_2,t)$ is a Markov jump process
with the jump rate density $\hat f^3(x_1,x_2,t,\rho^-,\rho^+)$.
\et

\ms
Naturally, we may wonder whether or not the law of the process
$(x_1,x_2)\mapsto \rho(x_1,x_2,t)$ is given by $\nu^{\ell^t,\Theta_t^H(f),\L}$, where $\ell^t$ represents the law of 
$\rho(a^-,t_0,t)$. To examine this possibility, let us switch to a more  symmetric notation and 
write $x_3$ for $t$, $\hat x$ for $(x_1,x_2,x_3)$, and $\rho_3$ for
$H(\rho_1)$. In this way,
\[
\hat f^i(\hat x, \rho^-,\rho^+)=[\rho^-,\rho^+]_i 
\hat f(\hat x, \rho^-,\rho^+),
\]
where
\[
[\rho^-,\rho^+]_i=\frac{\rho^-_i-\rho^+_i}{\rho^-_1-\rho^+_1}.
\]
Now if the law of the process
$(x_1,x_2)\mapsto \rho(x_1,x_2,t)$ is given by $\nu^{\ell^t,\Theta_t^H(f),\L}$, then we know that $x_i\mapsto \rho(\hat x)$ 
is a Markov jump process with the jump rate density $\hat f^i$.
Let us write $\cL^i$ for the infinitesimal generator of this jump process.
If $\ell (\hat x,\rho)$ is the law of $\rho(\hat x)$, then it must satisfy the forward equations
\[
\ell_{x_i}=\left(\cL^i\right)^*\ell,\ \ \ \ i\in\{1,2,3\},
\]
where $\left(\cL^i\right)^*$ denotes the adjoint of the operator 
$\cL^i$. From the compatibility of these equations, namely $\ell_{x_ix_j}=
\ell_{x_jx_i}$, we derive the kinetic equation
\be\la{eq1.21}
f^i_{x_j}-f^j_{x_i}=Q(f^i,f^j).
\ee 
In summary, the function $f$ must satisfy the kinetic equation
\eqref{eq1.21} for all $i,j\in\{1,2,3\}$. When these equations hold,
we may choose a Gibbsian measure for the boundary condition
at $x_1=a^+$, and we expect that the law of the marginals
$(x_1,x_2)\to \rho(x_1,x_2,t)$, and $(x_2,t)\to \rho(x_1,x_2,t)$
to be Gibbsian of the type we have constructed. We leave further investigation of the system \eqref{eq1.21} for future.

\ms\noi
{\bf Remark 1.4} Assume that $H(\rho_1,\rho_2)=H_1(\rho_1)+
H_2(\rho_2)$ with $H_1$ and $H_2$ convex and increasing. 
To display the dependence on the Hamiltonian function, we write
$\Phi_t^{H_i}$ for the Hamilton-Jacobi flow $\Phi_t$
(as was defined in Section 1.1) associated
with $H_i$. Writing $g^{x_2}(x_1):=g(x_1,x_2)$, and
$\hat g(x_1,x_2)=\hat g^{t,x_1}(x_2)=\Phi_t^{H_1}(g^{x_2})(x_1)$, then using \eqref{eq1.4}, it is not hard to show that solution $u$ of \eqref{eq1.18}
can be expressed as 
\[
u(x_1,x_2,t)=\Phi_t^{H_2}(\hat g^{t,x_1})(x_2).
\]
Now if $g$ is distributed according to $\nu^{\ell^0,f^0,\L}$, then
we may apply Theorem 1.3 to assert the marginal  $x_1
\mapsto \nabla \hat g(x_1,x_2)$ is a jump process
with the jump rate that are expressed  in terms of $\Theta_t^{H_1}f$. If the law of $\hat g$ is also a Gibbsian measure associated with 
$\Theta_t^{H_1}f$, then another application of Theorem 1.3 would allow us to assert that
the marginals of $\rho=u_x$ are jump processes
with the jump rates that are expressed  in terms of 
$\Theta_t^{H_2}\Theta_t^{H_1}f$.
\qed

So far we have described some of the challenges we encounter as we try to examine the invariance of the set $\cM_G$ under the Setting 1.1. Fortunately these 
challenges can be avoided when we examine the invariance of $\widehat \cM_G$.

\bs\noi
{\bf Setting 1.2(i)} {\em (Hamiltonian Function)}
 We  assume that the convex 
function $H(\rho_1,\rho_2)$ is increasing with respect to both $\rho_1$ and $\rho_2$.

\ms\noi 
{\bf(ii)} {\em (Initial Condition)} Let $\L=[a^-,a^+]\x[t_0,t_1]$,
and let $K:\bR\to\bR$ be an increasing continuous function.
We assume that $\nabla g=\rho(x,0)$ is distributed according to 
$\nu^{\ell^0,\tilde f,\L,K}$ as in Example 1.1. Recall 
\[
\tilde f(x_1,x_2,\rho^-_1,\rho_1^+)=f(x_1,x_2,\rho^-_1,K(\rho_1^-),
\rho_1^+,K(\rho_1^-)),
\]
with $f$ a kernel that satisfies the kinetic equation \eqref{eq1.12} in the set $\L$.

\ms\noi 
{\bf(iii)} {\em (Kernel)} We define a kernel $\hat f (x,t,\rho^-,\rho^+)$
as in the Setting 1.1(iii). We also continue to use our notations 
$\hat f^i$ and $\hat x=(x_1,x_2,x_3)$ as in Setting 1.1. 
With a slight abuse of notation, we define
$\Theta_t(\tilde f)$ as
\[
\Theta_t(\tilde f)(x_1,x_2,\rho^-_1,\rho_1^+)=\hat f(x_1,x_2,t,\rho^-_1,K(\rho_1^-),
\rho_1^+,K(\rho_1^-)).
\]

\ms\noi
{\bf(iv)} {\em (Boundary Condition)} We assign a boundary condition at $x=a^+$ so that the law of $(x_2,t)\mapsto u_x(a^+,x_2,t)=\rho(a^+,x_2,t)$
is again a Gibbsian measure $\mu$ of the type we defined
 in Example 1.1, but now
$(\rho_2,\rho_3)=(u_{x_2},u_{t})$ lies on the graph of 
$\hat K(m)=H(K^{-1}(m),m)$. In particular, $u_t=\hat K(u_{x_2})$ 
in the support of $\mu$. 
\qed

\bs\noi
{\bf Conjecture 1.1} Under the Setting 1.2, the law of 
$(x_1,x_2)\mapsto \rho(x_1,x_2,t)$ is given by the measure
$\nu^{\ell^0,\Theta_t(\tilde f),\L,K}$.
\qed

\bs
We have been able to partially verify this conjecture:

\bth\la{th1.4} Under the Setting 1.2, the process  
$x_i\mapsto \rho(\hat x),\ i\in\{1,2,3\},$ is a Markov jump process with the jump rate
density $\hat f^i$.
\et

\ms
In fact we can readily establish Theorem 1.4 with the aid of Theorem 1.3.
We explain this in three short steps:

\ms\noi
{\em (Step 1)} We first argue that the relationship
$u_{x_2}=K(u_{x_1})$ that is assumed at $t=0$, also holds at  later times. To see this, observe that 
if initially
\[
g(x_1,x_2)=\sup_{\rho_1\in R}(x_1\rho+x_2 K(\rho_1)-\a(\rho_1)),
\]
for a (discrete) set $R$ and a function 
$\a(m)=g^*(m,K(m))$, then on account of \eqref{eq1.3}, a similar formula is true at a later time $t$, where $\a(\rho_1)$ is replaced with
$\a(\rho_1)-t\tilde H(\rho_1)$, for
$\tilde H(\rho_1)=H(\rho_1,K(\rho_1))$.
 We can take advantage of this property to reduce the question of invariance
to  Case 1. After all  if  $u_{x_2}=K(u_{x_1})$ holds for
a solution $u$ of \eqref{eq1.18}, then such a solution also solves
the equation $u_t=\tilde H(u_{x_1})$.

\ms\noi
{\em (Step 2)}
We note that $\tilde H$ is convex if $H$ is convex. Let us present
a short proof of this when $H$ is $C^2$ and $K$ is differentiable:
\begin{align*}
\tilde H''(m)=&
H_{\rho_1\rho_1}(m,K(m))+2H_{\rho_1\rho_2}(m,K(m))K'(m)
+H_{\rho_2\rho_2}(m,K(m))K'(m)^2\\
=&(D^2H)(m,K(m))\begin{bmatrix}1\\K'(m)\end{bmatrix}
\cdot \begin{bmatrix}1\\K'(m)\end{bmatrix}\ge 0.
\end{align*}
Furthermore, if $K$ is increasing, and $H$ is increasing in both arguments, then $\tilde H$ is also increasing.
This allows us to apply Theorem 1.3 to assert that the process
$x_i\mapsto \rho(\hat x)$ is a Markov jump process with the jump
rate density $\hat f^i$, for $i=1$ and $i=3$.

\ms\noi
{\em (Step 3)} Since $K$ is increasing we can interchange the role of
$x_1$ with $x_2$. That is, from $u_{x_1}=K^{-1}(u_{x_2})$,
we learn that $u_t=\bar H(u_{x_2})$, where 
$\bar H (\rho_2)=H(K^{-1}(\rho_2),\rho_2)$. Since $\bar H$ is convex and increasing, and the boundary dynamics at $x_2=t_1$ is Markovian
with the jump rate density $\hat f^3$, we are at a position to apply
Theorem 1.3 once more to assert that the process $x_2\mapsto
\rho(\hat x)$, is Markov jump process
with the jump rate density $\hat f^2$.

\bs\noi
{\bf Remark 1.5} When $K$ is also concave or convex, then Conjecture 1.1 would follow from Theorem 1.4. For example, if $K$ is concave, then 
we can fully determine the function $(x_1,x_2)\to \rho(x_1,x_2,t)$ 
from its boundary values on $\p \L$. Simply because as we decrease $x_2$, there would be no (stochastic) fragmentation and the corresponding particle system involves free motion and collisions.
\qed

\subsection{Generalization to higher dimensions $d>2$.}

As we mentioned earlier, a Gibbsian measure $\nu^{f,\L}$ for the distribution of $\nabla g$ in \eqref{eq1.18} would lead to a 
Gibbsian measure $\hat \nu^{\hat f,\hat\L}$ for the distribution
of $\hat\rho(x,t)=(u_x,u_t)(x,t),$ where $\hat\L=\L\x [0,T]$, and
\[
\hat f(x,t,\rho^-,\rho^+)=:g\left(x,t,\rho^-,H(\rho^-), \rho^+,H(\rho^+)\right),
\]
represents a density rate at which  $\hat\rho^-=(\rho^-,H(\rho^-))$
changes to $\hat\rho^+=(\rho^+,H(\rho^+))$ at the point
$(x,t)\in\bR^3$. Indeed the function $u(\hat x)=u(x_1,x_2,x_3)
=u(x_1,x_2,t)$ is a piecewise linear convex function such that its gradient
$\hat\rho$ lies on the graph of $H$. The 
piecewise constant function $\hat\rho$ yields a Laguerre tessellation
\[
\widehat\bX=\left\{(\hat\rho,X(\hat\rho)):\ \hat\rho\in \widehat\cS\right\},
\]
of $\bR^3$. We may wonder whether or not we can apply our approach 
to build more general Gibbsian measure on the set of Laguerre tessellations of $\bR^3$. More specifically we wish to relax the restriction
$\hat\rho\in\{(a,H(a)):\ a\in\bR^2\}$. To describe a strategy for
achieving this, let us first discuss some of the features
of Laguerre tessellations in $\bR^3$. Generically the following statements are true:
\bi
\item If $\hat\rho^+\neq\hat\rho^-$, and 
$\widehat{X}(\hat\rho^-,\hat\rho^- ):=\widehat{X}(\hat\rho^-)
\cap \widehat{X}(\hat\rho^+)\neq\emptyset$, then 
$\widehat{X}(\hat\rho^-,\hat\rho^+ )$ is a 
($2$-dimensional) convex polygon.
The vector $\hat\rho^+-\hat\rho^-$
points from the
$\hat\rho^-$ side to the $\hat\rho^-$ side of 
$\widehat{X}(\hat\rho^-,\hat\rho^+ )$. 
\item If $\hat\rho^-,\hat\rho^+$, and $\hat\rho^*$ are distinct,
and  $\widehat{X}(\hat\rho^-,\hat\rho^+,\hat\rho^* ):=
\widehat{X}(\hat\rho^-)\cap \widehat{X}(\hat\rho^+)
\cap \widehat{X}(\hat\rho^*)\neq\emptyset$, then
$\widehat{X}(\hat\rho^-,\hat\rho^+,\hat\rho^* )$ is an edge
(a line segment).
We can uniquely determine a vector direction 
\[
\tau=\tau(\hat\rho^-,\hat\rho^+,\hat\rho^*)=
(v(\hat\rho^-,\hat\rho^+,\hat\rho^*),1),\ \ \ \ v\in\bR^2,
\]
of this edge by solving the system of linear equations 
$\tau\cdot (\hat\rho^\pm-\rho^*)=0.$
\item If $\hat\rho^i,\ i=1,\dots,4$ are distinct, and $
\widehat X(\hat\rho^1,\dots,\hat\rho^4):= \widehat X(\hat\rho^1)
\cap\dots\cap \widehat X(\hat\rho^4)\neq\emptyset,$
then $\widehat X(\hat\rho^1,\dots,\hat\rho^4)$ consists of a single point which is a vertex of our tessellation.
\ei

\ms
We wish to build a random tessellation $\widehat\bX$ so that 
its intersection with the plane $\{\hat x:\ x_3=t\}$ is a random planar
tessellation of the type we have constructed in the Section 1.3.
Let us ignore the lateral boundary dynamics and focus on our strategy for building such a tessellation inside the box $\hat \L$. In other words,
we start from a planar tessellation that represents the restriction
of $\widehat\bX$ to the plan $\{x_3=0\}$, and evolve it in a Markovian fashion as we increase $x_3$. The law of this planar tessellation is a suitable $\nu^{f^0,\L}$, except that the kernel
$f^0(x_1,x_2,\hat\rho^-,\hat\rho^+)$ is defined for $\hat\rho^\pm
\in\bR^3$. In other words, the cells of the initial tessellations
labeled/decorated by vectors in $\bR^3$. To build our 
random tessellation initially, we need to assume that the kernel $f^0$
satisfies the kinetic equation \eqref{eq1.12}, where only the first
two coordinates of $\hat\rho$ are used for determining the speed
$\a$. The evolution of our tessellation consists of the deterministic and the stochastic parts. As for the deterministic part, 
 a vertex associated with vector densities $\hat\rho^-,\hat\rho^+$, and $\hat\rho^*$ travels with the velocity 
$v(\hat\rho^-,\hat\rho^+, \hat\rho^*)$.  As $x_3$ increases, it is possible that an edge of vertices $a^-$ and $a^+$ collapses, or equivalently $a^-$ and $a^+$ collide. 
There would be two possibilities for the type of a collision that can occur. 
To explain this, let us write $C^-$ and $C^+$ for the cells which 
are sharing the edge $a^-a^+$. 
\bi
\item One of the cells $C^\pm$ is a triangle (generically not both $C^-$ and $C^+$ could be triangle). When this is the case, the whole triangular
cell collapses and becomes a vertex. The tessellation has lost 
a cell at time of such a collision.
\item Neither $C^-$ nor $C^+$ are triangle. If $a^-$ is a vertex associated with cells of labels $(\hat m^-,\hat\rho^-,\hat\rho^+)$, 
and $a^+$  is a vertex associated with cells of labels $(\hat m^+,\hat\rho^-,\hat\rho^+)$, then after the collision a new edge is created
with vertices $b^-$ and $b^+$. The new vertices $b^-$ and $b^+$
are now associated with
cells of marks $(\hat \rho^-,\hat m^-,\hat m^+)$, and $(\hat \rho^+,\hat m^-,\hat m^+)$, respectively (the role of $\hat m$ and $\hat \rho$ are swapped).
\ei

Our dynamics also involves a stochastic fragmentation; at a random time
$t$, a vertex $a$, associated with $(\hat\rho^1,\hat\rho^2,\hat\rho^3)$, 
can give birth to a triangle with vertices $b_1,b_2$ and $b_3$, and a
random label $\hat\rho^*$.
These vertices start 
their journey at the location $a$ at time $t$, and move away from each
others with velocities that are determined in terms of 
$(\hat\rho^1,\hat\rho^2,\hat\rho^3)$, and  $\hat\rho^*$.

In order to carry out our program in dimension $3$, we need to work out
the form of the fragmentation rate. We emphasis that when 
$\hat \rho$ lies on a graph of a convex function $H$, there would be no fragmentation. We conjecture that when the rate $\hat f$ satisfies 
a system of kinetic equations analogous to \eqref{eq1.21}, our outlined
strategy would yield a consistent family of Gibbsian measures on 
the set of tessellations of $\bR^3$.

\subsection{Bibliography and the outline of the paper}

Most of the earlier works on stochastic solutions of Hamilton-Jacobi PDEs have been carried out in the Burgers context. For example,
 Groeneboom [Gr] determined the statistics of
  solutions to Burgers equation ($H(p)=p^2/2$, $d=1$)
with white noise initial data.
  Burgers equation is not explicitly mentioned---the paper 
  discusses convex minorants of Brownian motion with parabolic
  drift---but these problems are connected by the Hopf-Lax-Oleinik
  solution formula \eqref{eq1.4}. 
Recently Ouaki [O] has extended this result to arbitrary convex Hamiltonian function $H$. The special cases of 
$H(p)=\i \1(p\notin[-1,1])$, and $H(p)=p^+$ were already studied in references Abramson-Evans [AE], Evans-Ouaki [EO], and Pitman-Tang [PW].

 Carraro and Duchon [CD1-2] considered
  \emph{statistical} solutions, which need not coincide with genuine
  (entropy) solutions, but realized in this context that L\'evy
  process initial data  should
  interact nicely with Burgers equation.  Bertoin [Be]
  showed this intuition was correct on the level of entropy solutions,
  arguing in a Lagrangian style.

Developing an alternative treatment to that given by Bertoin, which
relies less on particulars of Burgers equation and happens to be more
Eulerian, was among the goals of  Menon and Srinivasan [MS].
Most notably, [MS] formulates an interesting conjecture for the evolution
of the infinitesimal generator of the solution $\rho(\cdot,t)$ which is
equivalent with our kinetic equation \eqref{eq1.14} when there is no fragmentation, $f$ is independent of $x$, and
$\rho$ lies on a graph of a convex function (see Example 1.1{\bf(i)}).
When the initial data
$\rho(x,0)$ is allowed to assume values only in a fixed, finite set of
states, the infinitesimal generators of the processes $x\mapsto \rho(x,t)$ and $t\mapsto \rho(x,t)$ 
can be represented by triangular matrices.  
The integrability of this matrix evolution
has been investigated by Menon [M2] and Li [Li].  For generic matrices---where the genericity
assumptions unfortunately exclude the triangular case---this evolution
is completely integrable in the Liouville sense.  
The full treatment of Menon and Srinivasan's conjecture was achieved
in papers [KR1] and [KR2] (we also refer to [R1] for an overview).
The works of [KR1-2] have been recently extended in [R2] to allow 
inhomogeneous HJ equation of the type \eqref{eq1.1} in dimension one.

\bs\noi
The organization of the paper is as follows:
In Section 2 we give a precise construction of the particle system
$\bq(t)$ that we described in Definition 1.3{\bf(vii)}.
In Section 3, we derive a {\em forward equation} for the law of
$\bq(t)$ and  establish Theorem 1.1.  Section 4 is devoted to the proof of Theorem 1.2. In the Appendix we address the
question of well-posedness and the regularity of solutions of the kinetic equation.

\section{Construction of the Particle System}

As discussed in the introduction of the paper, we use $z$ for the variable $x_1$ and time $t$ for the variable $x_2$. 
We wish to build a probability measure on the space of 
\textit{Laguerre tessellations} with orientation $\tau : \{ (\rho^{-},\rho^{+}) \in (\bR^2)^2: \rho^{-} \prec \rho^{+} \} \to \bR^2$ that was given by \eqref{eq1.7}.
We have already given a rough description for this probability measure
in Section~1.3. In this section we give more thorough details and 
make some rudimentary preparations for the proof of our main results. Our strategy will be to build a consistent family of probability measures on the space of labeled interacting particle systems on fixed boxes.

Fix $t_0,t_1,a^-,a^+\in \bR$, with $a^-<a^+$ and $t_0<t_1$.
We will build a probability measure on the space of stochastic processes $(\bq(t))_{t \in[ t_0,t_1]}$ that take values in the space of particle systems of the form
\[
\bq(t):=\left((z_0,\rho^0(t)),(z_1(t),\rho^1(t)),\cdots,(z_{n(t)}(t),\rho^{n(t)}(t))\right),
\]
with $z_0=a^-$. To be more precise, let us introduce some notation. 

\bs\noi
{\bf Notation 2.1(i)}
Let $P^-<P^+$ be two fixed constants and define the state space $\Omega=\O_{a^-,a^+}$ of particle systems as the following disjoint union
\[
\Omega=\bigsqcup_{n=0}^{\i} \Omega^n
\]
where $\Omega^n := \overline{\Delta^n} \times R^n$, with
 $\Delta^n:=\big\{(z_1,\cdots,z_n) : a^-<z_1<\cdots<z_n<a^+\big\}$,
$\overline{\Delta^n}$ is the topological closure of 
$\Delta^n$ in $\bR^n$, and 
\[
R^n:=\left\{(\rho^0,\rho^1,\cdots,\rho^n) \in ([P^-,P^+]^2)^{n+1}:  \rho^0 \prec \rho^1 \prec \cdots \prec \rho^n\right\} .
\]

\ms\noi
{\bf (ii)}
For any $\bq \in \Omega$, let $\bn(\bq)$ be the unique integer $n$ such that $\bq \in \Omega^n$.

\ms\noi
{\bf(iii)}
  For any three labels $\rho,\rho',\rho'' \in [P^-,P^+]^2$ such that $\rho \prec \rho' \prec \rho''$, set
\[
[\rho,\rho',\rho'']=[\rho',\rho'']-[\rho,\rho']
\]

\ms\noi
{\bf(iv)}
 For any real number $r$, let $r^{+}=\max(r,0)$ and $r^{-}=\max(-r,0)$. For any two integers $m \le n$, we denote by $\llbracket m,n \rrbracket$ the set $\{m,m+1,\cdots,n\}$.

\ms\noi
{\bf(v)}
Let $\Gamma$ be the space of vector-valued right-continuous piecewise-constant functions on $[a^-,a^+]$, an define the map $\cV : \Omega \to \Gamma$ as follows.
For any $\bq \in \Omega$ of the form
\[
\bq=\left((a^-,\rho^0),(z_1,\rho^1),\cdots,(z_n,\rho^n)\right) \in \Omega^{n}
\]
Define
\be\la{eq2.1}
\cV(\bq)(z):=\sum_{i=0}^{n} \rho^i \ \ind(z_i \le z < z_{i+1})\ \
 \text{ for all } z \in [a^-,a^+],
\ee
with the convention that $z_0=a^-$ and $z_{n+1}=a^+$.

\ms\noi
{\bf(vi)}
Let $\Omega'$ to be the subset of $\Omega$ such that for any $\bq=((a^-,\rho^0),(z_1,\rho^1),\cdots,(z_n,\rho^n)) \in \Omega' \cap \Omega^{n}$, and if $z_i=z_{i+1}$ for some $i \in \llbracket 1,n \rrbracket$, then we have
\[
[\rho^{i-1},\rho^i,\rho^{i+1}] \ge 0.
\]
\qed

\bs

To construct our measure $\nu^{\ell,f,\L}$, 
we take a kernel $f\in \cF(\b,\L)$, 
\[
f(z,t,\rho^-,d\rho^+) =f(z,t,\rho^-,\rho^+)\ \b(d\rho^+),
\]
that satisfies the  kinetic equation \eqref{eq1.12} in $\L$. 
We then use $f$ to define a Markov process $\big(\bq(t):\ t\in[t_0,t_1]\big)$
in $\O$. The law of this Markov process for a suitable initial distribution $\ell$ is the desired measure $\nu^{\ell,f,\L}$.
For the reader's convenience, we recall the kinetic equation,
\be\la{kin-f}
\tau(\rho^-,\rho^+)\cdot\nabla f=Q(f)(z,t,\rho^{-},\rho^{+}),
\ee
 where $Q(f)=Q^+(f)-f\ Lf$,  with  
\begin{align*}
Q^{+}(f)(z,t,\rho^{-},\rho^{+})=&\int_{D(\rho^{-},\rho^{+})}  [\rho^{-},\rho^{*},\rho^{+}]f(z,t,\rho^{-},\rho^{*})f(z,t,\rho^{*},\rho^{+})\ \b(d\rho^{*}),\\
(Lf)(z,t,\rho^{-},\rho^{+})=&A(z,t,\rho^{+})-A(z,t,\rho^{-})
-[\rho^{-},\rho^{+}]\left(\lambda(z,t,\rho^{+})-\lambda(z,t,\rho^{-})
\right).
\end{align*}
Here, 
\begin{align*}
\lambda(z,t,\rho)=&\int_{R(\rho)} f(z,t,\rho,\rho^{*})\ 
\b(d\rho^{*}),\\
A(z,t,\rho)=&\int_{R(\rho)} [\rho,\rho^{*}]f(z,t,\rho,\rho^{*})\ \b(d\rho^{*}).
\end{align*}
Now, consider a nonnegative function
 $(z,t,\rho) \mapsto \ell(z,t,\rho)$,  verifying the following equations
\begin{align}\label{kin-x}
\ell_{z}(z,t,\rho)=&\int_{L(\rho)} f(z,t,\rho^{*},\rho)\ell(z,t,\rho^{*})\ \b(d\rho^{*})-\lambda(z,t,\rho)\ell(z,t,\rho),\\
\label{kin-t}
\ell_{t}(z,t,\rho)=&\int_{L(\rho)}[\rho^{*},\rho] f(z,t,\rho^{*},\rho)\ell(z,t,\rho^{*})\ \b(d\rho^{*})-A(z,t,\rho)\ell(z,t,\rho).
\end{align}
Moreover for any fixed $z$ and $t$ we assume
$
\int \ell(z,t,\rho)\ \b(d\rho)=1.
$
As we will see in Proposition 4.1 below,
the two equations \eqref{kin-x} and \eqref{kin-t} are compatible because of the kinetic equation verified by $f$. Mainly, the two flows generated in both direction $z$ and $t$ commute, as the corresponding vector fields commute in the sense of having a zero Lie Bracket. This last fact is exactly a reformulation of the equation \eqref{kin-f}.

Without loss of generality let us assume that $t_0=0$ to alleviate the notation. The purpose of this section is to construct the law of a stochastic process $(\bq(t))_{t \ge 0}$ that takes values in $\Omega$. 
Let us introduce first the following probability measure on 
$\Omega$. For any $t \ge 0$, let 
\be\la{eq2.5}
\mu(d\bq,t)=\sum_{n=0}^{\i} \ind(\bn(\bq)=n)\
 \mu^{n}(d\bq,t),
\ee
where $\mu^{n}(d\bq,t)=g^n(\bz,\brho)\ d\bz\ \b(d\brho)$ is a measure defined on $\Omega^n$, with
\begin{align*}
&\b(d\brho):=\prod_{i=0}^n\b(d\rho^i),\ \ \ \ 
d\bz=\prod_{i=0}^ndz_i,\\
&g^n(\bz,\brho):=\ell(a^-,t,\rho^0)\ \prod_{i=1}^{n} f(z_i,t,\rho^{i-1},\rho^i)\ \exp \left(-\sum_{i=0}^{n} \int_{z_i}^{z_{i+1}}
 \lambda(z,t,\rho_i)\ dz \right).
\end{align*}
Here $\bq=((z_0,\rho^0),(z_1,\rho^1),\cdots,(z_n,\rho^n))$,
and by convention $z_0=a^-$ and $z_{n+1}=a^+$. 

\subsection{The deterministic flow}
We first start by defining a deterministic flow on $\Omega$. Let $t \ge 0$, and $\bq \in \Omega^n$. We distinguish two cases:

\ms\noi
{\bf (1)}
If $\bq \in \text{Int}(\Omega^n)$, i.e.,
 $\bq=((a^-,\rho^0),(z_1,\rho^1),\cdots,(z_n,\rho^n))$ where 
$a^-<z_1<\cdots<z_n<a^+$. 
Let $(z_{i}(t))_{i=1}^{n}$ be defined as 
\[
\dot{z}_i(t)=-[\rho^{i-1},\rho^i],~~ \ \ z_i(0)=z_i,
\]
for $t \in [0,T^{*})$, where
\[
T^{*}=\inf\left \{t>0: z_i(t)=z_{i+1}(t) \text{ for some } i \in \{0,1,\cdots,n\} \right\},
\]
the time of the first collision. By convention, we always take $z_0(t)=a^-$ and $z_{n+1}(t)=a^+$ when we consider a particle system of $n$ particles. This defines a flow $\psi^{t}\bq$: 
\[
\psi^{t}\bq=\left((a^-,\rho^0),(z_1(t),\rho^1),\cdots,(z_n(t),\rho^n)\right) \in \Omega^n\ \text{ for }\ t \in [0,T^{*}).
\]
At time $t=T^*$, for any $i$ such that $z_i(T^{*}-)=z_{i+1}(T^{*}-)$, we remove the $(i+1)$-th particle $(z_{i+1}(T^{*}-),\rho^{i+1})$ from the particle system and relabel the particle $(z_{j}(T^{*}-),\rho^j)$ as $(z_{j-1}(T^{*}-),\rho^{j-1})$, for $j>i+1$ . We keep doing this procedure until we end up with an element in $\text{Int}(\Omega^{m})$ for some $m \le n-1$ (it might be possible that multiple collisions happen at the same time $T^{*}$), we then define $\psi^{T^{*}}\bq$ to be this element. Repeating the same process again starting from this new configuration, we get a sequence of collision times $0<T_1^{*}<T_2^{*}<\cdots$, where $t \mapsto \psi^{t} \bq$ evolves with free motion in each interval $[T_{k}^{*},T_{k+1}^{*})$. 

\ms
\noi
{\bf (2)} If $\bq \in \partial \Omega^{n}$, then for any $i$ such that $z_i=z_{i+1}$, there are two cases to consider:
\begin{itemize}
\item{ If $[\rho^{i+1},\rho^i] \ge [\rho^i,\rho^{i-1}]$, then we delete the particle $(z_{i+1},\rho^{i+1})$ and relabel the particle $(z_j,\rho^j)$ as $(z_{j-1},\rho^{j-1})$, for $j>i+1$.}
\item{If $[\rho^{i+1},\rho^i] < [\rho^i,\rho^{i-1}]$, we keep both (sticky) particles at $z=z_i=z_{i+1}$.}
\end{itemize}
Denote 
\[
\bar{\bq}:=\left((a^-,\bar\rho^0),(y_1,\bar\rho^1),\cdots,(y_m,\bar\rho^m)\right) \in \Omega^m
\]
to be the resulting configuration after doing this modification. Notice that $\bar{\bq}$ is not necessarily in $\text{Int}(\Omega^m)$, as we may still have an index $j$ such that $y_j=y_{j+1}$ and $[\bar\rho^{j-1},\bar\rho^j,\bar\rho^{j+1}]<0$. The flow $\psi^t \bq:=\psi^t \bar{\bq}$ for any $t>0$ is now defined in the same fashion as before, by letting each particle
 $y_j(t)$ have a free motion with the corresponding velocity 
$v_j:=-[\bar\rho^{j-1},\bar\rho^j]$. Instantaneously, for any small $t>0$,
the configuration  $\psi^t\bq$ belongs to $\text{Int}(\Omega^m)$. To see this, observe that for any $j$ such that $y_j=y_{j+1}$ we have
\[
\frac{d}{dt}(y_{j+1}(t)-y_{j}(t))=-[\bar\rho_{j-1},\bar\rho_j,\bar\rho_{j+1}]>0.
\]
The motion of the particle system then encounters collisions and behaves similarly to the first case when we start from a particle system in the interior of the state space.
\qed

\ms
More generally for any $s \le t$, we define the deterministic flow between time $s$ and time $t$ to be  $\psi_s^t \bq:=\psi^{t-s} \bq$ for any $\bq \in \Omega$. By construction, this flow verifies the semi-group property
\[
\psi_{t_1}^{t_3}\bq=\psi_{t_2}^{t_3}\psi_{t_1}^{t_2}\bq\ \text{ for any } \ t_1 \le t_2 \le t_3,\ \text{ and }\ \bq \in \Omega.
\]

\subsection{The stochastic flow and Markov process} {We will define a stochastic process $(\bq(t))_{t \ge 0}$, that takes values in $\Omega$, and such that $(t,\bq(t))_{t \ge 0}$ is strong Markov. Equivalently, this amounts to constructing a probability measure $\bP^{\bq}_{t}$ for any $\bq \in \Omega$ and $t \ge 0$. This probability measure should be understood as the law of $(\bq(\th))_{\th \ge t}$ conditionally on
$\bq(t)=\bq$. Naturally, the measure $\bP^{\bq}_t$ is concentrated on the set measurable maps $\bq:[t,+\i)\to\Omega$ such that
$\bq(t)=\bq$. Let us start with some notation. 

\bs\noi
{\bf Notation 2.2(i)}
For any $(\rho^{-},\rho^{+})\in [P^-,P^+]^2$ with $\rho^{-} \prec \rho^{+}$ and $t \ge 0$, we define the two quantities
\begin{align*}
\mf{c}_{-}(t,\rho^{-},\rho^{+}):=&
[\rho^{-},\rho^+]^{-}\ 
\frac{\ell(a^-,t,\rho^{-})f(a^-,t,\rho^{-},\rho^+)}{\ell(a^-,t,\rho^+)},\\
\mf{c}_{+}(t,\rho^{-},\rho^{+}):=&[\rho^{-},\rho^{+}]^{+}f(a^+,t,\rho^{-},\rho^{+}).
\end{align*}
They correspond to the rates of creation of particles respectively at $z=a^-$ and $z=a^+$.
For $\rho,\rho',\rho'' \in [P^-,P^+]^2$ with $\rho \prec \rho' \prec \rho''$,  $z \in [a^-,a^+]$ and $t \ge 0$, define the fragmentation rate at position $z$ and time $t$ as
\[
\mf{f}(z,t,\rho,\rho',\rho''):=[\rho,\rho',\rho'']^{-}\ 
\frac{\ f(z,t,\rho,\rho')f(z,t,\rho',\rho'')}{f(z,t,\rho,\rho'')}.
\]

\ms\noi
{\bf(ii)} For any $\rho \in [P^-,P^+]^2$ and $t \ge 0$, let
\[
\mf{C}_{-}(t,\rho):=\int_{L(\rho)} \mf{c}_{-}(t,\rho^{-},\rho)\ \b(d\rho^{-})\  \text{ and }\  \mf{C}_{+}(t,\rho):=\int_{R(\rho)} \mf{c}_{+}(t,\rho,\rho^{+})\ \b(d\rho^{+}),
\]
and for any $\rho^{-},\rho^{+} \in [P^-,P^+]^2$ such that $\rho^{-} \prec \rho^{+}$, and $t \ge 0$, let
\[
\mf{F}(z,t,\rho^{-},\rho^{+}):=\int_{D(\rho^{-},\rho^{+})} \mf{f}(z,t,\rho^{-},\rho^{*},\rho^{+})\ \b(d\rho^{*}).
\]
(See Definition 1.3{\bf(i)} for the definition of $R(\rho^-),\ L(\rho^+)$, and $D(\rho^{-},\rho^{+})$.)
For any particle system $\bq \in \Omega^n$ of the form $\bq=((a^-,\rho^0),(z_1,\rho^1),\cdots,(z_n,\rho^n))$, we define its \textit{particle rate} at time $t$ by
\[
\mf{r}(t,\bq):=\mf{C}_{-}(t,\rho^0)+\sum_{i=1}^{n} \mf{F}(z_i,t,\rho^{i-1},\rho^i) +\mf{C}_+(t,\rho^n).
\]

\ms\noi{\bf(iii)}
We introduce now a notation that corresponds to the state of the particle system after a creation or fragmentation. For any $\rho^{*} \prec \rho^0$, we define the new particle configuration 
$E_{-}^{\rho^{*}}\bq$ as
\[
E_{-}^{\rho^{*}}\bq:=\left((a^-,\rho^{*}),(a^-,\rho^0),(z_1,\rho^1),\cdots,(z_n,\rho^n)\right),
\]
where a new particle is added at $z=a^-$. Similarly for the barrier $z=a^+$, for any $\rho^*$ such that $\rho^{n} \prec \rho^{*}$, define the new particle configuration $E_{+}^{\rho^{*}}\bq$ by
\[
E_{+}^{\rho^{*}}\bq:=\left((a^-,\rho^0),(z_1,\rho^1),\cdots,(z_n,\rho^n),(a^+,\rho^{*})\right).
\]
Finally for $\rho^{*}$ such that $\rho^{i-1} \prec \rho^{*} \prec \rho^i$ for some $i \in \{1,\cdots,n\}$, let
\[
E_{i}^{\rho^{*}} \bq:=\left((a^-,\rho^0),(z_1,\rho^1),\cdots,
(z_{i-1},\rho^{i-1}),(z_i,\rho^{*}),(z_i,\rho^{i}),(z_{i+1},\rho^{i+1}),\cdots,(z_n,\rho^n)\right),
\]
denote the particle configuration we obtain after the fragmentation of the $i$-th particle.
\qed
\ms

\subsubsection{Construction of the Markov process} Consider now a probability space $(\Omega_0,\cF,\bP)$ on which is already defined an infinite i.i.d sequence $(\tau_i)_{i \ge 1}$ of standard exponential random variables. We will define a process $t\mapsto \bq(t)=\bq(t,\o)\in\O$ on this probability space with $\bq(0)=\bq \in \Omega$ and $\o\in\O_0$. The reader should keep in mind that our construction only works up to a time $T^{*}$ before the solution of the kinetic equation cease to be positive (see Appendix A), but to ease the notation we will assume that this state to be infinite. Thus, any temporal quantity $t$ in the future should be thought of implicitly as $t \wedge T^{*}$. Define now the stopping time $T_1$ as
\[
T_1=\inf \left\{ t \ge 0 : \int_{0}^{t} \mf{r}(\th,\psi_{0}^{\th}\bq)\ d\th \ge \tau_1 \right\}.
\]
For any $t \in [0,T_1)$, put $\bq(t)=\psi_{0}^{t}\bq$. Conditionally on $T_1$, write
\[
\psi_{0}^{T_1} \bq=\left((a^-,\rho^0),(z_1,\rho^1),\cdots,(z_n,\rho^n)\right).
\]
\begin{itemize}
\item{With probability $\displaystyle \frac{\mf{C}_{-}(T_1,\rho^0)}{\mf{r}(T_1,\psi_{0}^{T_1}\bq)}$, we sample $\rho^{*}$ with density $\displaystyle \frac{\mf{c}_{-}(T_1,\rho^{*},\rho^0)}{\mf{C}_{-}(T_1,\rho^0)}$ (with respect to the measure $\b(d\rho^{*})$) and put 
$\bq(T_1)=E_{-}^{\rho^*}\psi_{0}^{T_1}\bq$.}
\item{With probability $\displaystyle \frac{\mf{F}(z_i,T_1,\rho^{i-1},\rho^i)}{\mf{r}(T_1,\psi_{0}^{T_1}\bq)}$ for $i \in \{1,2,\cdots,n\}$, sample $\rho^{*}$ with density\\
 $\displaystyle \frac{\mf{f}(z_i,T_1,\rho^{i-1},\rho^{*},\rho^i)}{\mf{F}(T_1,\rho^{i-1},\rho^i)}$ and put $\bq(T_1)=E_{i}^{\rho^{*}}\psi_{0}^{T_1}\bq$.}
\item{With probability $\displaystyle \frac{\mf{C}_+(T_1,\rho^n)}{\mf{r}(T_1,\psi_{0}^{T_1}\bq)}$, sample $\rho^{*}$ with density $\displaystyle \frac{\mf{c}_+(T_1,\rho^{n},\rho^{*})}{\mf{C}_+(T_1,\rho^n)}$ and put $\bq(T_1)=E_{+}^{\rho^*}\psi_{0}^{T_1}\bq$.}
\end{itemize}
We repeat this process again by defining 
\[
T_2:=\inf\left \{t \ge T_1 : \int_{T_1}^{t} \mf{r}(\th,\psi_{T_1}^{\th}\bq(T_1))\ d\th \ge \tau_2\right\},
\]
putting $\bq(t)=\psi_{T_1}^{t}\bq(T_1)$ for $t\in [T_1,T_2)$ and resampling again a $\rho^{*}$ with the analogous above probabilities to define $\bq(T_2)$. This constructs a sequence of random times $T_1,T_2,\cdots$.
To ensure that our process $\bq(t)$ is defined for any time $t \in[0,T]$ and that only finitely many jumps happens, we would assume
\begin{equation}\label{bound-rates}
M:=\sup_{t \in[0,T]}\ \sup_{z \in [a^-,a^+]}\ \sup_{\rho^{-} \prec \rho^{+}} \max\left(\mf{C}_{-}(t,\rho^{-}),
\mf{F}(z,t,\rho^{-},\rho^{+}),
\mf{C}_{+}(t,\rho^{+})\right)<\infty.
\end{equation}
To this end, let us define $N^T(\bq)$ to be 
the number of stochastic jumps up to time $T>0$ when
we starts from $\bq$, i.e.,
\[
N^T(\bq)= \sup \{ n \ge 0 : T_n < T \}.
\]
The following lemma will be used in several occasions:
\begin{lemma}\label{numb-stochas} Assume that
\eqref{bound-rates} holds. Then 
there exists a constant $C_0=C_0(T,M)>0$, independent of $\bq$  such that 
\be\la{lem21} 
\bE[N^T(\bq)] \le C_0(\bn(\bq)+2)^2,\ \ \ \ \ 
\bP(T_{n}<T)\le C_0(\bn(\bq)+2)^2n^{-2}.
\ee 
In particular $N^T(\bq)$ is almost surely finite.
\end{lemma}

\ms\noi
{\bf Proof} For any $k \ge 0$ we have that 
\begin{align}\label{inequ}
\int_{T_k}^{T_{k+1}} \mf{r}(\th,\psi_{T_k}^{\th}\bq(T_k))\ d\th \ge \tau_{k+1},
\end{align}
with the convention that $T_0=0$, and where $(\tau_i)_{i \ge 1}$ is an i.i.d sequence of standard exponential random variables. The number of particles of $\bq(T_k)$ is at most $\bn(\bq)+k$, therefore
 for $\th \in [T_k,T_{k+1})$, we have that $\mf{r}(\th,\psi_{T_k}^{\th}\bq(T_k))\le M(\bn(\bq)+k+2)$. From this and
\eqref{inequ}, we deduce that
\[
T_{k+1}-T_k \ge \frac{\tau_{k+1}}{M(\bn(\bq)+k+2)}.
\]
From this we learn,
\[
\bP(N^T(\bq) > n)=\bP(T_{n+1}<T)\le \bP\left(\sum_{k=1}^{n} \frac{\tau_{k+1}}{\bn(\bq)+k+2} \le MT \right).
\]
Hence for any $\lambda>0$, by Markov inequality
\begin{align*}
\bP\left(\sum_{k=1}^{n} \frac{\tau_{k+1}}{\bn(\bq)+k+2} \le MT \right)&=\bP\left( \exp \left(-\lambda \sum_{k=1}^{n} \frac{\tau_{k+1}}{\bn(\bq)+k+2}\right) \ge \exp(-\lambda MT)\right)\\
&\le e^{\lambda MT} \prod_{k=1}^{n} \bE\left[\exp \left(- \frac{\lambda \tau_{k+1}}{\bn(\bq)+k+2}\right)\right]\\
&=e^{\lambda MT} \prod_{k=1}^{n} \left(1+\frac{\lambda}{\bn(\bq)+k+2}\right)^{-1}.
\end{align*}
Using the fact that $\log(1+x)\ge x -\frac{x^2}{2}$ for any positive 
$x $, we have 
\begin{align*}
\log\left(\prod_{k=1}^{n} \left(1+\frac{\lambda}{\bn(\bq)+k+2}\right)\right) \ge \lambda(H_{\bn(\bq)+n+2}-H_{\bn(\bq)+2})-\frac{\lambda^2}{2}\left(\sum_{k=\bn(\bq)+3}^{\bn(\bq)+n+2} \frac{1}{k^2}\right),
\end{align*}
where $H_n:=\sum_{k=1}^{n} \frac{1}{k}$ is the harmonic series. It is well-known however that
\[
\log(n+1) \le H_n \le \log n+1.
\]
Hence,
\[
\log\left(\prod_{k=1}^{n} \left(1+\frac{\lambda}{\bn(\bq)+k+2}\right)\right) \ge\lambda \log\left(\frac{\bn(\bq)+n+3}{\bn(\bq)+2}\right)-\lambda-\frac{\lambda^2 \pi^2}{12},
\]
which in turn implies
\[
\bP(N^{T}(\bq)>n) \le e^{\lambda MT-\lambda \log\left(\frac{\bn(\bq)+n+3}{\bn(\bq)+2}\right)+\lambda+\frac{\lambda^2 \pi^2}{12}}.
\]
In particular for $\lambda=2$, we get the bound
\[
\bP(T_{n+1}<T)=\bP(N^{T}(\bq)>n) \le e^{2MT+2+\frac{\pi^2}{3}}\left(\frac{\bn(\bq)+2}{\bn(\bq)+n+3}\right)^2.
\]
This certainly implies the second inequality in \eqref{lem21}. Finally
\[
\bE[N^{T}(\bq)] \le \frac{\pi^2}{6} e^{2MT+2+\frac{\pi^2}{3}} (\bn(\bq)+2)^2<\i,
\]
which implies the first inequality in \eqref{lem21}. \qed

\ms
By construction, the process $(t,\bq(t))_{t \ge 0}$ is a piecewise-deterministic process in the sense of Davis [Da], and thus we have the following proposition
\begin{prop}
The process $(t,\bq(t))_{t \ge 0}$ has the strong Markov property.
\end{prop}

\ms
It is not homogeneous because of the time dependence of the rates. Due to this Markovian property, we can talk about the stochastic flow 
$\Psi_{s}^{t}\bq$ for any $s \le t$ being defined as the realization of the particle system $\bq(t)$ at time $t$ conditioned to start at time $s$ at 
$\bq(s)=\bq$. The Markov property, ensures that this stochastic flow enjoys the semigroup property in distribution
\begin{align}\label{flow-stochastic}
\Psi_{t_1}^{t_3}\bq \ed \Psi_{t_2}^{t_3}\Psi_{t_1}^{t_2}\bq,
\end{align}
for any $t_1 \le t_2 \le t_3$, where on the right-hand side of \eqref{flow-stochastic} the stochastic flow $\Psi_{t_2}^{t_3}$ is independent of $\Psi_{t_1}^{t_2}$. We can also replace the times $t_i$'s by appropriate stopping times due the strong Markov property. The Markovian nature of our construction comes essentially from the memoryless property of exponential random variables.

\section{Forward Equation}
Recall the measure $\mu$ that was defined by \eqref{eq2.5}. 
The goal of this section is to prove the following theorem.

\begin{theorem}\label{main}
For any measurable function $G$ and $t \ge 0$, we have
\begin{equation}\label{id-dist}
 \int_{\Omega} \bE[G(\Psi_{0}^{t}\bq)]\ \mu(d\bq,0)
= \int_{\Omega} G(\bq)\ \mu(d\bq,t).
\end{equation}
\end{theorem}

\ms
Theorem 3.1 proves that if $\bq(0)$ is distributed according to the measure $\mu(d\bq,0)$, then $\bq(t)$ has the law $\mu(d\bq,t)$. By a density argument, it suffices to prove the equality \eqref{id-dist} for a suitable class of functions $G$ on $\Omega$, of the form 
\begin{equation}\label{expression-G}
G(\bq)=\exp \left(\int_{a^-}^{a^+} J(z)\cV(\bq)(z)\ dz \right),
\end{equation}
where $\cV$ was defined in \eqref{eq2.1}, 
and 
$J$ is a continuous function on $[a^-,a^+]$. For \eqref{id-dist},
it is enough to show 
\begin{align}\label{deriv}
\frac{d}{ds} \int_{\Omega} \bE[G(\Psi_{s}^{t}\bq)]\
\mu(d\bq,s)=0,\ \text{ for all } \ s \in [0,t].
\end{align}
From now on, we fix $t \ge 0$, and define the function
\[
G(\bq,s)=\bE[G(\Psi_{s}^{t}\bq)] \ \text{ for all } 0 \le s \le t .
\]

\subsection{Lipschitzness of $G$}
We will start by proving the following crucial theorem.
\begin{theorem}\label{lipschitz}
There exists a constant $C_1=C_1(P^-,P^+,V_{\i},J,a^-,a^+,t)>0$ such that
\[
|G(\bq,s)-G(\bq,s')| \le C_1(\bn(\bq)+2)^{2}\ |s-s'|,
\]
for all $s,s' \in [0,t]$ and for any $\bq \in \Omega$.
\end{theorem}
The proof of this Lipschitz property is carried out in two steps. 
These steps are formulated as Lemmas 3.1 and 3.2.
\begin{lemma}
Let $0\le s' \le s \le t$, and put $\theta:=s-s'$. There exists a constant $C_2=C_2(P^-,P^+,V_{\i},J,a^-,a^+,t)$ such that
\be\la{lem31}
\left|G(\bq,s')-\bE\left[G(\Psi_{s'}^{t-\theta}\bq)\right]\right| \le C_2(\bn(\bq)+2)^2 \theta.
\ee
\end{lemma}

\begin{lemma} There exists a constant $C'_2=C'_2(P^-,P^+,V_{\i},J,a^-,a^+,t)$ such that
\be\la{lem32}
\left|G(\bq,s)-\bE\left[G(\Psi_{s'}^{t-\theta}\bq)\right]\right| \le C'_2(\bn(\bq)+2)^2 \theta.
\ee
\end{lemma}

\bs\noi
{\bf Proof of Lemma 3.1} {\em (Step 1)}
In this step, we show that there exists a constant 
$C_3=C_3(P^-,P^+,\|J\|_{\i},a^-,a^+)$ such that 
\be\label{inequ-1}
\left|G(\bq,s')-\bE[G(\psi_{t-\theta}^{t}\Psi_{s'}^{t-\theta}\bq)\right|\le  C_3(\bn(\bq)+2)^2\theta.
\ee
We first use the Markov property to write
\[
G(\bq,s')=\bE[G(\Psi_{t-\theta}^{t}\Psi_{s'}^{t-\theta}\bq)].
\]
Let $\cE$ be the event  
\[
\cE:=\Big\{ \text{there exists a stochastic jump in } (t-\theta,t) \Big\}.
\]
Here by stochastic jump, we mean the creation of a new particle either at $z=a^-$ or at $z=a^+$, or the fragmentation of one of the particles.
We claim that there exists a constant $C_4=C_4(t)$ such that
\be\la{cE-bound}
\bP(\cE)\le C_4\left(\bn(\bq)+2\right)^2\theta.
\ee
 To see this, observe
\begin{align}\label{prob-jump}
\bP(\cE)=\bP\left(\int _{t-\theta}^{t} \mf{r}(u,\psi_{t-\th}^{u}\Psi_{s'}^{t-\theta}\bq)\ du \ge \tau \right),
\end{align}
where $\tau$ is a standard exponential random variable that is independent of $\{\Psi_{s'}^{v}\bq , s' \le v \le t-\theta \}$. Let $N_{s'}^{t-\theta}(\bq)$ be the number of stochastic jumps of the particle system started at time $s'$ at $\bq$ up to time $t-\theta$. By Lemma \ref{numb-stochas} we have that
\be\la{Nbound}
\bE[N_{s'}^{t-\theta}(\bq)] \le C_0(\bn(\bq)+2)^2,
\ee
where $C_0=C_0(t)$ is a constant that depends only on $t$. As the number of particles of $\psi_{t-\th}^{u}\Psi_{s'}^{t-\theta}\bq$ is at most $\bn(\bq)+N_{s'}^{t-\theta}(\bq)$ for all $u \in [t-\theta,t]$, then from \eqref{prob-jump}, it follows that
\begin{align*}
\bP(\cE) &\le \bP\Big( M\theta \left(\bn(\bq)+N_{s'}^{t-\theta}(\bq)+2\right) \ge \tau\Big)\\
&= \bE\left[1-e^{-M\theta(\bn(\bq)+N_{s'}^{t-\theta}(\bq)+2)}\right] \\
&\le M\theta \bE\left[\bn(\bq)+N_{s'}^{t-\theta}(\bq)+2\right]\\
 &\le C_4\left(\bn(\bq)+2\right)^2\theta,
\end{align*}
where the constant $M$ is the uniform bound on the rates defined in \eqref{bound-rates}, and $C_4= M(1+C_0)$. This completes the proof
of \eqref{cE-bound}.

 Evidently,
\begin{align}\label{decomp}
G(\bq,s')=\bE\left[\ind_{\cE^{c}}\ G(\psi_{t-\theta}^{t}\Psi_{s'}^{t-\theta}\bq)\right]+\bE\left[\ind_{\cE}\ G(\Psi_{s'}^{t}\bq)\right].
\end{align}
However, from the expression of $G$ in \eqref{expression-G}, we can find a constant $C_5=C(P^-,P^+,\|J\|_{\i},a^-,a^+)$ such that
\[
|G(\bq)| \le C_5 \ \text{ for all }\  \bq \in \Omega.
\]
From this and \eqref{decomp}, we learn
\[
\left|G(\bq,s')-\bE[G(\psi_{t-\theta}^{t}\Psi_{s'}^{t-\theta}\bq)\right|\le 2C_5\bP(\cE) \le C_3(\bn(\bq)+2)^2\theta,
\]
for a constant $C_3=2C_5C_4$.
This completes the proof of \eqref{inequ-1}.

\ms\noi
{\em (Step 2)} On account of \eqref{inequ-1}, it remains to show 
\begin{align}\label{inequ-2}
\left|\bE[G(\psi_{t-\theta}^{t}\Psi_{s'}^{t-\theta}\bq)]-\bE[G(\Psi_{s'}^{t-\theta}\bq)]\right|\le C_6(\bn(\bq)+2)^2\theta,
\end{align}
for some constant $C_6=C_6(P^-,P^+,V_{\i},J,t)$.
As a preparation, we first show that there exists a constant $C_7=C_7(P^-,P^+,V_{\i},J,a^-,a^+,t)$ such that
\begin{align}\label{lipschitz-flow}
\left|G(\psi_{0}^{r}\bq)-G(\bq)\right| \le C_7\bn(\bq) r, \
\text{ for all } \ r \in [0,t]\  \text{ and } \ \bq \in \Omega',
\end{align}
where the set $\O'$ is a set of full measure that was defined 
in Notation 2.1{\bf(vi)}. To prove \eqref{lipschitz-flow},
we fix $r>0$, and let $\rho:=\cV(\bq)$ and $\rho':=\cV(\psi_{0}^{r}\bq)$. As the exponential function is locally Lipschitz, for \eqref{lipschitz-flow} it suffices to show that there exists a constant $C_8=C_8(P^-,P^+,V_{\i},a^-,a^+,t),$ such that
\be\label{lipschitz-flow2}
\int_{a^-}^{a^+} |\rho'(z)-\rho(z)|\ dz \le C_8\bn(\bq) r.
\ee
Note that $\rho=\cV(\bar{\bq})$ where $\bar{\bq}$ is the particle 
configuration obtained from $\bq$ after deleting the redundant particles, so without loss of generality we can assume that $\bq \in \text{Int}(\Omega^n)$ for some $n \le \bn(\bq)$. Let 
\[
a^-<z_1<z_2<\cdots<z_n<a^+,
\]
be the discontinuity points of $\rho$. Let $\delta:=V_{\i}r$.
If $I$ is the interval
\[
I:=\bigcup_{i=1}^{n} [z_i-\delta,z_i+\delta],
\]
then for any $z \notin I$, we have $\rho'(z)=\rho(z)$, as the discontinuity 
points of the function $(z,v) \in [a^-,a^+] \times [0,r] \mapsto \cV(\psi_{0}^{v}\bq)(z)$ travel with speed at most $V_{\i}$. Therefore
\begin{align*}
\int_{a^-}^{a^+} |\rho'(z)-\rho(z)|\
dz &\le (\|\rho\|_{\i}+\|\rho'\|_{\i}) |I| \le 4\max(|P^-|,|P^+|)V_{\i}nr,
\end{align*}
which proves \eqref{lipschitz-flow2}. This in turn implies 
\eqref{lipschitz-flow}.

We are now ready to establish \eqref{inequ-2}.
We have that almost surely $\Psi_{s'}^{t-\theta}\bq \in \Omega'$, as this is equivalent to not having a stochastic jump at time $t-\theta$, an event that happens with probability one. This allows us to apply
\eqref{lipschitz-flow} to assert
\[
\bE \left[\left|G(\psi_{t-\theta}^{t}\Psi_{s'}^{t-\theta}\bq)-G(\Psi_{s'}^{t-\theta}\bq)\right| \right] \le C_7 \bE\left[\left(\bn(\Psi_{s'}^{t-\theta}\bq)+2\right)\right]\theta.
\]
This, the bound $\bn(\Psi_{s'}^{t-\theta}\bq)\le \bn(\bq)+N_{s'}^{t-\theta}(\bq)$, and  \eqref{Nbound} imply \eqref{inequ-2}.
 From \eqref{inequ-1} and \eqref{inequ-2}, we get \eqref{lem31}.
\qed

\bs
To finish the proof of Theorem \ref{lipschitz}, it remains to establish Lemma 3.2.
To achieve this, we will define a coupling of $(\Psi_{s'}^{t-\theta}\bq,\Psi_{s}^{t}\bq)$ for $t\ge s$.  Or equivalently, we define
a coupling of 
$(\Psi_{s'}^{t}\bq,\Psi_{s'+\th}^{t+\theta}\bq)$ for $t\ge s'$.

\subsubsection{Construction of the coupling} We fix $0 \le s' \le s $ and as before put $\theta:=s-s'$. We wish to construct two processes 
$(\bq(t))_{t \ge s'}$ and $(\bq'(t))_{t \ge s'}$ on the same probability space such that $(\bq(t))_{t \ge s'}$ has the law $(\Psi_{s'}^{t}\bq)_{t \ge s'}$ and $(\bq'(t))_{t \ge s'}$ has the law $(\Psi_{s'+\th}^{t+\theta}\bq)_{t \ge s'}$. \\

We start from a sequence of i.i.d exponential random variables of the form $(\tau_{i}^{j} , i \ge 1 , 1 \le j \le 3)$. We define first the \textit{coupling} rates as follows.\\

Let $t \ge s'$ and
for any $\rho \in [P^-,P^+]^2$, let
\begin{align*}
\mf{C}_{-}^{\text{coupling}}(t,\rho):=&\int_{L(\rho)} \mf{c}_{-}(t,\rho^{-},\rho) \wedge \mf{c}_{-}(t+\theta,\rho^{-},\rho)\ \b(d\rho^{-}),\\
\mf{C}_+^{\text{coupling}} (t,\rho) :=&\int _{R(\rho)} \mf{c}_+(t,\rho,\rho^{+}) \wedge \mf{c}_+(t+\theta,\rho,\rho^{+}) \ \b(d\rho^{+}),
\end{align*}
and for $\rho^{-},\rho^{+}\in [P^-,P^+]^2$ such that $\rho^{-} \prec \rho^{+}$, let
\[
\mf{F}^{\text{coupling}}(z,t,\rho^{-},\rho^{+}):=\int_{D(\rho^{-},\rho^{+})} \mf{f}(z,t,\rho^{-},\rho,\rho^{+})\wedge \mf{f}(z,t+\theta,\rho^{-},\rho,\rho^{+})\ \b(d\rho).
\]
The \textit{particle coupling rate} is defined as
\[
\mf{r}^{\text{coupling}}(t,\bq):=\mf{C}_{-}^{\text{coupling}}(t,\rho^0)+\sum_{i=1}^{n} \mf{F}^{\text{coupling}}(z_i,t,\rho^{i-1},\rho^i)+\mf{C}^{\text{coupling}}_{+}(t,\rho^n),
\] 
for any 
\[
\bq=((a^-,\rho^0),(z_1,\rho^1),\cdots,(z_n,\rho^n)) \in \Omega^n.
\]
Now, define
\begin{align*}
T_{1}^{1}:=&\inf \left\{ t \ge s' : \int_{s'}^{t} \mf{r}^{\text{coupling}}(u,\psi_{s'}^{u}\bq)\ du \ge \tau_1^1 \right\},\\
T_{1}^{2} := &\inf \left\{ t \ge s' : \int_{s'}^{t} (\mf{r}(u,\psi_{s'}^{u}\bq)-\mf{r}^{\text{coupling}}(u,\psi_{s'}^{u}\bq))\ du \ge \tau_{1}^{2} \right\},\\
T_{1}^{3} := &\inf \left\{ t \ge s' : \int_{s'}^{t} (\mf{r}(u+\theta,\psi_{s'}^{u}\bq)-\mf{r}^{\text{coupling}}(u,\psi_{s'}^{u}\bq))\ du \ge \tau_{1}^{3}\right \}.
\end{align*}
Put $T_1=\min(T_1^1,T_1^2,T_1^3)$, and for $t \in [s',T_1)$, set $\bq(t)=\bq'(t)=\psi_{s'}^{t}\bq=\psi_{s'+\th}^{t+\theta}\bq$. Now, write
\[
\psi_{s'}^{T_1}\bq = ((a^-,\rho^0),(z_1,\rho^1),\cdots,(z_n,\rho^n)) \in \Omega^n.
\]
Conditionally on $T_1$, we consider the following cases:
\begin{itemize}
\item{If $T_1=T_{1}^{2}$, we set $\bq'(T_1)=\psi_{s'+\th}^{T_1+\theta}\bq$, and define $\bq(T_1)$ by making a stochastic jump (either a creation of a particle at $z=a^{\pm}$ or a fragmentation of a particle) using the rate $\mf{r}(T_1,\psi_{s'}^{T_1}\bq)-\mf{r}^{\text{coupling}}(T_1,\psi_{s'}^{T_1}\bq)$. For $t \ge T_1$, conditionally on $(\bq(T_1),\bq'(T_1))$, we let the two processes $(\bq(t))_{t \ge T_1}$ and $(\bq'(t))_{t \ge T_1}$ evolve independently with respectively the law of $(\Psi_{T_1}^{t}\bq(T_1))_{t \ge T_1}$ and $(\Psi_{T_1+\theta}^{t+\theta}\bq'(T_1))_{t \ge T_1}$. }
\item{If $T_1=T_{1}^{3}$, we do the same as previously by switching the roles of $\bq$ and $\bq'$ and using the rates at $T_1+\theta$ instead of $T_1$, in a way that only $\bq'$ makes a jump at time $T_1$ and not $\bq$.}
\item{ If $T_1=T_{1}^{1}$, both $\bq$ and $\bq'$ make a stochastic jump at $T_1$ using the rate $\mf{r}^{\text{coupling}}(T_1,\psi_{s'}^{T_1}\bq)$.
As $\bq(T_1)=\bq'(T_1)$, we redo the same process again by using now random variables $(\tau_2^{1},\tau_2^2,\tau_2^3)$ and defining $(T_2^1,T_2^2,T_2^3)$, etc ...} 
\end{itemize}
We claim that this defines a coupling. Indeed, each one of the processes $(\bq(t))_{t \ge s'}$ and $(\bq'(t))_{t \ge s}$ is a piecewise-deterministic Markov process by construction, and the law of the first jump time and the value it takes at that time is easily verified to be the same as that of the two processes $(\Psi_{s'}^{t}\bq)_{t \ge s'}$ and $(\Psi_{s'+\th}^{t+\theta}\bq)_{t \ge s'}$.  Using this coupling, we can now prove
Lemma 3.2, which in turn finishes the proof of Theorem \ref{lipschitz}.

\bs\noi
{\bf Proof of Lemma 3.2 }
Keeping the notation from the construction of the coupling, we will exhibit an event $\mathcal{E}$ on which $\bq(u)=\bq'(u)$ for all $u \in [s',t-\theta]$. Writing $T_i$ for the time at which the coupled process experiences a jump for the $i$-th time,
we see that as long as $T_i^1=T_i$ for all $i$ such that $T_i<t-\th$, then we must have $\bq(t-\th)=\bq'(t-\th)$, 
as we jump using the same value $\rho^*$ at every step. Let
\[
\mathcal{E}=\{ T_i^{1}=T_i \text{ for all } i \le N_t \},
\]
where $N_t=\sup\{ n \ge 0 : T_n < t-\th\}$. Evidently,
\[
\left|\bE\left[G(\Psi_{s}^{t}\bq)]-
\bE[G(\Psi_{s'}^{t-\theta}\bq)\right]\right| =
\left|\bE\left[G(\bq'(t-\th))]-
\bE[G(\bq(t-\th))\right]\right|
\le 2\|G\|_{\i} \bP(\mathcal{E}^c).
\]
Hence for \eqref{lem32}, it suffices to show that there is a constant $C_{9}=C_9(P^-,P^+,V_{\i},t)>0,$ such that
\be\la{Ec-bound}
\bP(\mathcal{E}^c) \le C_{9} \theta(\bn(\bq)+2)^2.
\ee
To achieve this, first observe, 
\[
\mathcal{E}^{c} \subset \bigcup_{n=1}^{\i} \Big\{  T_i^{1}=T_i \text{ for all } i \le n-1 , \ T_{n-1}<t-\th,\  T_n^{1}>T_n\Big\}.
\]
As a result,
\be\la{eq3.15}
\bP(\mathcal{E}^{c}) \le \sum_{n=1}^{\i} \bP(T_n^{1}-T_{n-1}>T_n-T_{n-1} ,\ T_{n-1}<t-\th).
\ee
Recall that we can write 
\[
T_n^{1}-T_{n-1}=\inf \left\{ v \ge 0 : \int_{T_{n-1}}^{T_{n-1}+v} \mf{r}^{\text{coupling}}(u,\psi_{s'}^{u}\bq)\ du \ge \tau_{n}^1 \right \},
\]
and similarly for $T_n^{2}-T_{n-1}$ and $T_n^{3}-T_{n-1}$, with appropriately replacing the rates $\mf{r}$ used inside the integral. Note that the conditional probability
\[
P_n:=\bP\Big(T_n^{1}-T_{n-1}>T_n-T_{n-1}\  \Big\vert \ T_{n-1}\Big),
\]
satisfies
\begin{align}\nonumber
P_n&=\bE\left[1-\exp \left(-\int_{T_{n-1}}^{T_n^1} \left|\mf{r}(u,\psi_{s'}^{u}\bq)-\mf{r}(u+\theta,\psi_{s}^{u+\theta}\bq)\right|\ du \right) \bigg\vert\ T_{n-1}\right] \\
&\le \bE\left[\int_{T_{n-1}}^{T_n^1} \left|\mf{r}(u,\psi_{s'}^{u}\bq)-\mf{r}(u+\theta,\psi_{s'}^{u}\bq)\right|\ du \
\bigg\vert\ T_{n-1}\right].\la{eq3.16}
\end{align}
Now, for any fixed $\bq \in \Omega$, we have that
\[
\partial_t \mf{r}(t,\bq)=\mathcal{R}^-(t,\bq)+\sum_{i=1}^{\bn(\bq)} \mathcal{R}_i(t,\bq)+\mathcal{R}^+(t,\bq)
\]
where
\begin{align*}
\mathcal{R}^-(t,\bq):=&\int_{L(\rho^0)} [\rho^-,\rho_0]^{-} \
\partial_t \left(\frac{\ell(a^-,t,\rho^-)f(a^-,t,\rho^-,\rho^0)}{\ell(a^-,t,\rho^0)} \right)\ \beta(d\rho^-),\\
\mathcal{R}_i(t,\bq):=&\int_{D(\rho^{i-1},\rho^{i})} [\rho^{i-1},\rho,\rho^{i}]^-\ \partial_t \left( \frac{f(z_i,t,\rho^{i-1},\rho)f(z_i,t,\rho,\rho^i)}{f(z_i,t,\rho^{i-1},\rho^i)} \right)\ \beta(d\rho),\\
\mathcal{R}^+(t,\bq):=&\int_{R(\rho^{\bn(\bq)})} [\rho^{\bn(\bq)},\rho^+]^{+} \ \partial_t \left(f(a^+,t,\rho^{\bn(\bq)},\rho^+)\right)\ \beta(d\rho^+).
\end{align*}
Hence by using the uniform upper and lower bound of the kernel $f$ and $\ell$, and the uniform upper bound of
their first-order derivatives, 
there exists a uniform constant $M'>0$ such that
\[
\partial_t \mf{r}(t,\bq) \le M' \mathcal{A}(\bq),
\]
where $\mathcal{A}(\bq)$ is given by
\[
\int_{L(\rho^0)} [\rho^{-},\rho^0]^{-}\   \beta(d\rho^{-})+\sum_{i=1}^{\bn(\bq)} \int_{D(\rho^{i-1},\rho^i)} [\rho^{i-1},\rho,\rho^{i}]^{-}\ \beta(d\rho)+\int_{R(\rho^{\bn(\bq)})} [\rho^{\bn(\bq)},\rho^{+}]^{+}\ \beta(d\rho^{+}).
\]
From this and \eqref{eq3.16} we deduce, 
\be\la{eq3.17}
P_n \le M' \theta \ \bE \left[ \int_{T_{n-1}}^{T_n^{1}} \cA(\psi_{s'}^{u}\bq)\ du \ \Big | \ T_{n-1} \right].
\ee
Moreover, using the uniform upper and lower bound of the kernel $f$ and $\ell$, we can find  $\d_1>0$ such that
\be\la{eq3.18}
\mf{r}^{\text{coupling}}(u,\bq)\ge \d_1\cA(\bq).
\ee
On the other hand, since by definition,
\[
\tau_{n}^1=\int_{T_{n-1}}^{T_n^{1}}\mf{r}^{\text{coupling}}(u,\psi_{s'}^{u}\bq)\ du ,
\]
we use \eqref{eq3.18} to assert
\[
\tau_n^{1} \ge \d_1 \int_{T_{n-1}}^{T_n^1} \cA(\psi_{s'}^{u}\bq)\ du.
\]
This and \eqref{eq3.17} yield,
\[
P_n \le \d_1^{-1}M'\theta\   \bE[\tau_{n}^1 | T_{n-1}]=\d_1^{-1}M' \theta,
\]
as $\tau_{n}^1$ is independent of $T_{n-1}$. From this and \eqref{eq3.15} we learn
\[
\bP\left(\mathcal{E}^{c}\right) \le \d_1^{-1}M'\theta \sum_{n=1}^{\i} \bP(T_{n-1}<t).
\]
From this and \eqref{lem21} we deduce \eqref{Ec-bound}.
This completes the proof of \eqref{lem32}.
\qed

\subsection{Differentiation of $G$}
Recall that for Theorem 2.1, it suffices to verify \eqref{deriv}.
To carry out the differentiation in \eqref{deriv}, we first learn how to
differentiate the integrand $G$ with respect to $s$. Since
 $G(\bq,s)=\bE [G(\Psi_s^t\bq)]$, we expect a Kolmogorov type equation
of the form
\be\la{forward}
G_s(\bq,s)=-\bE [(\cL G)(\Psi_s^t\bq)],
\ee
 where $\cL$ denotes the generator of the process $\bq(\cdot)$. 
Since the deterministic part of the dynamics is discontinuous, the verification of \eqref{forward} poses some challenges that 
are handled in this subsection. The integrated version of \eqref{forward}
is our next result.

\begin{theorem}
For any $s \ge 0$ we have that the limit
\be\la{lim}
\lim_{s' \uparrow s} \frac{1}{s-s'} \int (G(\bq,s)-G(\bq,s'))\ \mu(d\bq,s),
\ee
equals to
\[
 - \sum_{n=0}^{\i}( \Gamma_n^{-}(s)+\Gamma_n^{+}(s))+
\sum_{n=1}^{\i}
(\Gamma_n^{\mf{f}}(s)+\Gamma_n^{t}(s)+\Gamma_n^{\text{d}}(s)),
\]
where for $n \ge 1$, we have
\begin{align*}
\Gamma_n^{-}(s)&=\int_{\Omega^n} \int_{L(\rho^0)} \mf{c}_{-}(s,\rho^*,\rho^0)(G(E_{-}^{\rho^*}\bq,s)-G(\bq,s))\ \b(d\rho^*)\  \mu^{n}(d\bq,s),\\
\Gamma_n^{+}(s)&=\int_{\Omega^n} \int_{R(\rho^n)} \mf{c}_+(s,\rho^n,\rho^*)(G(E_+^{\rho^*}\bq,s)-G(\bq,s))\ \b(d\rho^*)\mu^{n}(d\bq,s),\\
\Gamma_n^{\mf{f}}(s)&=\sum_{i=1}^{n} \int_{\Omega^n} \int_{D(\rho^{i-1},\rho^i)} \mf{f}(z_i,s,\rho^{i-1},\rho^*,\rho^i)(G(E_{i}^{\rho^*}\bq,s)-G(\bq,s))\ \b(d\rho^*)\mu^{n}(d\bq,s),\\
\Gamma_n^{t}(s)&= \sum_{i=1}^n\int_{\Omega^n}  [\rho^{i-1},\rho^i] \left(\frac{f_z(z_i,s,\rho^{i-1},\rho^i)}{f(z_i,s,\rho^{i-1},\rho^i)}+\lambda(z_i,s,\rho^i)-\lambda(z_i,s,\rho^{i-1})\right)G(\bq,s) \b(d\rho^*)\mu^{n}(d\bq,s),
\end{align*}
and we can write $\Gamma_n^{d}(s)=\Gamma_n^{d,-}(s)+\Gamma_n^{d,\mf{f}}(s)+\Gamma_n^{d,+}(s),$ with
\begin{align*}
\Gamma_n^{d,+}(s)=&-\int_{\Omega^{n-1}} \int_{R(\rho^{n-1})} [\rho^{n-1},\rho^*]f(a^+,s,\rho^{n-1},\rho^*)G(E_+^{\rho^*}\bq,s)\ \b(d\rho^*)\mu^{n-1}(d\bq,s),\\
\Gamma_n^{d,-}(s)=&\int_{\Omega^{n-1}} \int_{L(\rho^0)} [\rho^*,\rho^0]\frac{\ell(a^-,s,\rho^*)f(a^-,s,\rho^*,\rho^0)}{\ell(a^-,s,\rho^0)}G(E_{-}^{\rho^*}\bq,s)\ \b(d\rho^*)\mu^{n-1}(d\bq,s),\\
\Gamma_n^{d,\mf{f}}(s)=&\sum_{i=1}^{n-1} \int_{\Omega^{n-1}} \int_{D(\rho^{i-1},\rho^i)} [\rho^{i-1},\rho^*,\rho^i]
\frac{f(z_i,s,\rho^{i-1},\rho^*)f(z_i,s,\rho^*,\rho^i)}{f(x_i,s,\rho^{i-1},\rho^i)}\\
&\ \ \ \ \ \ \ \ \ \ \ \ \ \ \ \ \ \ \ \ \ \ \ \ \ \ \ \ \ \ \ \ \ \ \ \ \ \ 
\ \ \ \ \ \ \ \ \ \ \ \ \ \ \
G(E_i^{\rho^*}\bq,s)\ \b(d\rho^* )\mu^{n-1}(d\bq,s).
\end{align*}
For $n=1$, $\Gamma_n^{d,\mf{f}}(s)=0$, and without ambiguity the other terms have the same expression 
as for $n > 1$. 
\end{theorem}

\ms\noi
{\bf Proof} {\em (Step 1)} In this step, we use the Markov property
to derive a formula for $G(\bq,s')$ (see \eqref{decom1} at the end of this 
step).
To begin, observe that by the Markov property, 
\be\la{markov}
G(\bq,s')=\bE[G(\Psi_{s'}^{t}\bq)]=\bE[G(\Psi_{s}^{t}\Psi_{s'}^{s}\bq)]=\bE[G(\Psi_{s'}^{s}\bq,s)]
\ee
Using the notation that we used previously in our construction of the Markov process, let
\[
T_1=\inf \left\{ t \ge s' : \int_{s'}^{t} \mf{r}(u,\psi_{s'}^{u}\bq)\ du \ge \tau_1 \right\},
\]
where $\tau_1$ is an independent standard exponential random variable. Consider again $\mathcal{E}$ to be the event $\{ T_1 \in [s',s] \}$, then
\be\la{eq3.21}
\bE[G(\Psi_{s'}^{s}\bq,s)]=\bE[G(\Psi_{s'}^{s}\bq,s) \ind_{\mathcal{E}}]+\bP(\mathcal{E}^{c})G(\psi_{s'}^{s}\bq,s).
\ee
For the first term, we can use the strong Markov property at the stopping time $T_1$ to get
\be\la{eq3.22}
\bE\left[G(\Psi_{s'}^{s}\bq,s)\ind_{\mathcal{E}}\right]=
\bE\left[ G(\Psi_{s'}^{T_1}\bq,T_1)\ind_{\mathcal{E}}\right].
\ee
Using Theorem \ref{lipschitz}, we have
\begin{align*}
\left|\bE\left[(G(\Psi_{s'}^{T_1}\bq,T_1)-G(\Psi_{s'}^{T_1}\bq,s))\ind_{\mathcal{E}}\right]\right| &\le C_1\bE\left[\left(\bn\left(\Psi_{s'}^{T_1}\bq\right)+2\right)^2|T_1-s'|\ind_{\mathcal{E}}\right] \\
&\le C_1(\bn(\bq)+3)^2\  \bP(\mathcal{E})\ (s-s').
\end{align*}  
We certainly have
\begin{align}\nonumber
\bP(\mathcal{E})=&\bP \left(\tau_1 \le \int_{s'}^{s} \mf{r}(u,\psi_{s'}^{u}\bq)du \right)=1-\exp \left(-\int_{s'}^{s} \mf{r}(u,\psi_{s'}^{u}\bq)\ du \right)\\
&\le \int_{s'}^{s} \mf{r}(u,\psi_{s'}^{u}\bq)\ du\le M(\bn(\bq)+2)(s-s'),
\la{r-bound}
\end{align}
where $M$ is a uniform bound on the rates.
From this and the previous display we learn
\[
\left|\bE\left[\left(G\left(\Psi_{s'}^{T_1}\bq,T_1\right)-
G\left(\Psi_{s'}^{T_1}\bq,s\right)\right)\ind_{\cE}\right]\right| \le C_1M(\bn(\bq)+3)^3 (s-s')^2.
\]
From this, \eqref{markov}, \eqref{eq3.21}, and \eqref{eq3.22}
we deduce
\be\la{decom1}
G(\bq,s')=\bE[G(\Psi_{s'}^{T_1}\bq,s)\ind_{\cE}]+
\bP(\cE^{c})G(\psi_{s'}^{s}\bq,s)+R_1(\bq)(s-s')^2,
\ee
where $|R_1(\bq)| \le C_1M(\bn(\bq)+3)^3$.

\ms\noi
{\em (Step 2)} The main goal of this step is to use \eqref{decom1}
to establish the following decomposition:
\begin{align}\nonumber
G(\bq,s)-G(\bq,s')=&\cS'(\bq)+\cS^-(\bq)+\cS^+(\bq)
+\sum_{i=1}^{\bn(\bq)}\cS_i(\bq)\\
&+R_2(\bq)\left[(s-s')^2+(s-s')\1(\bq\in\hat\O)\right],\label{decompos}
\end{align}
where 
\begin{align*}
\cS'(\bq)=&\exp\left(-\int_{s'}^{s} \mf{r}(u,\psi_{s'}^{u}\bq)\ du\right)(G(\bq,s)-G(\psi_{s'}^{s}\bq,s)),\\
\cS^-(\bq)= &\int_{s'}^{s} \int_{L(\rho^0)} \mf{c}_{-}
(\th,\rho^*,\rho^0)(G(\bq,s)-G(E_{-}^{\rho^*}
\psi_{s'}^{\th}\bq,s))\ \b(d\rho^*)d\th,\\
\cS^+(\bq)=&\int_{s'}^{s}\int_{R(\rho^{\bn(\bq)})}
 \mf{c}_+(\th,\rho^{\bn(\bq)},\rho^*)(G(\bq,s)-G(E_{+}^{\rho^*}\psi_{s'}^{\th}\bq,s))\ \b(d\rho^*)d\th,\\
 \cS_i(\bq)=& \int_{s'}^{s}\int_{D(\rho^{i-1},\rho^i)} \mf{f}(z_i,\th,\rho^{i-1},\rho^*,\rho^i)(G(\bq,s)-G(E_{i}^{\rho^*}\psi_{s'}^{t}\bq,s))\ \b(d\rho^*)d\th,
\end{align*} 
and the term $R_2$ satisfies the bound 
\[
|R_2(\bq)| \le c_1(\bn(\bq)+3)^3,
\]
 for a constant
$c_1=c_1(P^-,P^+,J,a^-,a^+,t)>0$ that does not depend on $\bq$, with the set $\hat\Omega$ is defined as
\[
\hat \Omega=\hat \O_{s'}^{s}:=\Big\{ \bq \in \Omega :\ u \mapsto \psi_{s'}^u\bq \text{ experiences a collision in } [s',s]\Big \}.
\]
We can readily show that there exists a universal constant $c_2$ such that
\be\la{hato}
\mu(\hat\Omega_{s'}^{s},s) \le c_2(s-s').
\ee
and so upon integrating with respect to $\mu(d\bq,s)$, the error terms are all of order $O((s-s')^2)$.

To achieve  \eqref{decompos}, we first examine 
the first term on the right-hand side of \eqref{decom1}.
 From the boundedness of $G$ and \eqref{r-bound} we deduce
\be\la{decom3}
\bE\left[G\left(\Psi_{s'}^{T_1}\bq,s\right)\ind_{\mathcal{E}}\right]=
\bE\left[G\left(\Psi_{s'}^{T_1}\bq,s\right)\ind_{\mathcal{E}}\right]\1(\bq\notin\hat\O)+R_3(\bq)\1(\bq\in\hat\O)(s-s') ,
\ee
with $|R_3(\bq)| \le c_3(\bn(\bq)+2),$ for a constant
$c_3=c_3(t)$. Moreover for $\bq\notin\hat\O$,
\begin{align*}
\bE\left[G\left(\Psi_{s'}^{T_1}\bq,s\right)\ind_{\mathcal{E}}\right]=&\bE\left[\frac{\mf{C}_{-}(T_1,\rho^0)}{\mf{r}(T_1,\psi_{s'}^{T_1}\bq)}G\left(E_{-}^{\rho^*}\psi_{s'}^{T_1}\bq,s\right)
\ind_{\mathcal{E}}\right] \\
&+\bE\left[\frac{\mf{C}_+(T_1,\rho^{\bn(\bq)})}{\mf{r}(T_1,\psi_{s'}^{T_1}\bq)}G\left(E_{+}^{\rho^*}\psi_{s'}^{T_1}\bq,s\right)\ind_{\mathcal{E}}\right]\\
& +\sum_{i=1}^{\bn(\bq)} \bE\left[ \frac{\mf{F}(z_i-[\rho^{i-1},\rho^i](T_1-s'),T_1,\rho^{i-1},\rho^i)}{\mf{r}(T_1,\psi_{s'}^{T_1}\bq)}G\left(E_{i}^{\rho^*}\psi_{s'}^{T_1}\bq,s\right)\ind_{\mathcal{E}}
\right]
\\
&\ \ \ =:\cT_-+\cT_++\sum_{i=1}^{\bn(\bq)}\cT_i,
\end{align*}
where each $\rho^*$ is distributed according to the density previously described in the construction of the stochastic flow. 
Note that the distribution function of $T_1$ is given by
\[
\bP(T_1 \ge \th) = \exp \left(-\int_{s'}^{\th} \mf{r}(u,\psi_{s'}^{u}\bq)\ du \right),
\]
or equivalently,
\be\la{T1}
\bP(T_1 \in d\th)=\mf{r}(\th,\psi_{s'}^{\th}\bq)
\exp \left(-\int_{s'}^{\th} \mf{r}(u,\psi_{s'}^{u}\bq)\ du \right)d\th
\ee
We can certainly write
\begin{align*}
\cT_-=&\bE\left[\ind_{\mathcal{E}}\  \int_{L(\rho^0)} 
\frac{\mf{c}_{-}(T_1,\rho^*,\rho^0)}{\mf{r}(T_1,\psi_{s'}^{T_1}\bq)}G\left(E_{-}^{\rho^*}\psi_{s'}^{T_1}\bq,s\right)
\ \b(d\rho^*)\right]\\
=&\int_{s'}^{s} \int_{L(\rho^0)} \mf{c}_{-}(\th,\rho^*,\rho^0)
\exp \left(-\int_{s'}^{\th} \mf{r}(u,\psi_{s'}^{u}\bq)\ du \right) 
G\left(E_{-}^{\rho^*}\psi_{s'}^{\th}\bq,s\right)\ \b(d\rho^*)d\th,
\end{align*}
where we have used \eqref{T1} for the second equality.
On the other hand, from the boundedness of 
$G$, we learn that  the expression
\[
\left|\int_{s'}^{s} \int_{L(\rho^0)} \mf{c}_{-}(\th,\rho^*,\rho^0)
\left(\exp \left(-\int_{s'}^{\th} \mf{r}(u,\psi_{s'}^{u}\bq)du \right)-1\right)
G\left(E_{-}^{\rho^*}\psi_{s'}^{\th}\bq,s\right)\ \b(d\rho^*)d\th\right|, 
\]
is bounded above by $c_4(s-s')^2$
for some constant $c_4=c_4(P^-,P^+,J,a^-,a^+)$. This in turn implies
\begin{align*}
\cT_-\ \1(\bq\notin\hat\O)=&\int_{s'}^{s} \int_{L(\rho^0)} 
\mf{c}_{-}(\th,\rho^*,\rho^0) G\left(E_{-}^{\rho^*}\psi_{s'}^{\th}\bq,s\right)\ \b(d\rho^*)d\th\
\ \1(\bq\notin\hat\O)+R_4(\bq)(s-s')^2\\
=&\int_{s'}^{s} \int_{L(\rho^0)} 
\mf{c}_{-}(\th,\rho^*,\rho^0) G\left(E_{-}^{\rho^*}\psi_{s'}^{\th}\bq,s\right)\ \b(d\rho^*)d\th\\
&+R_5(\bq)\left[(s-s')^2+(s-s')\1(\bq\in\hat\O)\right],
\end{align*}
where $R_4(\bq),R_5(\bq)\le c_5$, for a constant $c_5$.
We treat the terms $\cT_+$ and $\cT_i$ in the same fashion. For example,
\begin{align*}
\cT_i=&\int_{s'}^{s}\int_{\rho^{i-1}}^{\rho^i} \mf{f}(z_i-[\rho^{i-1},\rho^i](\th-s'),\th,\rho^{i-1},\rho^*,\rho^i)G(E_{i}^{\rho^*}\psi_{s'}^{\th}\bq,s)\ d\th d\b(\rho^*)+R_5(\bq)(s-s')^2\\
=&\int_{s'}^{s}\int_{\rho^{i-1}}^{\rho^i} \mf{f}(z_i,\th,\rho^{i-1},\rho^*,\rho^i)G(E_{i}^{\rho^*}\psi_{s'}^{\th}\bq,s)\ d\th \b(d\rho^*)+R_6(\bq)(s-s')^,
\end{align*}
where $R_5(\bq),R_6(\bq)\le c_6$, for a constant $c_6$.
 Here for the last equality, we have used the Lipschitzness of
the rate $\mf{f}$. This in turn implies
\begin{align*}
\cT_i\ \1(\bq\notin\hat\O)=&\int_{s'}^{s}\int_{\rho^{i-1}}^{\rho^i} \mf{f}(z_i,\th,\rho^{i-1},\rho^*,\rho^i)G(E_{i}^{\rho^*}\psi_{s'}^{\th}\bq,s)\ d\th \b(d\rho^*)\\
&+R_7(\bq)\left[(s-s')^2+(s-s') \1(\bq\notin\hat\O)\right],
\end{align*}
 where $R_7(\bq)\le c_7$, for a constant $c_7$.
From these representations of $\cT_\pm\ \1(\bq\notin\hat\O)$ and 
$\cT_i\ \1(\bq\notin\hat\O)$, \eqref{r-bound},
\eqref{decom1} and \eqref{decom3} we deduce
\begin{align}\nonumber
G(\bq,s')=&\int_{s'}^{s} \int_{L(\rho^0)} \mf{c}_{-}(\th,\rho^*,\rho^0)
G\left(E_{-}^{\rho^*}\psi_{s'}^{t}\bq,s\right)\ d\th d\rho^*\\
\la{decom2}
&+\sum_{i=1}^{\bn(\bq)} \int_{s'}^{s}\int_{\rho^{i-1}}^{\rho^i} \mf{f}(z_i,\th,\rho^{i-1},\rho^*,\rho_i)G(E_{i}^{\rho^*}\psi_{s'}^{\th}\bq,s)\ d\th d\rho^*\\
\nonumber&+\int_{s'}^{s}\int_{\rho^{\bn(\bq)}}^{P^+} \mf{c}_+(\rho^{\bn(\bq)},\rho^*,\th)G(E_{+}^{\rho^*}\psi_{s'}^{\th}\bq,s)\ d\th d\rho^*\\
&+e^{-\int_{s'}^{s} \mf{r}(\th,\psi_{s'}^{t}\bq)\ d\th}\
G\left(\psi_{s'}^s\bq,s\right)+R_1(\bq)(s-s')^2\nonumber\\
&+R_8(\bq)\left[(s-s')^2+(s-s') \1(\bq\notin\hat\O)\right],\nonumber
\end{align}
where $|R_8(\bq)| \le c_8$, for a constant $c_8$. 
On the other-hand,
\begin{align*}
1=&e^{-\int_{s'}^{s} \mf{r}(\th,\psi_{s'}^{t}\bq)\ d\th}+\int_{s'}^{s} \mf{r}(\th,\psi_{s'}^{\th}\bq)\ d\th+R_9(\bq)(s-s')^2\\
=&e^{-\int_{s'}^{s} \mf{r}(\th,\psi_{s'}^{\th}\bq)\ d\th}
+\1(\bq\notin\hat\O)\int_{s'}^{s} \mf{r}(\th,\psi_{s'}^{\th}\bq)\ d\th\\
&+R_{10}(\bq)\left[(s-s')^2+(s-s')\1(\bq \in \hat\O)\right],\\
\end{align*}
with $|R_9(\bq)|,|R_{10}(\bq)|\le c_9(\bn(\bq)+2)^2$, for 
a constant $c_9$ that is independent of $\bq$. We now
use the Lipschitzness of our rate
$\mf{f}$ to replace $z_i-[\rho^{i-1},\rho^i](t-s')$ with $z_i$.
As a result,
\begin{align}\nonumber
1-e^{-\int_{s'}^{s} \mf{r}(\th,\psi_{s'}^{\th}\bq)\ d\th}=&\int_{s'}^{s} \int_{L(\rho^0)} 
\mf{c}_{-}(\th,\rho^*,\rho^0)\ \b(d\rho^*)d\th\\
\la{1}
&+\sum_{i=1}^{\bn(\bq)} 
\int_{s'}^{s}\int_{D(\rho^{i-1},\rho^i)} 
\mf{f}(z_i,\th,\rho^{i-1},\rho^*,\rho^i)\ \b(d\rho^*)d\th\\
\nonumber&+\int_{s'}^{s}\int_{R(\rho^{\bn(\bq))} }
\mf{c}_+(z_i,\th,\rho^{i-1},\rho^*,\rho^i)\ \b(d\rho^*)d\th\\
&+
R_{11}(\bq)\left[(s-s')^2+(s-s')\1(\bq \in \hat\O)
\right],\nonumber
\end{align}
with again $|R_{11}(\bq)|$ bounded by a constant multiple
of $(\bn(\bq)+2)^2$.
Here again we have used the Lipschitzness of our rate
$\mf{f}$ to replace $z_i-[\rho^{i-1},\rho^i](t-s')$ with $z_i$.
We now multiply both sides of \eqref{1} by $G(\bq,s)$ and subtract 
the outcome from  \eqref{decom2} to arrive at \eqref{decompos}.

\ms\noi
{\em(Step 3)}
Fix now $n \ge 1$, and let us analyze each term of the sum
in \eqref{decompos}, integrated against the probability measure 
$\mu^{n}(d\bq,s)$. We start from $\cS'$, and focus on the spatial integration. To prepare for this, we need some definitions.
Let us write 
$u_i,\ i\in \llbracket 0,n \rrbracket$ for the relative velocities
of the particles:
\[
u_0:=v_1, \ \ \ \ u_n:=-v_n,\ \ \ \ u_i:=v_{i+1}-v_i.
\]
We also write $U$ for the set of particle configurations $\bz\in
\D^n$ that do not experience any collision in the interval $[s',s]$, and define
the sets $A_i$ and $B_i$, $i\in \llbracket 0,n \rrbracket$ by:
\begin{align*}
A_0=&\Big\{\bz \in {\Delta^n} : z_1<a^--u_0(s-s') \Big\},
\ \ \ \ \ \ \ \ \ 
B_0=\Big\{ \bz \in {\Delta^n} : \ z_1<a^-+|u_0|(s-s')
\Big\},\\
A_n=&\Big\{\bz \in{\Delta^n} : \
a^+-u_n(s-s')<z_n\Big\},\  \ \ \ \ \ \ \ 
B_n=\Big\{\bz \in {\Delta^n} : \
a^+-|u_n|(s-s')<z_n\Big\},\\
A_i=&\Big\{ \bz \in {\Delta^n} :\ z_{i+1}-z_i<-u_i(s-s')
\Big\},   \ \ \ \ \
B_i=\Big\{ \bz \in {\Delta^n} :\ z_{i+1}-z_i<|u_i|(s-s')
\Big\}.
\end{align*}
for $i \in \llbracket 1,n-1 \rrbracket$. Note that
the action of the flow 
$\psi_{s'}^{s}$ on the set $U$ is simply a translation in the $\bz$-space:
\[
\psi_{s'}^{s}( \bz,\brho)=:( \phi_{s'}^{s}( \bz),\brho),\ \ \ \ 
\phi_{s'}^{s}( \bz)=(z_i+v_i(s-s'))_{i=1}^{n},
\]  
where $v_i=-[\rho^{i-1},\rho^{i}]$. 
Then the set $U$ can be expressed as
$U={\Delta^n} \setminus \bigcup_{i=0}^{n} A_i$.
Moreover, writing $|A|$ for the Lebesgue measure of the set $A$, it is not 
hard to show that there exists a constant $c_{10}$ such that
\begin{align}\la{leb1}
&|A_i|\le|B_i|\le c_{10}(s-s')\ \frac{(a^+-a^-)^{n}}{n!},\\
&|A_i\cap A_j|\le|B_i\cap B_j|\le c_{10}(s-s')^2\ \frac{(a^+-a^-)^{n}}{n!}.\la{leb2}
\end{align}

 In the present step, we  fix  
$\brho:=(\rho^0,\rho^1,\cdots,\rho^n)$, satisfying 
$\rho^0\prec\rho^1\prec\cdots\prec\rho^n$,
and focus on the integration with respect to the space variable 
$\bz:=(z_1,\cdots,z_n) \in {\Delta^n}$. More specifically
we will show
\begin{align}\nonumber
\int_{{\Delta^n}} \cS'(\bq) g^n(\bz,\brho)\ d\bz=&
\sum_{i=0}^{n} \ind_{\{u_i>0\}} \int_{ B_i} G(\bq,s)g^n(\bz,\brho)\ d\bz-\sum_{i=0}^{n} \ind_{\{ u_i<0\}} \int_{ B_i} G(\psi_{s'}^{s}\bq,s)g^n(\bz,\brho)\ d\bz\\
\label{combination1}
&+\int_{U}G(\bq,s)\left(g^n(\bz,\brho)-g^n(\phi_{s}^{s'}(\bz),\brho)
\right)\ d\bz+\g(s-s')R(\brho),
\end{align}
where $R$ satisfies the bound $|R(\brho)|\le c_{12}c_{11}^n/n!$,
for  positive constants $c_{11}$ and $c_{12}$ and $\g$ 
is an increasing non-negative function such that $\g(\th)/\th\to 0$ 
as $\th\to 0$.

To prove \eqref{combination1}, we use
 \eqref{leb2},  and the boundedness of $f$ and $\ell$, to assert
\be\label{splitting}
\int_{\overline{\Delta^n}} \cS'(\bq) g^n(\bz,\brho)\ d\bz
=\cS''(\brho)+\sum_{i=0}^{n}\cS'_i(\brho)+
(s-s')^2 R_0(\brho),
\ee
with the term $R_0$ satisfying $|R_0| \le c_{13} c_{11}^n/n!$,
for  constants $c_{11}$ and  $c_{13}$,
and the other terms given by
\begin{align*}
 \cS'_i(\brho)=&\int_{A_i} \exp \left(-\int_{s'}^{s} \mf{r}(u,\psi_{s'}^{u}\bq)\ du \right)
\left(G(\bq,s)-G(\psi_{s'}^{s}\bq,s)\right)g^n(\bz,\brho)\ d\bz, \\
\cS''(\brho)=&\int_{U} \exp \left(-\int_{s'}^{s} \mf{r}(u,\psi_{s'}^{u}\bq)\ du \right)\left(G(\bq,s)-G(\psi_{s'}^{s}\bq,s)\right)
g^n(\bz,\brho) \ d\bz .
\end{align*}
Observe that the replacement of the exponential 
\[
\exp \left(-\int_{s'}^{s} \mf{r}(u,\psi_{s'}^u\bq)\ du \right),
\]
with
\[
{\text{either}}\ \ \ \exp \left(-\int_{s'}^{s} \mf{r}(u,\bq)\ du \right)
\ \ \ {\text{or}} \ \ \ \exp \left(-\int_{s'}^{s} \mf{r}(u,\psi_{s'}^s\bq)\ du \right)
\]
 results in an error of size $O((n+2)(s-s')^2)$ by the Lipschitzness of the rates with respect to the space variable. Because of this, we can write,
\begin{align*}
\cS''(\brho)=&\int_{U} \exp \left(-\int_{s'}^{s} \mf{r}(u,\bq)\ du \right) G(\bq,s)g^n(\bz,\brho)\ d\bz  \\
&-\int_{U} \exp \left(-\int_{s'}^{s} \mf{r}(u,\psi_{s'}^{s}\bq)\ du\right)G(\psi_{s'}^{s}\bq,s)g^n(\bz,\brho)\ d\bz+(s-s')^2 R_1(\brho),
\end{align*}
with the term $R_1$ satisfying 
$|R_1(\brho)| \le c_{14}(n+2) c_{11}^n/n!$,
for  a constant $c_{14}$.
By a change of variables, 
\begin{align*}
\cS''(\brho)=&\int_{U} \exp \left(-\int_{s'}^{s} \mf{r}(u,\bq)\ du \right) G(\bq,s)g^n(\bz,\brho)\ d\bz  \\
&-\int_{\phi_{s'}^{s}(U)} \exp \left(-\int_{s'}^{s} \mf{r}(u,\bq)du \right)G(\bq,s)g^n(\phi_{s}^{s'}(\bz),\brho)\ d\bz+(s-s')^2 R_1(\brho),
\end{align*}
where  
\[
\phi_{s}^{s'}(\bz)=\left(\phi_{s'}^{s}\right)^{-1}(\bz)=\left(z_i-v_i(s-s')\right)_{i=1}^{n},\ \ \ \ 
\psi_{s}^{s'}(\bz)=\left(\phi_{s}^{s'}(\bz),\brho\right),\ \ \ \
\psi_{s}^{u}:=\psi_{s'}^{u}\circ \psi_{s}^{s'},
\]
are the reverse flows. 
This allows us to assert
\be\la{S''}
\cS''(\brho)=\cS''_1(\brho)+\cS''_2(\brho)+(s-s')^2R_1(\brho),
\ee
where 
\begin{align*}
\cS''_1(\brho):=&\int_{\phi_{s'}^{s}(U)} \exp \left(-\int_{s'}^{s} \mf{r}(u,\bq)\ du \right) G(\bq,s)
\left(g^n(\bz,\brho)-g^n(\phi_{s}^{s'}(\bz),\brho)\right)
 d\bz , \\
\cS''_2(\brho):=&\left(\int_U-\int_{\phi_{s'}^{s}(U)}\right) \exp \left(-\int_{s'}^{s} \mf{r}(u,\bq)\ du \right)G(\bq,s)g^n(\bz,\brho)\ d\bz.
\end{align*}
We now compare the set $U$ with its translate
$\phi_{s}^{s'}(U)$.
Observe,
\begin{align*}
U=&\Big\{ \bz \in \overline{\Delta^n} : \ z_1 > 
\max(a^-,a^--u_0(s-s')) , \ z_n<\min(a^+,a^++u_n(s-s')),\\
&\ \ \ \ \ \ \ \ \ \ \ \ \ \ \ \ \ \ \ \ \ \ \ \ \ \ \ \ \ \ \ \ \ \ \ \ \ \ \ \ z_{i+1}-z_i>\max(0,-u_i(s-s')) \ \text{ for }\  i \in \llbracket 1,n-1 \rrbracket \Big \},\\
\phi_{s'}^{s}(U) =&\Big\{ \bz \in \overline{\Delta^n} :\  z_1>
\max(a^-,a^-+u_0(s-s')),\ z_n<\min(a^+,a^+-u_n(s-s'))
\\ 
&\ \ \ \ \ \ \ \ \ \ \ \ \ \ \ \ \ \ \ \ \ \ \ \ \ \ \ \ \ \ \ \ \ \ \ \ \ \ \ \ \ \ \ z_{i+1}-z_i > \max(0,u_i(s-s')) \ \text{ for }\  i \in \llbracket 1,n-1 \rrbracket\Big \}.
\end{align*}
Hence  the symmetric difference of the sets 
$U$ and $\phi_{s'}^{s}(U)$ can be represented as
\be\la{S''1}
U\D\ \phi_{s'}^{s}(U)= \cup_{i=0}^n B_i.
\ee
Note that we can replace the exponential with $1$ at a cost of $O((n+1)(s-s'))$. From this and  \eqref{leb1} we deduce,
\[
\cS_2''(\brho)=\left(\int_{U}   -\int_{\phi_{s'}^{s}(U)} \right)G(\bq,s)g^n(\bz,\brho)\ d\bz+(s-s')^2R_2(\brho),
\]
with the term $R_2$ satisfying 
$|R_2(\brho)| \le c_{15}(n+1) c_{11}^n/n!$,
for  a constant $c_{15}$. 
 This, \eqref{S''1} and \eqref{leb2} allow us to ignore the overlaps of the sets $B_i, \ i\in \llbracket 1,n \rrbracket$, so that we can
write
\be\la{S''3}
\cS_2''(\brho)=
\sum_{i=0}^{n} \text{sign}(u_i) \int_{B_i}G(\bq,s)
g^n(\bz,\brho)\ d\bz+(s-s')^2R_3(\brho),
\ee
with the term $R_3$ satisfying 
$|R_3(\brho)| \le c_{16}(n+1) c_{11}^n/n!$,
for  a constant $c_{16}$.

We now turn our attention to $\cS''_1$. By the Taylor expansion
 we have that
\be\la{taylor}
\frac{g^n(\phi_{s}^{s'}(\bz),\brho)}{g^n(\bz,\brho)}=1+(s'-s)\sum_{i=1}^{n} v_i \frac{g^n_{z_i}(\bz,\brho)}{g^n(\bz,\brho)}+R_4(\bz,\brho,
s-s')\g(s-s'),
\ee
where $|R_4|\le c_{17}(n+2)$, for a constant $c_{12},$
and $\g$ is a function satisfying $\g(\th)/\th\to 0$ as $\th\to 0$.
In particular, there exists a constant $c_{18}$ such that 
\be\la{gbound}
\frac{|g^n(\phi_{s}^{s'}(\bz),\brho)-g^n(\bz,\brho)|}{g^n(\bz,\brho)}
\le c_{18}(n+2)(s-s').
\ee
 This allows us to make two 
changes in $\cS''_1$ at a cost of a constant multiple of $\g(s-s')$, namely
 replacing the exponential with $1$, and
 replacing the set $\phi_{s'}^s(U)$ with $U$.
As a result, 
\be\label{S''4}
\cS_1''(\brho)= 
\int_{U}G(\bq,s)\left(g^n(\bz,\brho)-g^n(\phi_{s}^{s'}(\bz),\brho)\right)\ d\bz
+R_5(\brho)\g(s-s'),
\ee
with again  the term $R_5$ satisfying 
$|R_5(\brho)| \le c_{20}(n+1) c_{11}^n/n!$,
for  a constant $c_{20}$.
Coming back to the second
 term $\cS'_i$ in \eqref{splitting}, we see use \eqref{leb1} to replace the exponential term of the integrand with 1:
\begin{align}\label{S''5}
\cS'_i(\brho)&= \int_{A_i} \left(G(\bq,s)-G(\psi_{s'}^{s}\bq,s))g^n(\bz,\brho\right)\ d\bz+R_6(\brho)\g(s-s'),
\end{align}
with the term $R_6$ satisfying 
$|R_6(\brho)| \le c_{21}(n+1) c_{11}^n/n!$,
for  a constant $c_{21}$.  Now, notice that for $i \in \llbracket 0,n \rrbracket$ 
if $u_i>0$, then $A_i= \emptyset$, otherwise if $u_i<0$, then $A_i=B_i$.  From this,  \eqref{splitting}, \eqref{S''},
 \eqref{S''3}, \eqref{S''4}, and \eqref{S''5}, we obtain
\eqref{combination1}

\ms\noi
{\em (Step 4)} In this step, we use \eqref{combination1} to show
\begin{align}
\int_{\Omega^n}  \cS'(\bq)  \  \mu^n(d\bq,s)
=&-(s-s')\left(\G_n^d(s)+\G_n^t(s)\right)+\hat R^n(s,s')\g(s-s'),\label{deterministic}
\end{align}
with $\hat R^n$ satisfying $|\hat R^n| \le c_{22}(n+1)^2 c_{11}^n/n!$,
for  a constant $c_{22}$.  Put
\[
\hat B_i=B_i\setminus\cup_{j\neq i}B_j,
\]
so that the sets $(\hat B_i:\ i\in\llbracket 0,n\rrbracket)$ are mutually disjoint. The bound \eqref{leb2} allows us to replace
$B_i$ with $\hat B_i$ in \eqref{combination1} at a small cost:
\be\label{combination}
\int_{\overline{\Delta^n}} \cS'(\bq) g^n(\bz,\brho)\ d\bz=
\sum_{i=0}^{n}\left(\widehat \cS_i^+(\brho)-\widehat \cS_i^-(\brho)\right)
+\cS'''(\brho)+R_8(\brho)\g(s-s'),
\ee
where $R_7$ satisfies $|R_7(\brho)| \le c_{23}(n+1) c_{11}^n/n!$, for a constant $c_{23}$, and
\begin{align*}
\widehat \cS_i^+(\brho)=&
 \ind_{\{u_i>0\}} \int_{\hat B_i} G(\bq,s)g^n(\bz,\brho)\ d\bz,\\
\widehat \cS_i^-(\brho)=&
 \ind_{\{ u_i<0\}} \int_{\hat B_i} G(\psi_{s'}^{s}\bq,s)g^n(\bz,\brho)\ d\bz,\\
\cS'''(\brho)=&\int_{U}G(\bq,s)\left(g^n(\bz,\brho)-g^n(\phi_{s}^{s'}(\bz),\brho)\right)\ d\bz.
\end{align*}
We establish \eqref{combination} by proving
\begin{align}\la{hatS}
&\int\sum_{i=0}^{n}\left(\widehat \cS_i^+(\brho)
-\widehat \cS_i^-(\brho)\right)\ \b(d\brho)=-(s-s')\G_n^d(s)
+\hat R_{0}(s,s')(s-s')^2,\\
&\int \cS'''(\brho)\ \b(d\brho)
=-(s-s')\G_n^t(s)
+\hat R_{1}(s,s')(s-s')^2,\la{hatS'}
\end{align}
with $|\hat R_0| ,|\hat R_1|\le c_{24}(n+1) c_{11}^n/n!$, for a constant $c_{24}$.

To prove \eqref{hatS}, first observe that if 
$\sigma:=\sigma(\bq,s')$ denotes the first collision time of the deterministic flow starting at time $s'$ at $\bq=(\bz,\brho)$,
and $\s'=\s\wedge T_1$, 
then for sure $\s=\s'<s$ provided that $u_i<0$, $\bz\in\hat B_i$, and no stochastic jump  occurs in the interval $[s',s]$
(equivalently, $T_1> s$).
From this, the strong Markov property, \eqref{r-bound},
and Lemma \ref{lipschitz}  we deduce 
\begin{align*}
G(\psi_{s'}^{s}\bq,s)&=\bE[G(\Psi_{s}^{t}\psi_{s'}^{s}\bq)]
=\bE[G(\Psi_{s}^{t}\psi_{s'}^{s}\bq)\ \ind_{\cE^c}]+O(s-s')\\
&=\bE[G(\Psi_{\sigma'}^{t}\psi_{s'}^{\sigma'}\bq)\ \ind_{\cE^c}]+O(s-s')=\bE[G(\Psi_{\sigma'}^{t}\psi_{s'}^{\sigma'}\bq)]+O(s-s')\\
&=\bE[G(\psi_{s'}^{\sigma'}\bq,\sigma')]+O(s-s')
=G(\psi_{s'}^{\sigma'}\bq,s)+O(s-s'),
\end{align*}
 provided that $u_i<0$, and $\bz\in\hat B_i$.
On $\hat B_i$, and when $u_i<0$, we have that
\[
G(\psi_{s'}^{\sigma}\bq,s)=G(E_{i}^{\rho^i}\bq_i,s),
\]
for $i \in \llbracket 0,n \rrbracket$, with $\bq_i=(\bz_i,\brho_i)$ the configuration $\bq$ with the particle $i$ removed, and
where $E^{\rho^0}_0$ and $E^{\rho^n}_{n}$ should be understood as 
$E^{\rho^0}_{-}$ and $E^{\rho^n}_{+}$.  
After replacing $G(\psi_{s'}^{s}\bq,s)$ with $G(E_{i}^{\rho^i}\bq_i,s)$
at a cost of $O(s-s')$, we replace back the set $\hat B_i$ with $B_i$.
This allows us to write
\be\la{hatS2}
\widehat\cS^-_i(\brho)=(s-s')\cT^-_i(\brho)+R_9(\rho)(s-s')^2,
\ee
where $R_9$ satisfies $|R_9| \le c_{25}(n+1) c_{11}^n/n!$, for a constant $c_{25}$, and
\begin{align*}
 \cT^-_0(\brho)=& u_0^-\frac{\ell(a^-,s,\rho^0)f(a^-,s,\rho^0,\rho^1)}{\ell(a^-,s,\rho^1)} \int_{{\Delta^{n-1}}} G(E_{-}^{\rho^0}\bq_0,s)g^{n-1}(\bz_0,\brho_0)\ d\bz_0,\\
  \cT^-_i(\brho)=&
 u_i^- \frac{f(z_i,s,\rho^{i-1},\rho^i)f(z_i,s,\rho^i,\rho^{i+1})}{f(z_i,s,\rho^{i-1},\rho^{i+1})}\int_{{\Delta^{n-1}}} G(E_i^{\rho^i}\bq_i,s)
g^{n-1}(\bz_i,\brho_i)\ d\bz_i,\\
 \cT^-_n(\brho)=&
 u_n^- f(a^+,s,\rho^{n-1},\rho^n)\int_{{\Delta^{n-1}}} G(E_i^{\rho^n}\bq_n,s)
g^{n-1}(\bz_n,\brho_n)\ d\bz_n,
\end{align*}
for $i\in\llbracket 1,n-1 \rrbracket$,
where $u^-=\ind_{\{ u<0\}}|u|$.

The terms $\widehat\cS_i^+$ can be treated likewise:  Fix $i \in \llbracket 0,n \rrbracket$, and  define the time $\sigma':=s-\frac{z_{i+1}-z_i}{u_i}$, with the convention that $z_0=a^-$ and $z_{n+1}=a^+$. When $\bz \in \hat B_i$ and $u_i>0$,  there would be no collision in the interval
$[\s',s]$, and 
the reverse flow $\psi_{s}^{\sigma'}\bq$ is well-defined
for $\bq=(\bz,\brho)$. We then 
define $\tilde \bq$ to be $\tilde{\bq}=\psi_{s}^{\sigma'}\bq$. 
By similar arguments we have 
\[
G(\bq,s)=G(\tilde{\bq},\sigma')+O(s-s')=G(\tilde{\bq},s)+O(s-s').
\]
 Again, we see that $G(\tilde{\bq},s)=G(E_{i}^{\rho^i}\bq_i,s),$
and hence we get the analog of \eqref{hatS2}, namely
\be\la{hatS3}
\widehat\cS^+_i(\brho)=(s-s')\cT^+_i(\brho)+R_{10}(\rho)(s-s')^2,
\ee
where $R_{10}$ satisfies $|R_{10}| \le c_{25}(n+1) c_{11}^n/n!$, 
 and the expression for $\cT_i^+$ is the same as $\cT_i^-$, except
that $u_i^-$ is replaced with $u_i^+$.
We integrate both sides of \eqref{hatS2} and \eqref{hatS3}
 against $\b$, and take the difference to arrive at \eqref{hatS}.
After integrating out \eqref{hatS} and \eqref{S'''}
with respect to $\brho$ and by relabeling 
$\rho^i$ as $\rho^*$ and $\rho^j$ for $j>i+1$ by $\rho^{j-1}$, we 
arrive at \eqref{deterministic}.

 We now focus on \eqref{hatS'}.
From \eqref{taylor} and the straightforward computation 
\[
\frac{g^n_{z_i}(\bz,\brho)}{g^n(\bz,\brho)}=\frac{f_z(z_i,s,\rho^{i-1},\rho^i)}{f(z_i,s,\rho^{i-1},\rho^i)}+\lambda(z_i,s,\rho^i)-\lambda(z_i,s,\rho^{i-1}),
\]
 we deduce 
\be\la{S'''}
\cS'''(\brho)
=-(s-s')\sum_{i=1}^{n}\cS'''_i(\brho)+R_{10}(\brho)\g(s-s'),
\ee
where $R_{10}$ satisfies 
$|R_{10}(\brho)| \le c_{26}(n+1) c_{11}^n/n!$,
 for a constant $c_{26}$, and
\begin{align*}
\cS'''_i(\brho)=& \int_{U} G(\bq,s)[\rho^{i-1},\rho^i]
\left(\frac{f_z(z_i,s,\rho^{i-1},\rho^i)}{f(z_i,s,\rho^{i-1},\rho^i)}+\lambda(z_i,s,\rho^i)-\lambda(z_i,s,\rho^{i-1})\right)g^n(\bz,\brho)\
d\bz.
\end{align*}
Note that \eqref{leb1} provides us a bound on the Lebesgue measure  of the set $U^c$. This bound and \eqref{gbound} allow us to replace the domain of integration from $U$ to the whole simplex
${\Delta}^n$ at a cost of replacing $R_{10}$
with $R_{11}$, that satisfies a similar bound. After such a replacement, we integrate both sides with respect to $\b$ to deduce \eqref{hatS'}.

\ms\noi
{\em(Final Step)} Note that our error terms are bounded by  constant multiples of $\bar c_n\g(s-s')$, with 
$\bar c_n=(n+1)^3c_{11}^n/n!$. Since
$\sum_n \bar c_n<\i$, and $\g(\th)/\th\to 0$ as $\th\to 0$, these error terms can be ignored as we calculate the limit in \eqref{lim}.
On account of this, \eqref{decompos}, and \eqref{deterministic}, 
it remains to verify
\begin{align}\label{stochastic}
\int_{\Omega^n} \cS^\pm(\bq)\
\mu^{n}(d\bq,s)=&-(s-s')\G_n^\pm(s)
+(s-s')^2R^\pm(s),\\
\sum_{i=1}^n\int_{\Omega^n} \cS_i(\bq)\
\mu^{n}(d\bq,s)=&-(s-s')\G_n^{\mf{f}}(s)
+(s-s')^2R^{\mf{f}}(s),\nonumber
\end{align}
with $R^\pm,R^f$ satisfying $|R^\pm|,|R^{\mf{f}}|\le c_{27}(n+1) c_{11}^n/n!$,
 for a constant $c_{27}$. 

We only verify \eqref{stochastic} in the case of $\cS^-$, as the other cases can be treated in the same fashion. For this, first observe
that because of \eqref{hato}, we may write
\be\la{S-}
\int_{\Omega^n} \cS^-(\bq)\ \mu^{n}(d\bq,s)=\cS^-_1-\cS^-_2+
R_n^-(s) (s-s')^2,
\ee
where $R_0^-$ satisfies $|R_0^-|\le c_{28}(n+1) c_{11}^n/n!$,
for a constant $c_{28}$, and
\begin{align*}
\cS^-_1=&\int_{\Omega^n} \int_{s'}^{s} \int_{L(\rho^0)} \mf{c}_{-}(\th,\rho^*,\rho^0)G(\bq,s)\ \b(d\rho^*)d\th\  \mu^{n}(d\bq,s),\\
\cS^-_2=&\int_{\Omega^n} \int_{s'}^{s} \int_{L(\rho^0)} \mf{c}_{-}(\th,\rho^*,\rho^0)G(E_0^{\rho^*}\psi_{s'}^\th\bq,s)\ \1(\bq\notin  \hat\O)\ \b(d\rho^*)d\th\ \mu^{n}(d\bq,s).
\end{align*}
On the other-hand, for $\bq\notin \hat\O$, the 
  the flow $\psi_{s'}^{t}$ experiences no collision, and is just a translation that preserves the volume. Again focusing on the spatial integration first, we certainly have
\[
\int_{\D^n} G(E_0^{\rho^*}\psi_{s'}^\th\bq,s) \1(\bq\notin  \hat\O)
 g^{n}(\bz,\brho)\ d\bz
=\int_{\D^n} G(E_0^{\rho^*}\bq,s) \1(\psi_{s'}^\th\bq\notin  \hat\O)
 g^{n}(\phi_\th^{s'}\bz,\brho)\ d\bz.
\]
We can then use \eqref{gbound} (with $s$ replaced with $\th$)
to replace $g^{n}(\phi_\th^{s'}\bz,\brho)$ with $g^{n}(\bz,\brho)$
at a cost that is bounded by a constant multiple of $\bar c_n (s-s')^2$. 
As in \eqref{hato}, we can readily show that at a cost of $O((s-s')^2)$,
we can now drop $\1(\psi_{s'}^\th\bq\notin  \hat\O)$. From all this we conclude 
\[
\cS^-_2=\int_{\Omega^n} \int_{s'}^{s} \int_{L(\rho^0)} \mf{c}_{-}(\th,\rho^*,\rho^0)G(E_0^{\rho^*}\bq,s)\ \1(\bq\notin  \hat\O)\ \b(d\rho^*)d\th\ \mu^{n}(d\bq,s)+\hat R_n^-(s) (s-s')^2,
\]
with $|\hat R_n^-|\le c_{29}\bar c_n$, for a constant $c_{29}$.
From this and \eqref{S-} we deduce \eqref{stochastic}, completing the
proof of our theorem.
\qed

\subsection{Proof of Theorem \ref{main}}

\noi
 {\em(Step 1)} We need to check that
\[
\lim_{s' \to s} \frac{1}{s-s'}\left(\int G(\bq,s)\mu(d\bq,s)-\int G(\bq,s')
\ \mu(d\bq,s')\right)=0.
\]
We can certainly write
\be\la{three}
 \int G(\bq,s)\ 
\mu(d\bq,s)-\int G(\bq,s')\mu(d\bq,s')=\cX_1(s',s)+\cX_2(s',s)
-\cX_3(s',s),
\ee
where
\begin{align*}
\cX_1(s',s)=& \int (G(\bq,s)-G(\bq,s'))\ \mu(d\bq,s),\\
\cX_2(s',s)=&\int G(\bq,s) \ (\mu(d\bq,s)-\mu(d\bq,s')),\\
\cX_3(s',s)=&\int((G(\bq,s)-G(\bq,s'))\ (\mu(d\bq,s)-\mu(d\bq,s')).
\end{align*}
We work out $\cX_2$ by differentiating $\mu$ with respect to the time $s$. Evidently
\be\la{X2}
\lim_{s' \uparrow s} (s-s')^{-1} \cX_2(s',s)=\int G(\bq,s) \ \dot\mu(d\bq,s),
\ee
where $\dot \mu$ represents the $s$-derivative of $\mu$. We may write  $\dot\mu^{n}=X^{n}\mu^{n}$, where
\[
X^n(\bq,s)=\frac{\dot{\ell}(a^-,s,\rho^0)}{\ell(a^-,s,\rho^0)}+\sum_{i=1}^{n} \frac{f_s(z_i,s,\rho^{i-1},\rho^i)}{f(z_i,s,\rho^{i-1},\rho^i)}-\sum_{i=0}^{n} \int_{z_i}^{z_{i+1}} \lambda_s(z,s,\rho^i)dz .
\]
Now using the kinetic equations verified by both $\ell$ and $f$, let us find an explicit expression of the Radon-Nidokym derivative $X^{n}$.
For the first term we have 
\[
\frac{\dot{\ell}(a^-,s,\rho^0)}{\ell(a^-,s,\rho^0)}=\int_{L(\rho^0)} [\rho^*,\rho^0]\ \frac{\ell(a^-,s,\rho^*)f(a^-,s,\rho^*,\rho^0)}{\ell(a^-,s,\rho^0)}\ \b(d\rho^*)-A(a^-,s,\rho^0).
\]
For the second term we get for $i \in \llbracket 1,n \rrbracket$,
\begin{align*}
\frac{f_s(z_i,s,\rho^{i-1},\rho^i)}{f(z_i,s,\rho^{i-1},\rho^i)}=&[\rho^{i-1},\rho^i]\ \frac{f_z(z_i,s,\rho^{i-1},\rho^i)}{f(z_i,s,\rho^{i-1},\rho^i)}\\
&+\int_{D(\rho^{i-1},\rho^i)} [\rho^{i-1},\rho^*,\rho^i]\
\frac{f(z_i,s,\rho^{i-1},\rho^*)f(z_i,s,\rho^*,\rho^i)}{f(z_i,s,\rho^{i-1},\rho^i)}\ \b(d\rho^*)\\
&+[\rho^{i-1},\rho^i](\lambda(z_i,s,\rho^i)-\lambda(z_i,s,\rho^{i-1}))-(A(z_i,s,\rho^{i})-A(z_i,s,\rho^{i-1})).
\end{align*}
For the last term, let us first show $\l_s=A_z$.
To see this, observe
\begin{align*}
\lambda_s(z,s,\rho)=&\int_{R(\rho)} f_s(z,s,\rho,\rho^{+})\ \b(d\rho^{+})
=\int_{R(\rho)} [\rho,\rho^{+}]f_z(z,s,\rho,\rho^{+})\ \b(d\rho^{+})
\\
&+\iint _{\{\rho \prec \rho^* \prec \rho\}} [\rho,\rho^*,\rho^{+}]f(z,s,\rho,\rho^*)f(z,s,\rho^*,\rho^{+})\ \b(d\rho^*)\ \b(d\rho^{+})\\
&+\int_{R(\rho)} [\rho,\rho^{+}](\lambda(z,s,\rho^{+})-\lambda(z,s,\rho))+(A(z,s,\rho)-A(z,s,\rho^{+}))f(z,s,\rho,\rho^{+})\ \b(d\rho^{+}).
\end{align*}
By integrating out $\rho^+$, we get that the double integral equals to
\begin{align*}
\int_{R(\rho)} \left(f(z,s,\rho,\rho^*)A(z,s,\rho^*)-[\rho,\rho^*]
f(z,s,\rho,\rho^*)\lambda(z,s,\rho^*)\right)\ \b(d\rho^*).
\end{align*}
However $\rho^*$ and $\rho^+$ are both just dummy variables in our integrals, so summing over the terms, we get 
\begin{align*}
\lambda_s(z,s,\rho)=&\int_{R(\rho)} [\rho,\rho^+]f_z(z,s,\rho,\rho^+)\ \b(d\rho^+)-\lambda(z,s,\rho)\int_{R(\rho)} [\rho,\rho^+]f(z,s,\rho,\rho^+)\ \b(d\rho^+)\\
&+A(z,s,\rho)\int_{R(\rho)}f(z,s,\rho,\rho^+)\ \b(d\rho^+)=
\int_{R(\rho)} [\rho,\rho^+]f_z(z,s,\rho,\rho^+)\ \b(d\rho^+),
\end{align*}
confirming our claim $\l_s=A_z$. As a result,
\begin{align*}
\int_{z_i}^{z_{i+1}} \lambda_s(s,z,\rho^i)\
dz&=\int_{R(\rho)} [\rho^i,\rho^+] \ \b(d\rho^+) 
\int_{z_i}^{z_{i+1}} f_z(z,s,\rho^i,\rho^+)\ dz\\
&=A(z_{i+1},s,\rho^i)-A(z_i,s,\rho^i).
\end{align*}
Summing over everything, we get that
\begin{align}\nonumber
X^{n}(\bq)=&-A(a^+,s,\rho^n)+\int_{L(\rho^0)} [\rho^*,\rho^0]\frac{\ell(a^-,s,\rho^*)f(a^-,s,\rho^*,\rho^0)}{\ell(a^-,s,\rho^0)}\ \b(d\rho^*)\\ \nonumber
&+\sum_{i=1}^{n} [\rho^{i-1},\rho^i]\left(\frac{f_z(z_i,s,\rho^{i-1},\rho^i)}{f(z_i,s,\rho^{i-1},\rho^i)}+\lambda(z_i,s,\rho^i)-\lambda(z_i,s,\rho^{i-1})\right)\\
&+\sum_{i=1}^{n} \int_{D(\rho^{i-1},\rho^i)} [\rho^{i-1},\rho^*,\rho^i]\frac{f(z_i,s,\rho^{i-1},\rho^*)f(z_i,s,\rho^*,\rho^i)}{f(z_i,s,\rho^{i-1},\rho^i)}\ \b(d\rho^*).\la{Xn}
\end{align}

\ms\noi
{\em(Step 2)}
From \eqref{Xn}, we learn that
 there is  a uniform constant $C>0$ such that 
$|X^n(\bq)| \le C(n+1)$. Hence using this observation and Theorem \ref{lipschitz} we get 
\[
\left|\int_{\Omega} (G(\bq,s)-G(\bq,s'))(\mu(d\bq,s)-\mu(d\bq,s')) \right| \le C(s-s')^2 \int_{\Omega} (\bn(\bq)+2)^3 \ \mu(d\bq,s) .
\]
As a result,
\[
\lim_{s' \uparrow s} (s-s')^{-1} \cX_3(s',s)=0.
\]
Because of this, \eqref{three}, and \eqref{X2} we are done if we can 
show
\be\la{last}
\lim_{s' \to s} (s-s')^{-1}\cX_1(s',s)=-
\sum_{n=0}^{\i} \int_{\Omega^n} X^{n}(\bq)G(\bq,s)
\ \mu^{n}(d\bq,s).
\ee
By \eqref{Xn},
\[
\int_{\Omega^n} X^{n}(\bq)G(\bq,s)\
\mu^{n}(d\bq,s)=\Lambda_n^{-}(s)+\G_n^{t}(s)+\Lambda_n^{\mf{f}}(s)+\Lambda_n^{+}(s),
\]
where
\begin{align*}
\Lambda_n^{-}(s)=&\int_{\Omega^n}\int_{L(\rho^0)} [\rho^*,\rho^0]\frac{f(a^-,s,\rho^*,\rho^0)\ell(a^-,s,\rho^*)}{\ell(a^-,s,\rho^0)}G(\bq,s)\ \b(d\rho^*)\ \mu^{n}(d\bq,s),\\
\Lambda_n^{+}(s)=&-\int_{\Omega^n} \int_{R(\rho^n)} [\rho^n,\rho^*]f(a^+,s,\rho^n,\rho^*)G(\bq,s)\ \b(d\rho^*)\
 \mu^{n}(d\bq,s),\\
\Lambda_n^{\mf{f}}(s)=&\int_{\Omega_L^n} \sum_{i=1}^{n} \int_{D(\rho^{i-1},\rho^i)}[\rho^{i-1},\rho^*,\rho^i] \frac{f(z_i,s,\rho^{i-1},\rho^*)f(z_i,s,\rho^*,\rho^i)}{f(z_i,s,\rho^{i-1},\rho^i)}
G(\bq,s)\ \b(d\rho^*)\mu^{n}(d\bq,s),
\end{align*}
where of course $\Lambda_0^{\mf{f}}(s)=\G_0^{t}(s)=0$. 
On account of Theorem 3.3, 
we will be done if we show that for any $n \ge 0$, the following equalities hold:
\begin{align}\la{lamgam}
&\Lambda_n^{-}(s)=\Gamma_n^{-}(s)
+\Gamma_{n+1}^{d,-}(s),\\ \nonumber
&\Lambda_n^{+}(s)=\Gamma_n^{+}(s)
+\Gamma_{n+1}^{d,+}(s),\\ \nonumber
&\Lambda_n^{\mf{f}}(s)=\Gamma_n^{\mf{f}}(s)+\Gamma_{n+1}^{d,\mf{f}}(s).
\end{align}
We will verify this only for the first equality, as the others are done in a similar fashion. By the construction of our Markov process, we know that when $[\rho^*,\rho^0]>0$, then $G(E_{-}^{\rho^*}\bq,s)=G(\bq,s)$ (as in this case the particle at $z=a^-$ corresponds to an exit of the interval $[a^-,a^+]$ and becomes irrelevant instantaneously), thus by splitting the two cases whether $[\rho^*,\rho^0]$ is negative or nonnegative we get
\begin{align*}
\Gamma_{n+1}^{d,-}(s)=&\int_{\Omega^{n}} \int_{L(\rho^0)}  [\rho^*,\rho^0]^{+}\frac{\ell(a^-,s,\rho^*)f(a^-,s,\rho^*,\rho^0)}{\ell(a^-,s,\rho^0)}G(\bq,s)\ \b(d\rho^*)\mu^{n}(d\bq,s)\\
&-\int_{\Omega^n} \int_{L(\rho^0)} \mf{c}_{-}(s,\rho^*,\rho^0)G(E_{-}^{\rho^*}\bq,s)\ \b(d\rho^*)\mu^{n}(d\bq,s).
\end{align*}
This immediately implies the first identity in \eqref{lamgam}.

For the terms corresponding to the fragmentation i.e.,
 $\Lambda_n^{\mf{f}}$ and $\Gamma_n^{\mf{f}}$, we use now the observation that when $[\rho^{i-1},\rho^*,\rho^{i}]>0$, 
then $G(E_{i}^{\rho^*}\bq,s)=G(\bq,s)$ to get similarly the desired equality.
\qed

\subsection{Proof of the genericity of the tessellation $\bX_\L$}
 
We now state and prove a Proposition that guarantees the genericity of the tessellation $\bX_\L$, which is induced by the process
$\bq(t)=(q_i(t):\ 1\le i\le \bn(t))$, where $\bn(t):=\bn(\bq(t))$.

\bp\la{pro3.1} Under the assumptions of Theorem 1.1, the probability of the occurrence of the event $z_{i-1}(t)=z_i(t)=z_{i+1}(t)$ is zero.
In words, no three particles arrive at the same location almost surely.
\ep

\ms\noi
{\bf Proof} {\em(Step 1)} The main idea is that since particles have
uniformly bounded velocities, the probability of the occurrence of two collisions in a time interval of size $\d$ is of order $O(\d^2)$. The reason for this is that if we trace back the colliding particles to the boundary of the box, we will have a configuration in which two pairs of particles have distances of order $O(\d)$. 

Let us write $\s(t)$ for the smallest time $\s>t$ such that
$z_{i-1}(\s)=z_i(\s)=z_{i+1}(\s)$ for some index 
$i\in \{2,\dots, \bn(t)-1\}$. 
We first claim that there exists a constant $c_1$ such that
\be\la{eq3.57}
\bP\left(\s(t)\in (t,t+\d)\right)\le c_1 \d^2,
\ee
 for every $\d>0$. Once this is established, we can then choose
$\d=T/n$, and argue
\[
\bP\left(\s(0)\in [0,T]\right)\le n 
\sup_i\bP\left(\s(t_i)\in [t_i,t_{i+1}]\right)\le c_1 T^2 n^{-1},
\]
for $t_i=iT/n,\ i=0,\dots, n-1$. We then send $n\to\i$ to deduce
that there is no triple collision almost surely.

It remains to prove \eqref{eq3.57}. Let us write $\cE_r(t_1,t_2)$ 
for the event
that at least $r$ stochastic jumps occur in the interval $[t_1,t_2]$.
We claim that there exists a constant $c_2=c_2(T)$ such that
\be\la{eq3.58}
\bP\left(\cE_2(t,t+\d)\right)\le c_2 \d^2,
\ee
for every $t\in[0,T]$, and $\d\in (0,1)$. 
This is an immediate consequence of the strong Markov property of 
the process $\bq(t)$, and the bound \eqref{cE-bound}: Indeed if $\s_1$ and
$\s_2$ denote the times of the first and second stochastic jumps after time $t$, then 
\begin{align*}
\bP\left(\cE_2(t,t+\d)\right)&\le \bE \ \1_{\cE_1(0,t+\d)}\ \bP^{\bq(\s_1)}
\left(\cE_1(\s_1,t+\d)\right)\\
&\le c_3\d  \ \bE \left( \bn(\bq(\s_1))+2\right)^2\ \1_{\cE_1(t,t+\d)} \\
&\le 
 c_3\d  \bE \left (\bn(\bq(t))+3\right)^2\ \1_{\cE_1(t,t+\d)}  \\
&\le c_4\d^2 \bE \left (\bn(\bq(t))+3\right)^4= c_5\d^2,
\end{align*}
for some constants $c_3, c_4$, and $c_5$. 
Because of \eqref{eq3.58}, the bound
\eqref{eq3.57} would follow if we can show 
\be\la{eq3.59}
\bP\left(\s(t)\in (t,t+\d),\ \s_2> t+\d \right)\le c_6 \d^2,
\ee
for some constant $c_6$. 

\ms\noi
{\em(Step 2)}
It remains to establish \eqref{eq3.59}. Let us write $\th_1$ and $\th_2$
for the first and the second times of particle collisions after $t$.
By convention, $\th_1=\th_2=\th$ when at time $\th$ there is a double 
collisions (which includes the case $z_{i-1}(t)=z_i(t)=z_{i+1}(t)$).
We claim 
\be\la{eq3.60}
\bP\left(\th_1,\th_2\in (t,t+\d),\ \th_2\le \s_1 \right)\le c_7 \d^2.
\ee
To see this, observe that the particle $z_i$ travels with the velocity
$v_i=-[\rho^{i-1},\rho^i]\in[-V_\i,V_\i]$. Hence, if there is a collision between
$z_i$ and $z_{i+1}$, and between $z_j$ and $z_{j+1}$ at time $\th_1$, then $|z_i(t)-z_{i+1}(t)|=O(\d)$ and $|z_j(t)-z_{j+1}(t)|=O(\d)$.
Since $i\neq j$, the probability of such event is of order
$O(\d^2)$ because the law of $\bq(t)$ has a bounded density with respect to the Lebesgue
measure. In summary,
\[
 \bP\left(\th_1,\th_2\in (t,t+\d),\ \th_2\le \s_1 \right)\le c_8 
\d^2\bE (\bn(\bq(t)))^2\le c_7\d^2,
\]
for constants $c_7$ and $c_8$, proving \eqref{eq3.60}.

\ms\noi
{\em(Final Step)} On account of \eqref{eq3.60}, the bound \eqref{eq3.59} would follow if we can show
\be\la{eq3.63}
 \bP\left( \s_1\le \th_1=\th_2=\s(t)\le t+\d\right)\le c_9 \d^2,
\ee
for a constant $c_9$. At $\s_1$ a new particle is created either at boundary point, or as a result of a fragmentation. We only treat the latter  because the former can be treated likewise. Let us assume that 
the particle $(z_i,\rho^i)$ is replaced with two particles 
$(z_{i,1},\rho^*) $ and $(z_{i,2},\rho^i) $. At the time $\s(t)$
a double collision occurs. If none of $z_{i,1}$ or $z_{i,2}$ are involved
in this collision, then can be treated as in {\em Step 2}. If, for example
the particle $z_{i,2}$ is to the right of $z_{i,1}$, and 
is involved in the double collision
at $\th:=\th_1=\th_2=\s(t)$, 
then $z_{i,2}(\th)=z_{i+1}(\th)=z_{i+2}(\th)$. In spite of a change 
of velocity of the $i$-th particle, we still must have that 
$|z_i(t)-z_{i+1}(t)|,|z_{i+1}(t)-z_{i+2}(t)|=O(\d)$, which results in a 
bound of order $O(\d^2)$ for the probability of such event. This completes the proof of \eqref{eq3.63}.
\qed

\bs

\section{Proof of Theorem 1.2}

As we stated in the Introduction, our assumption on the support of $f$ allows us to deduce Theorem 1.2 from Theorem 1.1. To explain our strategy, observe that by Theorem 1.1 we already know that the process
$x_1\mapsto \rho(x_1,x_2)$ is a Markov process for every $x_2\in[t_0,t_1]$. Under the assumptions of Theorem 1.2, no new particle is created on the left side
of $\L$, and the process $x_2\mapsto\rho(a^+,x_2)$ is an
independent Markov jump process with the jump rate density
$[\rho^-,\rho^+]f\big((a^+,x_2),\rho^-,\rho^+\big)$.
Though its initial state $\rho(a^+,t_0)$ depends on the dynamics
of the lower side of $\L$.  
 We may deduce Theorem 1.2 from 
Theorem 1.1 by interchanging
$x_1$ with $x_2$, and reversing direction on both axes. We explain the 
consequences of  
such operations on our variables in steps {\bf(1)} and {\bf(2)} below. 
In {\bf(1)} we verify the compatibility of the forward
equations \eqref{kin-x} and \eqref{kin-t} that are satisfied by the marginal $\ell$.  We use $\ell$ in {\bf(2)} to give a recipe for
the jump rates of the reversed processes, and  see how the time 
reversal and variables swap operations are compatible
with the kinetic equation \eqref{eq1.12}.

\bs\noi
{\bf(1)} We first address the effect of a time reversal on our Markov jump processes on the lower and the right sides of $\L$. 
Given a kernel $h(x,\rho^-,d\rho^+)$, define the linear operator
$\cL(x,h)$ by
\[
\big(\cL(x,h) F\big)(\rho^-):=\int \big(F(\rho^+)-F(\rho^-)\big)\ h(x,\rho^-,d\rho^+).
\]
If $\ell(x,d\rho)$ is the law of $\rho(x)$ with respect to the measure 
$\nu^{f,\L}$, then it satisfies the forward equation associated with the operator $\cL (x,f)=\cL(x,{f^1})$:
\be\la{eq5.1}
\ell_{x_1}=\cL(x,{f})^*\ \ell=\ell*f-A(f)\ \ell.
\ee
where $\cL(x,{f})^*$ denotes the adjoint of the operator 
$\cL(x,{f})$, and 
\[
(\ell*h)(x,d\rho^+):=\int \ell(x,d\rho^-)\ h(x,\rho^-,d\rho^+).
\]
The equation \eqref{eq5.1}
is an immediate consequence of Theorem 1.1. This and the kinetic equation \eqref{eq1.12} imply a similar equation for the second partial derivative, namely

\bp\la{pro1.1} Assume that $f$ and $\ell$ are bounded functions,
 $f$ is $C^1$ in $x$-variable and that
$\ell$ is $C^2$ in $x$-variable.
Also assume that $f$ satisfies \eqref{eq1.12}, $\ell$ satisfies \eqref{eq5.1}, and that the equation
\be\la{eq5.2}
\ell_{x_2}=\cL(x,{f^2})^*\  \ell=\ell*f^2-A(f^2)\ \ell,
\ee
holds when $x_1=a^+$.  Then \eqref{eq5.2} holds 
for $x_1\in[a^-,a^+]$.
\ep

\ms\noi
{\bf Proof}
It is not hard to show that 
the right-hand side of \eqref{eq1.12},
integrated with respect to $\rho^+$ is $0$. As a result, 
\[
A(f^1)_{x_2}-A(f^2)_{x_1}=0.
\]
 Let us set
\[
\xi=\ell_{x_2}-\ell*f^2+A(f^2)\ \ell.
\]
From differentiating both sides of \eqref{eq5.1} with respect
to $x_2$ we learn
\begin{align*}
\ell_{x_1x_2}=&\ell_{x_2}*f^1+\ell *f^1_{x_2}-A(f^1)_{x_2}\ \ell
-A(f^1)\ \ell_{x_2}\\
=&\xi*f^1+\ell*f^2*f^1-\big(A(f^2)\ell\big)*f^1
+\ell*f^1_{x_2}\\
&-A(f^1)_{x_2}\ \ell-A(f^1)\xi-A(f^1)(\ell*f^2)+A(f^1)A(f^2)\ell\\
=&\ell*\big[f^2*f^1-A(f^2)\otimes f^1-f^2
\otimes A(f^1)+f^1_{x_2}\big]\\
&-\ell\big[A(f^1)_{x_2}-A(f^1)A(f^2)\big] +\xi*f^1-A(f^1)\xi,
\end{align*}
where the operation $h\otimes k$ was defined in Remark 1.1.
Similarly,
\begin{align*}
\ell_{x_2x_1}=&\xi_{x_1}+\ell_{x_1}*f^2+\ell *f^2_{x_1}-A(f^2)_{x_1}\ \ell
-A(f^2)\ \ell_{x_1}\\
=&\xi_{x_1}+\ell*f^1*f^2-\big(A(f^1)\ell\big)*f^2
+\ell *f^2_{x_1}\\
&-A(f^2)_{x_1}\ \ell-A(f^2)(\ell*f^1)+A(f^2)A(f^1)\ell\\
=&\ell*\big[f^1*f^2-A(f^1)\otimes f^2-f^1
\otimes A(f^2)+f^2_{x_1}\big]\\
&-\ell\big[A(f^2)_{x_1}-A(f^1)A(f^2)\big] +\xi_{x_1}.
\end{align*}
From $\ell_{x_2x_1}=\ell_{x_1x_2}$, $A(f^2)_{x_1}=A(f^1)_{x_2}$,
and \eqref{eq1.12} we deduce
\be\la{uni}
\xi_{x_1}=\xi*f-A(f)\xi.
\ee
This means that $\xi(\cdot,x_2,\rho)$ satisfies the forward equation for
the Markov jump process $x_1\mapsto \rho(x_1,x_2)$ associated with the kernel $f$.
We wish to use the condition $\xi(a^+,x_2,\rho)=0$, to deduce that 
$\xi(x_1,x_2,\rho)=0$ for $x_1\in[a^-,a^2]$. This being true for
every $x_2\in[t_0,t_1]$ yields the desired result. Indeed
if $\var:\bR\to[0,\i)$ is a $C^1$ Lipschitz function such that 
$\var(0)=0$, and $\var(r)\ge |r|-c_0$, for some constant $c_0$, then
\begin{align*}
\frac{d}{dx_1}\int \var(\xi(x,\rho))\ \b(d\rho)
=&\int \var'(\xi(x,\rho))\ (\xi*f-A(f)\xi)(x,\rho)\ \b(d\rho)\\
\le &c_1\int |\xi(x,\rho)|\ \b(d\rho)\le 
c_1\int \var(\xi(x,\rho))\ \b(d\rho)+c_2,
\end{align*}
for constants $c_1$ and $c_2$. From this, Gronwall's inequality,
and the condition  $\xi(a^+,x_2,\rho)=0$, we deduce that $\var(\xi)=0$.
This completes the proof because we can approximate $|\xi|$ by functions of the form $\var(\xi)$, with $\var$ as above.
\qed

\bs\noi
{\bf(2)} As it is well-known, a time reversal of a 
Markov process can be realized as a Markov process with a generator that can be 
described in terms of original process and its marginals.
Indeed if we decrease $x_1\in[a^-,a^+]$, the process 
$x_1\mapsto \rho(x_1,x_2)$ is a Markov process with the jump rate
\[
\hat f(x,\rho^+,d\rho^-):=\eta(x,\rho^-,\rho^+)\ f(x,\rho^-,d\rho^+),\ \ \ {\text{ where }}\ \ \ \ \eta(x,\rho^-,\rho^+)=\frac{\ell(x,d\rho^-)}{\ell(x,d\rho^+)}.
\]
Similarly, as we decrease $x_2\in[t_0,t_1]$, the process 
$x_2\mapsto \rho(a^+,x_2)$ is a Markov process with the jump rate
$\hat f^2(a^+,x_2,\rho^+,d\rho^-)$, where 
\[
\hat f^2(x,\rho^+,d\rho^-):=\eta(x,\rho^-,\rho^+)\ f^2(x,\rho^-,d\rho^+).
\]
We also define
\be\la{eq5.3}
\tilde f(x,\rho^+,d\rho^-):=\hat f(-x,\rho^+,d\rho^-),
\ee
to represent the jump rate density of the process $x\mapsto\rho(-x)$.

 If $g$ is a convex function such that 
$\rho=(\rho_1,\rho_2)=\nabla g$
is distributed according to $\nu^{f,\L}$, and 
$\var(x_1,x_2)=(-x_2,-x_1)$, then $\hat g=g\circ \var$ is a convex
function that is defined on 
\[
\hat\L:=\var(\L)=[-t_1,-t_0]\x [-a^+,-a^-],
\]
 and $\hat\rho:=\nabla\hat g=(-\rho_2,-\rho_1)\circ \var$ is distributed
according to a probability measure that is denoted by $\hat\nu$.
Note $[\hat\rho^-,\hat\rho^+]=[\rho^-,\rho^+]^{-1}$.
According to $\hat\nu$, the process $x_1\mapsto\hat\rho(x_1,-a^+)$ 
is a jump process with the jump rate density
\[
 \tilde f^2(-a^+,x_1,-\rho^-,-\rho^+),
\]
with respect to the measure $\hat\b$ which is the push-forward 
of $\b$ under the map $\th(\rho):=-\rho$. Similarly, the process 
$x_2\mapsto\hat\rho(-t_0,x_2)$ 
is a jump process with the jump rate density
\[
 \tilde f^1(x_2,-t_0,-\rho^-,-\rho^+).
\]
We are now in a position to apply Theorem 1.1 to assert that the process
$x_1\mapsto \hat\rho(x_1,x_2)$ is a Markov jump process
for every $x_2\in[-a^+,-a^-]$. This in turn implies that
as we decrease $x_2$, the process
$x_2\mapsto \rho(x_1,x_2)$ is a (reversed) Markov jump process
for every $x_1\in[a^-,a^+]$, completing the proof of Theorem 1.2. 
To apply Theorem 1.1 though, we need to make sure that our candidate  the jump rate density of the jump
process $x_1\mapsto \hat\rho(x_1,x_2)$, namely 
\be\la{eq5.4}
  \bar f\left(x_1,x_2, \rho_1^-,\rho^-_2, \rho_1^+,\rho^+_2\right):= 
  \tilde  f^2\left(x_2,x_1, -\rho_2^-,-\rho^-_1, -\rho_2^+,-\rho^+_1\right),
  \ee
satisfies the kinetic equation. This will be carried out in Proposition 4.2.

\bp\la{pro4.2} Let $f=f^1$ be a solution of \eqref{eq1.12}.
 Then the following statements are true:
 
  \ms\noi
  {\bf(i)} The reversed kernel $\tilde f$, given by \eqref{eq5.3}
  satisfies \eqref{eq1.12}.

\ms\noi
 {\bf(ii)} The kernel $\bar f$ given by \eqref{eq5.4}
satisfies \eqref{eq1.12} where $\a$ is replaced with 
$\hat\a:=\a^{-1}$.
  \ep

\ms\noi
{\bf Proof (i)} Observe that  \eqref{eq5.1} and \eqref{eq5.2} can be rewritten as
\begin{align*}
\frac{\ell_{x_1}}{\ell}&=\frac{\ell*f^1}{\ell}-A(f^1)
=A(\hat f^1)-A(f^1),\\
\frac{\ell_{x_2}}{\ell}&=\frac{\ell*f^2}{\ell}-A(f^2)=A(\hat f^2)-A(f^2),
\end{align*}
As a consequence, 
\be\la{eq5.5}
\left(\tau\cdot\frac{\nabla \eta}{\eta}\right)(x,\rho^-,\rho^+)=
\frac {Q^-( f)}{ f}(x,\rho^-,\rho^+)+
\frac {Q^-(\hat f)}{\hat f}(x,\rho^+,\rho^-),
\ee
because
\[
\frac{Q^-( f)}{f}(x,\rho^-,\rho^+)=
A( f^2)(x,\rho^+)-A( f^2)(x,\rho^-)
-[\rho^-,\rho^+]\big(A( f^1)(x,\rho^+)-A( f^1)(x,\rho^-)\big).
\]
On the other hand, we can readily show
\[
\frac {\hat f^1*\hat f^2}{\hat f}(x,\rho^+,\rho^-)=
\frac {f^2*f^1}{ f}(x,\rho^-,\rho^+),\ \ \ \ 
\frac {\hat f^2*\hat f^1}{\hat f}(x,\rho^+,\rho^-)=
\frac {f^1*f^2}{ f}(x,\rho^-,\rho^+),
\]
which is an immediate consequence of $\eta(x,\rho^-,\rho^+)=\eta(x,\rho^-,\rho^*)\eta(x,\rho^*,\rho^+).$
From this, \eqref{eq5.5} and our assumption on $f$ we deduce
\begin{align*}
\left(\tau\cdot \frac{\nabla \hat f}{\hat f}\right)(x,\rho^+,\rho^-)
=&\left(\tau\cdot \frac{\nabla  f}{f}\right)(x,\rho^-,\rho^+)
+\left(\tau\cdot \frac{\nabla  \eta}{\eta}\right)(x,\rho^-,\rho^+)\\
=&\frac{Q^+( f)}{f}(x,\rho^-,\rho^+)-\frac{Q^-( f)}{f}(x,\rho^-,\rho^+)+\tau\cdot \frac{\nabla  \eta}{\eta}(x,\rho^-,\rho^+)\\
=&-\frac{Q^+(\hat f)}{\hat f}(x,\rho^+,\rho^-)
+\frac{Q^-( \hat f)}{\hat f}(x,\rho^+,\rho^-)
=-\frac{Q(\hat f)}{\hat f}(x,\rho^+,\rho^-),
\end{align*}
which is the reversed kinetic equation. This implies that
$\tilde f$ satisfies \eqref{eq1.12} because 
\[
\nabla \tilde f(x,\rho^-,\rho^+)=- \nabla\hat f(-x,\rho^-,\rho^+).
\]

\ms\noi
{\bf(ii)} Observe that if $\bar f^2=\hat\a \bar f$, then
\[
  \bar f^2\big(x_1,x_2, \rho_1^-,\rho^-_2, \rho_1^+,\rho^+_2\big)= 
  =\tilde  f\big(x_2,x_1, -\rho_2^-,-\rho^-_1, -\rho_2^+,-\rho^+_1\big).
  \]  
By {\bf(i)}, we know that  $\a \tilde f_{x_1}-\tilde f_{x_2}=-Q(\tilde f).$ 
After swapping $x_1$ with $x_2$ we deduce
\[
\bar f_{x_2}-\hat \a \bar f_{x_1}=-Q(\tilde f).
\]  
Finally observe that $-Q(\tilde f)=Q(\bar f)$ because when $\tilde f^1$
is swapped with $\tilde f^2$, the sign of $Q$ changes.
\qed

\bs\noi
{\bf Remark 4.2} An alternative strategy for completing the proof of Theorem 1.2 is to use
Proposition 1.1. We already know that $\ell$ satisfies \eqref{eq5.2}. On the other hand, since we also know that the process $x_2\mapsto
\rho(x_1,x_2)$ is a Markov jump process, the measure $\ell$ also satisfies
\[
\ell_{x_2}=\ell*h-A(h)\ell,
\]
where $h$ is its jump rate.
One should be able to deduce from this that $h=f^2$. 
\qed

\appendix
\section{The Kinetic Equation}
The purpose of this section is to prove the existence of a solution of the kinetic equation. To have a more conventional notation,
we write $(x,t)$ for $(x_1,x_2)$ throughout this section.
We start first with the following notation:

\ms\noi
{\bf Notation A.1 (i)}
 We fix $P^-<P^+$ two real numbers, such that the range of our piecewise constant function $\rho$ is in the box $[P^-,P^+]^2$.

\ms\noi
{\bf (ii)} For any measure space $\cE$, let $\cF_b(\cE)$ be the space of real-valued bounded measurable functions defined on $\cE$. 

\ms\noi
{\bf (iii)}
We introduce the function space  $\cX$ to be the set kernels
$h \in \cF_b(\bR \times ([P^-,P^+]^2)^2)$ such that 
 $x \mapsto h(x,\rho^-,\rho^+)$ is $C^1$ and Lipschitz 
 for  all 
$\rho^-$ and $\rho^+$.

\ms\noi
{\bf(iv)}
We equip $\cX$ with the following norm
\[
\|h\|_{\cX}:= \sup_{x \in \bR} \ \sup_{\rho^-,\rho^+}\
\left[ |h(x,\rho^-,\rho^+)| + |\partial_x h(x,\rho^-,\rho^+)|
\right].
\]
It is standard that $(\cX,\|\cdot\|_{\cX})$ is a Banach space.

\ms\noi
{\bf (v)} For any $v \ge 0$, let $\G^{v}$ and $\G^v_+$ be the sets
\begin{align*}
\G^{v}:=&\Big \{(\rho^-,\rho^+) \in ([P^-,P^+]^2)^2 :
\  \rho^-\prec\rho^+, \ |[\rho^-,\rho^+]| \le v \Big \}\\
=&\Big \{(\rho^-,\rho^+) \in ([P^-,P^+]^2)^2 :
\  \rho^+-\rho^-\in C^v \setminus \{0\}\Big \},\\
\G^{v}_+:
=&\Big \{(\rho^-,\rho^+) \in ([P^-,P^+]^2)^2 :
\  \rho^+-\rho^-\in C^v_+ \setminus \{0\}\Big \},
\end{align*}
where $C^v$ and $C^v_+$ are the cones
\begin{align}\la{cone}
C^v&=\Big\{m=(m_1,m_2)\in\bR^2:\ m_1\ge 0, \ |m_2|
\le vm_1\Big\},\\
C_+^v&=\Big\{m=(m_1,m_2)\in\bR^2:\ m_1,m_2\ge 0, \ m_2
\le vm_1\Big\}.\nonumber
\end{align}

\ms\noi
{\bf (vi)} Let $V_{\infty} \ge 0$, and $\d_0>0$. We write 
$\cX(V_\i,\d_0)$ for the set of
$h \in \cX$ with the following properties:
\ms

{\bf (1)} The function $h(\cdot,\rho^-,\rho^+)$ is zero for $(\rho^-,\rho^+) \notin \G^{V_{\i}}$. \ms

{\bf (2)} There exists a constant $\d_0>0$ such that $\inf_{x \in \bR} \inf_{(\rho^-,\rho^+) \in \G^{V_{\i}}} h(x,\rho^-,\rho^+) \ge \d_0$. \ms

{\bf (3)} For all $\rho^-,\rho^+$, the function $x \mapsto h(x,\rho^-,\rho^+)$ is $C^2$ such that 
\[ 
\sup_{x \in \bR} \sup_{\rho^-,\rho^+} |\partial^2_x h(x,\rho^-,\rho^+)| < \i.
\]
Likewise, we write 
$\cX_+(V_\i,\d_0)$ for the set of
$h \in \cX_+(V_\i,\d_0)$ with the similar properties, except that
the set $\G^{V_{\i}}$  in {\bf(i)} and 
{\bf(iii)} is replaced with $\G^{V_{\i}}_+$.
\qed

\bs
The following theorem proves the existence of a local solution of the kinetic equation.
\ms

\begin{theorem} Given $h\in\cX(V_\i,\d_0)$,
denote by $M_0:=\sup_{x \in \bR} \sup_{\rho^-,\rho^+} h(x,\rho^-,\rho^+)$, and define the time
\[
T^{*}:=\min \left(\frac{1}{12V_{\i}M_0}\ ,\ 
\frac{\d_0}{48V_{\i}M_0^2}\right).
\]
Then, there exists a unique solution 
\[
f: \bR \times [0,T^{*}] \times ([P^-,P^+]^2)^2 \to \bR,
\]
of the kinetic equation 
\[
f_t-[\rho^-,\rho^+]f_x=Q(f)=:Q^{+}(f)-Q^{-}(f),
\]
where
\begin{align*}
Q^{+}(f)(x,t,\rho^-,\rho^+)=&\int \left([\rho^+,\rho^*]-[\rho^*,\rho^-]\right)f(x,t,\rho^-,\rho^*)f(x,t,\rho^*,\rho^+)\ \b(d\rho^*),\\
Q^{-}(f)(x,t,\rho^-,\rho^+)= &\left(\int ([\rho^+,\rho^*]-[\rho^-,\rho^+])f(x,t,\rho^+,\rho^*)\
\b(d\rho^*) \right. \\& \left. \ \ \ \ -
\int ([\rho^-,\rho^*]-[\rho^-,\rho^+])
f(x,t,\rho^-,\rho^*)\ \b(d\rho^*) \right) f(x,t,\rho^-,\rho^+),
\end{align*}
with $f(\cdot,0,\cdot,\cdot)=h$.
The function $f$ is $C^1$ in the variables $(x,t)$ for all fixed $\rho^-,\rho^+$ and $f(\cdot,\cdot,\rho^-,\rho^+) \equiv 0$ for all $(\rho^-,\rho^+) \notin \G^{V_{\i}}$. Furthermore, we have that
\[
\sup_{t \in [0,T^{*}]} \|f(\cdot,t,\cdot,\cdot)\|_{\cX} < \i,\ \ \ \
\inf_{t \in [0,T^{*}]} \inf_{x \in \bR} \inf_{(\rho^-,\rho^+) \in \G^{V_{\i}}} f(x,t,\rho^-,\rho^+) \ge \frac{\d_0}{2}.
\]
Moreover if $h\in\cX_+(V_\i,\d_0)$, then there exists a unique solution $f$ with similar properties except that the set $\G^{V_{\i}}$ must be replaced with $\G^{V_{\i}}_+$.
\end{theorem}

\ms\noi
{\bf {Proof}}  {\em(Step 1)} We assume here without loss of generality that the measure $\beta$ has total mass $1$ on the box $[P^-,P^+]^2$. By the following standard change of variables, we transform the previous PDE to an ODE. For instance, define the function $g$ as
\begin{equation}\label{change-var}
g(x,t,\rho^-,\rho^+)=f(x-[\rho^-,\rho^+]t,t,\rho^-,\rho^+) .
\end{equation}
Then by the chain rule, $g$ must verify the following ODE
\[
g_t=\tilde{Q}^{+}(g)-\tilde{Q}^{-}(g)
=\tilde{Q}^{+}(g)-\tilde L(g)g,
\]
where
\begin{align*}
\tilde{Q}^{+}(g)(x,t,\rho^-,\rho^+)=&
\int [\rho^-,\rho^*,\rho^+]\ g\left(x-([\rho^-,\rho^+]-[\rho^-,\rho^*])t,t,\rho^-,\rho^*\right)  \\
& \ \ \ \ \ \ \ \ \ \ \ \ \ \ \ \ \ \ \ \ \ \ \
g(x-([\rho^-,\rho^+]-[\rho^+,\rho^*])t,t,\rho^*,\rho^+)\ \ \beta(d\rho^*),\\
 \tilde{L}(g)(x,t,\rho^-,\rho^+)= &\int 
([\rho^+,\rho^*]-[\rho^-,\rho^+])\ g(x-([\rho^-,\rho^+]-[\rho^+,\rho^*])t,t,\rho^+,\rho^*)\ \beta(d\rho^*)\\ 
&\ \ \ \  -
\int ([\rho^-,\rho^*]-[\rho^-,\rho^+])\ g(x-([\rho^-,\rho^+]-[\rho^-,\rho^*])t,t,\rho^-,\rho^*)\ \beta(d\rho^*) .
\end{align*}
We will prove the existence of  a solution $g$ by an approximation scheme and then recover the desired $f$ via the equation \eqref{change-var}. Define the functional $\cH: \cX \times \bR \to \cX$,
$\cH(g,t):=\cH^+(g,t)-\cK(g,t)g$,  by
\begin{align*}
\cH^+(h,t)(x,\rho^-,\rho^+) :=&\int 
[\rho^+,\rho^*,\rho^-]\ h(x-([\rho^-,\rho^+]-[\rho^-,\rho^*])t,\rho^-,\rho^*) \\
& \ \ \ \ \ \ \ \ \ \ \ \ \ \ \ \ \ \ \ \ \
h(x-([\rho^-,\rho^+]-[\rho^+,\rho^*])t,\rho^*,\rho^+)\ \beta(d\rho^*),\\
\cK^+(h,t)(x,\rho^-,\rho^+) :=& \int 
([\rho^+,\rho^*]-[\rho_-,\rho^+])\ h(x-([\rho^-,\rho^+]-[\rho^+,\rho^*])t,\rho^+,\rho^*)\ \beta(d\rho^*)  \\ 
& \ \ \ \ \ \ \ \ \ -
\int ([\rho^-,\rho^*]-[\rho^-,\rho^+])\
h(x-([\rho^-,\rho^+]-[\rho^-,\rho^*])t,\rho^-,\rho^*)\ \beta(d\rho^*) .
\end{align*}
Our goal is to prove the existence of a \textit{local} solution $g: [0,T^{*}] \mapsto \cX$ to the inhomogeneous ODE 
\begin{equation}\label{ode}
\dot{g}(t)=\cH(g(t),t)
\end{equation}
under the initial condition $g(0)=h$. As the function space $(\cX,\|\cdot\|_{\cX})$ is clearly Banach, we will construct a Cauchy sequence $(g_n)_{n \in \bN}$ of elements in $C([0,T^{*}],\cX)$ that will converge to our desired solution $g$.

\ms\noi
{\em(Step 2)}
For any fixed $n \in \bN$, we define the polygonal function $g_n$ such that $g_n(0)=h$ and 
\[
\dot{g}_n(t)=\cH \left(g_n\left(\frac{j}{n}\right),\frac{j}{n} \right) \text{ for all } t \in \left(\frac{j}{n},\frac{j+1}{n}\right)
\]
for all $j \ge 0$. Let us denote $g_n^{j}=g_n \left(\frac{j}{n} \right)$, then it is clear that all $g_n^{j}$ are $C^2$ in the variable $x$. We have that
\begin{equation}\label{induct}
n(g_n^{j+1}-g_n^{j})=\cH \left(g_n^{j},\frac{j}{n} \right).
\end{equation}
Let us prove first that $g_n^{j}(\cdot,\rho^-,\rho^+) \equiv 0$ for all $(\rho^-,\rho^+) \notin \G^{V_{\i}}$ by induction on $j$. Suppose this is true for $j$ and we wish to prove it for $j+1$.
Take  $x \in \bR$, 
$(\rho^-,\rho^+) \notin \G^{V_{\i}}$, 
and take any $\rho^*$ such that $\rho^- \prec \rho^* \prec \rho^+$ . Since $C^v$ of \eqref{cone} is a cone, we have that
either $(\rho^-,\rho^*) \notin \G^{V_{\i}}$ or
$(\rho^*,\rho^+) \notin \G^{V_{\i}}$. In either cases
\[
g_n^{j}\left(x-\left([\rho^-,\rho^+]-[\rho^-,\rho^*]\right)\frac{j}{n},\rho^-,\rho^*\right) 
g_n^{j}\left(x-([\rho^-,\rho^+]-[\rho^+,\rho^*])\frac{j}{n},\rho^*,\rho^+\right)= 0,
\]
by the induction hypothesis. As a result, 
$g_n^{j+1}(x,\rho^-,\rho^+) = 0$, as desired.

 Next, let us define 
\begin{align*}
m_j:&=\inf_{x \in \bR} \inf_{(\rho^-,\rho^+) \in \G^{V_{\i}}} g_n^{j}(x,\rho^-,\rho^+),\ \ \ \ \ \ \ \ \
M_j:=\sup_{x \in \bR} \sup_{(\rho^-,\rho^+) \in \G^{V_{\i}}} |g_n^{j}(x,\rho^-,\rho^+)|,\\
 M'_j:&=\sup_{x \in \bR} \sup_{(\rho^-,\rho^+) \in \G^{V_{\i}}} |\partial_x g_n^{j}(x,\rho^-,\rho^+)| ,
\ \ \ \ 
M''_j:=\sup_{x \in \bR} \sup_{(\rho^-,\rho^+) \in \G^{V_{\i}}} |\partial^2_x g_n^{j}(x,\rho^-,\rho^+)|,
\end{align*}
It is clear from the expression of $\cH$ that we have for all $j \ge 0$,
\begin{equation}\label{inequ-m}
n(m_{j+1}-m_j) \ge -6V_{\i}M_j^2,\ \ \ \ 
n(M_{j+1}-M_j) \le 6V_{\i} M_j^2.
\end{equation}
Let us prove first by induction on $j$ the following inequality,
\begin{equation}\label{inequ-M}
rj<1\ \ \ \implies\ \ \ M_j \le M_0 (1-rj)^{-1}.
\ee
where $r=\frac{6V_{\i}M_0}{n}$.
The verification for $j=0$ is trivial. Assume it is true for $j$, then from the second inequality in \eqref{inequ-m}, it suffices to prove that
\[
(1-rj)^{-1} \left(1+r(1-rj)^{-1}\right) \le \left(1-r(j+1) \right)^{-1}.
\]
This inequality is equivalent to
\[
(1-(j-1)r)(1-(j+1)r) \le (1-jr)^2,
\]
which is clearly true. 
As an immediate consequence of \eqref{inequ-M} we have that
\[
\sup_{t \in [0,T^{*}]} \|g_n(t)\|_{L^{\i}} \le M_0 \sup_{t \in [0,T^{*}]} \left(1-\frac{6 \lf nt \rf V_{\i} M_0}{n}\right)^{-1} \le 2M_0.
\]

By differentiating the identity \eqref{induct}, we also have that 
\[
n(M'_{j+1}-M'_j) \le 12V_{\i} M'_j M_j,
\]
So for all $j$ such that $\frac{j+1}{n} \le T^*$, we have that
\[
M'_{j+1} \le M'_j\left(1+\frac{24V_{\i}M_0}{n}\right),
\]
from which it follows that
\[
M'_j \le M'_0\left(1+\frac{24V_{\i}M_0}{n}\right)^j\le M'_0 e^{\frac{2j}{nT^*}}\le M'_0 e^2,
\]
and hence
\[
\sup_{t \in [0,T^{*}]} \|\partial_x g_n(t)\|_{L^{\i}} \le M'_0e^2 .
\]
Likewise, by differentiating twice the identity \eqref{induct}, we get that
\[
n(M''_{j+1}-M''_j) \le12V_{\i}M_jM''_j+12V_{\i}(M'_j)^2\le 24V_{\i}M_0 M''_j+12V_{\i}(M'_0)^2e^4.
\]
From this, it follows by similar arguments as before that
\[
M''_j \le \left(M''_0+\frac{(M'_0)^2e^4}{2M_0}\right)e^{\frac{24V_{\i}M_0j}{n}},
\]
and hence
\[
\sup_{t \in [0,T^{*}]} \|\partial^2_x g_n(t)\|_{L^{\i}} \le \left(M''_0+\frac{(M'_0)^2e^4}{2M_0}\right)e^{2}.
\]
Now since we have for every $j$ such that $\frac{j+1}{n} \le T_{*}$, 
\[
m_{j+1} \ge m_j -\frac{24V_{\i}M_0^2}{n},
\]
it follows easily that
\[
\inf_{t \in [0,T_{*}]} \inf_{x \in \bR} \inf_{(\rho^-,\rho^+) \in \G^{V_{\i}}} g_n(t)(x,\rho^-,\rho^+) \ge \frac{\d_0}{2} >0.
\]
We have hence proved that all the approximating functions
 $(g_n)_{n \in \bN} \in \cC([0,T^{*}],\cX)$ are supported on 
$\G^{V_{\i}}$ in the $(\rho^-,\rho^+)$ variables,
and are uniformly bounded from above and below by positive constants 
in their supports.

\ms\noi
 {\em (Step 3)} To finish the proof, we shall show that 
the sequence $\{g_n\}$  is Cauchy. This is achieved by obtaining Lipschitz estimates on $g_n$. Observe that
for any $s<t$, and $k_1,k_2 \in \cX$ that are $C^2$ in the $x$-variable and supported on $\G^{V_{\i}}$ such that
\[
\max(\|\partial^2_x k_1\|_{L^{\i}},\|\partial^2_x k_2\|_{L^{\i}})<\i,
\]
 it is straightforward to show 
\begin{align*}
\| \cH(k_1,t)-\cH(k_2,t)\|_{\cX} \le &6V_{\i}\left(\|k_1\|_{\cX}
+\|k_2\|_{\cX}\right)\|k_1-k_2\|_{\cX},\\
\|\cH(k_1,t)-\cH(k_1,s)\|_{\cX} \le &72V_{\i}^2\|k_1\|_{\cX}
\left(\|k_1\|_{\cX}+\|\partial^2_x k_1\|_{L^{\i}}\right)(t-s).
\end{align*}
Let us denote 
\[
M:=\max \left(2M_0,M'_0e^2, \left(M''_0+\frac{(M'_0)^2e^4}{2M_0}\right)e^{2}\right)
\]
The constant $M$ is a uniform upper bound on the supremum norm of $g_n(t),\partial_x g_n(t),\partial^2_x g_n(t)$ for all $t \in [0,T^{*}]$ and $n \in \bN$. We have that
\begin{align*}
\|\dot{g}_n(t)-\dot{g}_m(t)\|_{\cX} =&\left \|\cH\left(g_n \left(\frac{\lf nt \rf}{n}\right),\frac{\lf nt \rf}{n}\right)-\cH\left(g_m \left(\frac{\lf mt \rf}{m}\right),\frac{\lf mt \rf}{m}\right) \right\|_{\cX}\\
  \le& \left\|\cH\left(g_n \left(\frac{\lf nt \rf}{n}\right),\frac{\lf nt \rf}{n}\right)-\cH\left(g_n \left(\frac{\lf nt \rf}{n}\right),t\right) \right\|_{\cX}\\
&+ \left\|\cH\left(g_n \left(\frac{\lf nt \rf}{n}\right),t\right)-\cH\left(g_n(t),t\right)\right\|_{\cX}\\
&+ \left\|\cH(g_n(t),t)-\cH(g_m(t),t) \right\|_{\cX}\\
& + \left\|\cH(g_m(t),t)-\cH\left(g_m \left(\frac{\lf mt \rf}{m}\right),t\right) \right\|_{\cX}\\
& + \left\|\cH\left(g_m \left(\frac{\lf mt \rf}{m}\right),t\right)-\cH\left(g_m \left(\frac{\lf mt \rf}{m}\right),\frac{\lf mt \rf}{m}\right) |\right\|_{\cX} \\
 \le &144V_{\i}^2M^2 \left(\frac{1}{n}+\frac{1}{m}\right)+12V_{\i}M \left( \left\|g_n\left(\frac{\lf nt \rf}{n}\right)-g_n(t)\right\|_{\cX}\right. \\ &+  \left. \left\|g_m\left(\frac{\lf mt \rf}{m}\right)-g_m(t)\right\|_{\cX}\right)+12V_{\i}M\|g_n(t)-g_m(t)\|_{\cX}.
\end{align*}
On the other hand,
\[
\left\|g_n \left(\frac{\lf nt \rf}{n} \right)-g_n(t)\right\|_{\cX} \le \frac{1}{n}\left\| \cH \left(g_n \left(\frac{ \lf nt \rf}{n} \right)\right) \right\|_{\cX} \le \frac{6V_{\i} M^2}{n}
\]
and similarly for the term concerning $m$. Hence there exist
 two positive constants $C_1,C_2$ that only depend on 
$V_{\i}$ and $M$, 
 such that for all $t \in [0,T^{*}]$, 
\[
\|\dot{g}_n(t)-\dot{g}_m(t)\|_{\cX} \le C_1\left(\frac{1}{n}+\frac{1}{m}\right)+C_2\|g_n(t)-g_m(t)\|_{\cX},
\]
which implies 
\[
\left|\left|g_n(t)-g_m(t)\right|\right|_{\cX} \le C_1 \left(\frac{1}{n}+\frac{1}{m}\right)t+C_2\int_{0}^{t}\|g_n(s)-g_m(s)\|_{\cX}\ ds.
\]
This, and the Gronwall's inequality give
\[
\sup_{t \in [0,T^{*}]} \|g_n(t)-g_m(t)\|_{\cX} \le C_1\left(\frac{1}{n}+\frac{1}{m} \right)T^{*}\left(1+c_2T^{*}e^{C_2T^{*}}\right)
\]
which implies that $(g_n)_{n \in \bN}$ is a Cauchy sequence and therefore admits a limit $g_{\i} \in \cC([0,T^{*}],\cX)$. The function $g_{\i}$ (that we now regard as a function of the four variables $(x,t,\rho^-,\rho^+)$) is $C^1$ in the variables $x$ and $t$, and verify the inhomogeneous ODE \eqref{ode} and is bounded uniformly from below by $\frac{\d_0}{2}$ and is such that
\[
\sup_{t \in [0,T^*]} \sup_{x \in \bR} \sup_{\rho^-,\rho^+} g_{\i}(x,t,\rho^-,\rho^+) \le M.
\]
Moreover, for any fixed $x$ and $t$ in its domain of definition, the function $(\rho^-,\rho^+) \mapsto g_{\i}(x,t,\rho^-,\rho^+)$ is supported on $\G^{V_{\i}}$. Now defining
\begin{equation}\label{ode-l}
f(x,t,\rho^-,\rho^+)=g_{\i}(x+[\rho^-,\rho^+]t,t,\rho^-,\rho^+)
\end{equation}
$f$ is again $C^1$ in $x$ and $t$, verify the same properties as $g_{\i}$ and verifies the desired kinetic equation. 

Finally we remark that in the above proof, we may replace the set $\G^{V_{\i}}$ with $\G^{V_{\i}}_+$.
\qed

\bs
For the second part of this section, we will prove the existence of the solution to the Kolmogorov forward equation both in space $x$ and time $t$. More precisely, we wish to address the existence of a unique 
uniformly positive solution $\ell$ of the equations \eqref{kin-x} and
\eqref{kin-t}, provided that the kernel $f$ is uniformly positive. 
We remark that these equations are consistant by
Proposition 4.1. Because of this, we only need to solve \eqref{kin-x}
in $[a^-,a^+]$ for an initial condition $\ell (a^-,t,\cdot)$ that solves
\eqref{kin-t}. The existence of a solution to \eqref{kin-x} can be carried
out by standard arguments. However, we need to ensure the constructed
solution is uniformly positive in $\L$, if the initial $\ell^0(\rho)=\ell
(a^-,t_0,\rho)$ is uniformly positive. Observe that if $\ell$ solves \eqref{kin-x}, then
\[
\frac{d}{dx}\int\ell (x,t,\rho)\ \b(d\rho)=0,
\]
because the $\b$-integral of the right-hand side of  \eqref{kin-x} is
$0$. This means 
\be\la{int}
\int\ell(x,t)\ \b(d\rho)=1,
\ee
if this is the case for $x=a^-$. On the other hand, if the total integral
of $\b$ is one, $f\ge \d_1$
for some positive constant $\d_1$,
and $\ell$ is a solution of \eqref{kin-x} satisfying \eqref{int}, then
\[
\ell_x(x,t,\rho)\ge \d_1-\l(x,t,\rho)\ell(x,t,\rho),
\]
which leads to the lower bound
\[
\ell(x,t,\rho)\ge \ell(a^-,t,\rho)e^{-\int_{a^-}^x \l(\th,t,\rho)\ d\th}
+\d_1\int_{a^-}^xe^{-\int_{y}^x \l(\th,t,\rho)\ d\th}\ dy.
\]
From this we learn that $\ell$ is uniformly positive in $\L$
if this is the case on the left boundary side of $\L$. By assumption,
$\ell$ is uniformly positive at $(a^-,t_0)$, and as $t$ varies,
the function $t\mapsto\ell(a^-,t,\rho)$ satisfies \eqref{kin-t}.
If the kernel $f$ is supported in $\G^{V_{\i}}_+$, then 
$[\rho^-,\rho^+]f\ge 0$, and a repetition of the above reasoning 
guarantees 
\[
\ell(a^-,t,\rho)\ge \ell(a^-,t_0,\rho)e^{-\int_{t_0}^x A(a^-,\th,\rho)\ d\th}.
\]
In summary, when the kernel $f$ is supported in $\G^{V_{\i}}_+$, and is
uniformly positive on its support, we can construct a unique uniformly positive solution $\ell$ to forward equations \eqref{kin-x} and
\eqref{kin-t} by standard arguments. However some care is needed
if $[\rho^-,\rho^+]$ can change sign in the support of our kernel $f$.
In this case, we can guarantee the existence of a uniformly positive
solution to \eqref{kin-x} and
\eqref{kin-t} if we either replace the time interval $[0,T^*]$ with a shorter interval, or assume that the initial $\ell(a^-,t_0,\rho)$
is sufficiently positive.  As an example, we demonstrate how a lower 
bound of $1/6$ on the initial $\ell$ can guarantee the positivity of the solution.

\ms

\begin{theorem} Fix $a^{-}<a^{+}$. Let $\ell^{0} : [P^-,P^+]^2 \to [0,+\infty)$ be a measurable function such that there exists two constants $c,C>0$ with
\[
c \le \ell^{0}(\rho) \le C \text{ for all } \rho
\]
and $\int \ell^{0}(\rho)\ \b(d\rho)=1$. Moreover, assume that $c \ge \frac{1}{6}$. 
Then there exists a $C^1$ solution $\ell : [a^{-},a^{+}] \times [0,T^{*}] \times [P^-,P^+]^2 \to [0,+\infty)$ to the equations \eqref{kin-x} and \eqref{kin-t} such that
$\ell(a^{-},0,\cdot)=\ell^{0}$, and such that $\ell$ is uniformly bounded below by a positive constant and 
\[
\int \ell(x,t,\rho)\ \b(d\rho)=1,
\]
for all $(x,t) \in [a^{-},a^{+}] \times [0,T^{*}]$.
\end{theorem}

\ms\noi
{\bf Proof } Without loss of generality let us assume that $a^{-}=0$ and denote $a^{+}=a$. We will construct a two-parameter function $\ell :[-V_{\i}T^{*},a] \times [0,T^{*}] \to \cF_b([P^-,P^+]^2)$. The reason why we extend the space domain to $[-V_{\i}T^{*},a]$ instead of $[0,a]$ will be made clear later. 
Let us define first $\ell(\cdot,0)$ on $[-V_{\infty}T^{*},a]$ using the first ODE in the $x$-direction. The utility of the condition $c \ge \frac{1}{6}$ is to ensure the non-negativity of $\ell$ as we run the ODE backwards from $0 \rightarrow -V_{\i}T^{*}$. Our strategy for proving the existence of the solution of the ODE at $t=0$ is done in a similar fashion as the kinetic equation via an approximation scheme. In other words, we construct a polygonal approximating $\ell_n : [-V_{\i}T^{*},a] \to \cF_b([P^-,P^+]^2)$ by putting $\ell_n(0)=\ell^{0}$, and for any $k \in \bZ$ by the inductive relation.  More precisely,
we put $f^k(\rho,\rho_*):=f(k/n,0,\rho,\rho_*)$, and require that functions $\ell_n^{k}:=\ell_n \left(\frac{k}{n} \right) \in \cF_b([P^-,P^+]^2)$ to satisfy
 \begin{align*}
n \left(\ell^{k+1}_n(\rho)-\ell^k_n(\rho)\right)=&\int f^k \left(\rho^*,\rho\right) \ell_n^k(\rho^*)\ \b(d\rho^*)  -\left( \int f^k\left(\rho,\rho^+\right)\ \b(d\rho^+) \right) \ell_n^k(\rho ), 
\end{align*}
for $k\ge 0$, and 
\begin{align*}
-n \left(\ell^{k-1}_n(\rho)-\ell^k_n(\rho)\right)=&\int f^k \left(\rho^*,\rho\right) \ell_n^k(\rho^*)\ \b(d\rho^*)  -\left( \int f^k\left(\rho,\rho^+\right)\ \b(d\rho^+) \right) \ell_n^k(\rho ), 
\end{align*}
for $k \le 0$.
The intermediate values $\ell_n(x)$ for $x \in (\frac{k}{n},\frac{k+1}{n})$ are obtained by linear interpolation. As an initial observation, remark that
\[
\int \ell_n^{k}(\rho)\ \b(d\rho)=\int \ell_n^{k \pm 1}(\rho)\ \b(d\rho),
\]
and hence
\[
\int \ell_n^{k}(\rho)\ \b(d\rho)=1,\ \ \  \text{ for all } \ k \in \bZ.
\]
Now, if we take $n \ge M_0$ where 
$M_0=\|f(\cdot,0,\cdot,\cdot)\|_{L^{\i}}$, then by induction it follows that $\ell_n^{k} \ge 0$ for all $k \ge 0$, as we have that
\begin{align*}
\ell_n^{k+1}(\rho)&=\ell_n^{k}(\rho)+\frac{1}{n} \left(\int f^k \left(\rho^*,\rho\right) \ell_n^k(\rho^*)\ \b(d\rho^*) -  
 \left( \int f^k\left(\rho,\rho^+\right)\ \b(d\rho^+) \right) \ell_n^k(\rho) \right)  \\
& \ge \ell_n^{k}(\rho) - \frac{M_0}{n} \ell_n^{k}(\rho), 
\end{align*}
which in turn implies the following lower bound 
\[
\ell_n^{k}(\rho) \ge \ell^0(\rho)\left(1 -\frac{M_0}{n} \right)^k \ge \ell^{0}(\rho) e^{-\frac{M_0 k}{n}}\ \  \text{ for all } \ k \ge 0.
\]
On the other hand, for $k \le 0$ we have
\begin{align*}
\ell_n^{k-1}(\rho)&=\ell_n^{k}(\rho)-\frac{1}{n} \int f^k (\rho^*,\rho) \ell_n^k(\rho^*)\ \b(d\rho^*) +\frac 1n   
\left( \int f^k\left(\rho,\rho^+\right)\ \b(d\rho^+) \right) \ell_n^k(\rho) ,
\end{align*}
which leads to 
\[
 \ell_n^{k}(\rho) \ge \ell^{0}(\rho)-\frac{M_0k}{n},
\]
because
\[
\int f^k \left(\rho^*,\rho\right) \ell_n^k(\rho^*)\ \b(d\rho^*) \le M_0 \int \ell_n^{k}(\rho^*)\ \b(d\rho_*) =M_0.
\]
In particular, if $\frac{k}{n} \ge -V_{\i}T^{*}$, then 
$\frac{M_0k}n\ge -\frac 12$, and 
\[
\inf_{\rho} \ell_n^{k}(\rho) \ge c -\frac{1}{12} \ge \frac{1}{12}.
\]
We have therefore constructed the polygonal approximating function $\ell_n : [-V_{\i}T^{*},a] \to \cF_b([P^-,P^+]^2)$ such that it is uniformly bounded from below by $\min({1}/{12},c e^{-M_0 a })$.
The sequence $(\ell_n)_{n \in \bN}$ is a Cauchy sequence in the space $C([-V_{\i}T^{*},a],\cF_b([P^-,P^+]^2)$ where $\cF_b([P^-,P^+]^2)$ is viewed as a Banach space equipped with the uniform norm. We obtain that the limit $\ell_{\i}:=\lim_{n \to \i} \ell_n$ is a solution to the ODE
\[
(\ell_{\i})_x(x,\rho)=\int f(x,0,\rho^*,\rho)\ell_{\i}(x,\rho^*)\ \b(d\rho^*) - \left( \int f(x,0,\rho,\rho^+)\ \b(d\rho^+) \right) \ell_{\i}(x,\rho)
\]
We define $\ell(\cdot,0)=\ell_{\i}$. We will move on now to prove the existence of the solution $\ell$ as an ODE in the time variable $t$. In order to preserve the non-negativity of $\ell$, we have taken advantage in the ODE in the $x$-direction of the positivity of the kernel $f(x,t,\rho^-,\rho^+)$, however in the $t$-direction the kernel is equal to $[\rho^-,\rho^+]f(x,t,\rho^-,\rho^+)$. To circumvent this difficulty, we take advantage of the finite speed propagation (this also explains why we have constructed $\ell(0,\cdot)$ on $[-V_{\i}T^{*},a]$ instead of just $[0,a]$). For any $x \in \bR$ and $t \in [0,T^{*}]$ we define
\[
\tilde{f}(x,t,\rho^-,\rho^+)=f(x+V_{\i}t,t,\rho^-,\rho^+) \text{ for all } \rho^-,\rho^+.
\]
We define a function $\tilde{\ell} :[-V_{\i}T^{*},a] \times [0,T^{*}] \to \cF_b([P^-,P^+]^2)$ that satisfies the initial condition
 $\tilde{\ell}(x,0)=\ell_{\i}(x)=\ell(x,0)$. Now, for $x=-V_{\i}T^{*}$ we define $\tilde{\ell}(-V_{\i}T^{*},\cdot) : [0,T^{*}] \to \cF_b([P^-,P^+]^2)$ by solving the ODE
\begin{align*}
\tilde{\ell}_t(-V_{\i}T^{*},t,\rho)=&\int \left([\rho^*,\rho]+V_{\i} \right) \tilde{f}(-V_{\i}T^{*},t,\rho^*,\rho)\tilde{\ell}(-V_{\i}T^{*},t,\rho^*)\ \b(d\rho^*)  \\
&-\left( \int \left([\rho,\rho^+]+V_{\i}\right)\tilde{f}(-V_{\i}T^{*},t,\rho,\rho^+)\ \b(d\rho^+) \right)\tilde{\ell}(-V_{\i}T^{*},t,\rho),
\end{align*}
with initial condition 
$\tilde{\ell}(-V_{\i}T^{*},0)=\ell_{\i}(-V_{\i}T^{*})$.
Now, for any fixed $t \in (0,T^{*}]$ we define 
$\tilde{\ell} (\cdot,t) : [-V_{\i}T^{*},a] \to \cF_b([P^-,P^+]^2)$ by solving the ODE on $[-V_{\i}T^{*},a]$ 
\[
\tilde{\ell}_x(x,t,\rho)=\int \tilde{f}(x,t,\rho^*,\rho)\tilde{\ell}(x,t,\rho^*)\ \b(d\rho^*) - \left( \int \tilde{f}(x,t,\rho,\rho^+)\ \b(d\rho^+) \right) \tilde{\ell}(x,t,\rho)
\]
with initial condition determined by $\tilde{\ell}(-V_{\i}T^{*},t)$. The existence of these solutions is done by exactly the same approximation scheme than before, and the function $\tilde{\ell}$ is bounded uniformly from below on the box $[-V_{\i}T^{*},a]\times [0,T^{*}]$ due to the non-negativity of the kernels
$([\rho^-,\rho^+]+V_{\i})f(x,t,\rho^-,\rho^+)$ and $f(x,t,\rho^-,\rho^+)$. Moreover, if we assume that initially $f(0,\cdot)$ is $C^3$ then we get that $f$ is $C^2$ in the variables $(x,t)$, it follows that $\ell$ is also $C^2$ and thus from Proposition 5.1, the ODE in $t$ is verified for all $x \in [-V_{\i}T^{*},a]$, i.e
\begin{align*}
\tilde{\ell}_t(x,t,\rho)=&\int \left([\rho^*,\rho]+V_{\i} \right) \tilde{f}(x,t,\rho^*,\rho)\tilde{\ell}(x,t,\rho^*)\ \b(d\rho^*)  \\
&-\left( \int \left([\rho,\rho^+]+V_{\i}\right)\tilde{f}(x,t,\rho,\rho^+)\ \b(d\rho^+) \right)\tilde{\ell}(x,t,\rho)
\end{align*}
Now, it suffices to define
\[
\ell(x,t,\rho)=\tilde{\ell}(x-V_{\i}t,t,\rho) \text{ for all } (x,t) \in [0,a] \times [0,T^{*}] \text { and } \rho \in [P^-,P^+]^2
\]
then $\ell$ is $C^1$ in $(x,t)$ and verify the desired ODEs. Moreover, the total of mass of $\ell$ is conserved through space and time. 
\qed

\bs\noi
{\bf Declarations}
\begin{itemize}
\item The authors did not receive support from any organization for the submitted work.
\item The authors have no relevant financial or non-financial interests to disclose.
\item  The authors have no conflicts of interest to declare that are relevant to the content
of this article.
\item Data sharing not applicable to this article as no datasets were generated or analysed
during the current study
\end{itemize}

\bs\noi
{\bf Acknowledgments} The authors wish to express their sincere gratitude to Govind Menon for his very helpful comments on the first draft of this article.

\end{document}